\documentclass[11pt] {article} 
  \usepackage{amsmath}
    \usepackage{amssymb}
    
\usepackage{pdfsync}

\newtheorem{example}{Example}[section]
\newtheorem{theorem}{Theorem}[section]
\newtheorem{lemma}{Lemma}[section]
\newtheorem{corollary}{Corollary}[section]

\newtheorem{remark}{Remark}[section]

\newcommand{\eqnsection}{
   \renewcommand{\theequation}{\thesection.\arabic{equation}}
   \makeatletter
   \csname @addtoreset\endcsname{equation}{section} 
   \makeatother}
\def \ov{\overline}

\def \be{\begin{equation}}
\def \ee{\end{equation}}
\def \bt{\begin{theorem}} 
\def \et{\end{theorem}}
\def \bl{\begin{lemma}} 
\def \el{\end{lemma}}
\def \bea{\begin{eqnarray}}
\def \eea{\end{eqnarray}}
\def \bas{\begin{eqnarray*}}
\def \eas{\end{eqnarray*}}

\def \al{\alpha}
\def \bb{\beta}
\def \ga{\gamma}
\def\Ga{\Gamma}
\def \de{\delta}
\def \De{\Delta}
\def \ep{\epsilon}

\def \la{\lambda} 
\def \La{\Lambda}

\def \om{\omega}
\def \Om{\Omega}

\def \si{\sigma}

\def \th{\theta}

\def \ff{\infty}
\def \wh{\widehat}
\def \wt{\widetilde}

\def \cd{\,\cdot\,}

\def\stl{\stackrel{law}{=}}

\def \BB{{\cal B}}

\def \DD{{\cal D}}

\def \FF{{\cal F}}

\def \KK{{\cal K}}

\def \QQ{{\cal Q}}

\def \TT{{\cal T}}
\def \UU{{\cal U}}
\def \VV{{\cal V}}
\def \WW{{\cal W}}
\def \XX{{\cal X}}
\def \YY{{\cal Y}}
\def \ZZ{{\cal Z}}
\def \mN{{\overline {\mathbb N}}}
\def \mb { \mathbf  }

\def\b1{\mathbf 1}

\def \({\left(}
\def \){\right)}

\def \nn{\nonumber}
\def \Proof{\noindent{\bf Proof $\,$ }}

\def \bc{\begin{center} }
\def \ec{\end{center} }
\def \bs{\begin{slide} }
\def \es{\end{slide} }

\def\square{{\vcenter{\vbox{\hrule height.3pt
        \hbox{\vrule width.3pt height5pt \kern5pt
           \vrule width.3pt}
        \hrule height.3pt}}}}
\def\qed{{\hfill $\square$ \bigskip}}

\eqnsection

\begin{document}

\title{  Asymptotic properties of permanental   sequences}

  \author{  Michael B. Marcus\,\, \,\, Jay Rosen \thanks{Research of     Jay Rosen was partially supported by  grants from the Simons Foundation.   }}
\maketitle

 \begin{abstract}     
  Let $U=\{U_{j,k},j,k\in\mN\}$ be the potential  of  a transient symmetric Borel right process $X$ with state space $\mN$. For any excessive function   $f=\{f_{k,k\in \mN}\}$ for $X$ , $\wt  U=\{\wt U_{j,k},j,k\in\mN\}$, where
\begin{equation}
   \wt U_{j,k}= U_{j,k} +f_{ k},\qquad j,k\in\mN,\label{a.1}
   \end{equation}
is the kernel of an $\al$-permanental sequence  $\wt X_{\al}=(\wt X_{\al, 1} ,\ldots)$ for all $\al>0$.  The symmetric potential $U$ is also the covariance of a mean zero Gaussian sequence 
 $\eta=\{\eta_{j},j\in \mN\}$. 
 Conditions are given on the potentials $U$   and excessive functions $f$ under which,  
 \be
\limsup_{j\to \ff}\frac{  \eta_{j}}{( 2\,\phi_{j})^{1/2} }=1 \quad 
a.s. \quad \implies \quad \limsup_{n\to \ff}\frac{\wt  X_{\al, j}}{\phi _{j} }=1\quad 
a.s.,\label{a.2}
\ee
  for all $\al>0$, and sequences $\phi=\{\phi_{j}\}$ such that   $f_{j}=o(\phi _{j})$.
  
\medskip	  The function $\phi$ is determined by $U$. Many examples are given in which $U$ is the potential of symmetric birth and death processes with and without emigration, first and higher order Gaussian autoregressive sequences and L\'evy processes on $\mathbf Z$.

  \end{abstract}  
 
\maketitle

\section{Introduction}\label{sec-1}  

 We   define an $R^{n}$ valued   $\al$-permanental random variable  $ (\wt X_{\al,1 } ,\ldots,\wt X_{\al,n})$ to be  a non-negative random variable with Laplace transform,  
\begin{equation}
   E\(e^{-\sum_{i=1}^{n}s_{i}\wt X_{\al,i} }\) 
 = \frac{1}{ |I+K S|^{ \al}},   \label{int.1} 
 \end{equation}
for some  $n\times n$ matrix $K$  and diagonal matrix $S$ with positive entries $ s_{1},\ldots,s_{n} $, and $\al>0$.     We refer to the matrix $K$ as the kernel of $(\wt X_{\al ,1},\ldots,\wt X_{\al,n})$.     

An $\al$-permanental process $\wt X_{\al}=\{\wt X_{\al}(t),t\in \TT\}$   is a stochastic process that has finite dimensional distributions that are $\al$-permanental random variables. In this paper we   take $\TT=\mN$, the strictly positive integers, and refer to $\wt X_{\al}=\{\wt X_{\al,j},j\in\mN\}$
as an infinite dimensional  $\al$-permanental sequence. 

Eisenbaum and Kaspi, \cite[Theorem 3.1]{EK} show that the right hand side of (\ref{int.1}) is the Laplace transform of a  non-negative $n$-dimensional random variable for all $\al>0$ if and only if   $gKg=\{ g_{i}K_{i,j}g_{j}, i,j\in [1,n]\}$
 is the potential density of a transient Markov chain with state space $[1,n]$, for some strictly positive sequence $ \{g_{i}\}_{i=1}^{n} $.   In this paper we  combine the  $ \{g_{i}\}_{i=1}^{n} $ with $K$ and consider $ \wt  U =gKg$, which is the potential density of a transient Markov chain.

  The matrix   $\wt U$ is not necessarily  symmetric. When it is, it is the covariance of a Gaussian process. 
 Let $\eta=\{\eta_{j},j\in 1,\ldots, n\}$ be   a mean zero Gaussian vector with covariance   $ C=\{   C_{j,k},,j,k\in [1,n]\}$. It is well known that
\begin{equation}
   E\(\exp\(-\sum_{j=1}^{n}s_{j}\frac{\eta_{j}^{2}}{2} \) \)
 = \frac{1}{ |I+ CS|^{1/2}}. \label{int.1mm} 
 \end{equation}

The challenge is to find examples of $\wt U$ that are not symmetric. In this case the corresponding permanental processes are really something new. We obtain examples of   kernels $\wt U$ that are not symmetric by modifying symmetric kernels. 
	Let $X$ be a symmetric transient Markov process with potential density $U$ with respect to  counting measure and let $f=(f_{1},\ldots )$ be an excessive function for   $ X $. We  consider kernels $\wt U$ of the form,
 \be
 \wt  U_{j,k}=  U_{j,k} +f_{ k},\qquad j,k\in\mN.\label{1.10}
   \end{equation}

 Clearly, $\wt U$ is not symmetric. However,  	the kernels   of  $\al$-permanental random variables are not unique. For example, if $K$ satisfies (\ref{int.1}) so does 
 $\La K\La^{-1}$ for any $\La\in \DD_{n,+}$, the set of  $n\times n$  diagonal matrices with strictly positive diagonal entries.  We say that an $n\times n$  matrix $K$ is equivalent to a symmetric matrix, or symmetrizable,  if there exists an $n\times n$  symmetric matrix $W$ such that, 
\begin{equation}
    |I+KS| = |I+WS| \quad\mbox{for all $S\in \DD_{n,+}$}\label{1.3qq}.
   \end{equation}
 Nevertheless, it follows from  \cite[Theorem 1.1]{MRns} that in   Theorem \ref{theo-mchain}  below   we can always find excessive functions $f$ such that      $\{\wt U_{j,k}; l\le j,k\le n\}$ is not symmetrizable for all sufficiently large $l$ and $n$. In fact we show in \cite{MRns} that it is only in highly structured situations that the kernel of a permanental process is symmetrizable.

\medskip	 The fact that $\wt  U$ is the kernel of $\al$-permanental processes is given by the next theorem,  which is   part of  \cite[Theorem 1.11]{MRejp}.   
\bt\label{theo-borelN}
Let   ${  
  X}\!=\!
(\Om,  \FF_{t}, X_t,\th_{t},P^x
)$ be a  symmetric   transient Borel right process with state space $ \mN$, and  strictly positive potential  density  $ U$. Then for any   finite excessive function   $f$ for  $  X$ and $\al>0$, $\wt  U$  
is the kernel of an $\al$-permanental sequence  $\wt X_{\al} $.
\et

Recall that a non-negative function $f$ is excessive for   $   X$,  if $ P_{t}f(x)\uparrow  f(x)$ as $t\to 0$, for all $x$.
 The function $f$ is a   potential function of  $  X$, or simply a potential of $   X$, if $f= U h$ for some $h\ge 0$.  Since $ U h(x)=\int_{0}^{\ff}P_{t}h(x)\,dt$, it is easy to check that all potential functions are excessive. The potentials that play  a major role is this paper 
are $f= U h$ where $h\in\ell_{1}^{+}$  or $c^{+}_{0}$. 
Note that since $ U_{j,k}\leq  U_{j,j}\wedge  U_{k,k}$, (see \cite[(13.2)]{book}),  when $h\in\ell_{1}^{+}$,   $f_{j}=( U h)_{j}<\ff$  for all $j\in\mN$.

  In Theorem \ref{theo-borelN} we consider  two families of $\al$-permanental sequences; $\wt  X_{\al}$ with kernels $\wt  U$ and $   X_{\al}$ with kernels $ U$. Furthermore, $   X_{1/2}$ is a sequence of Gaussian squares as defined in (\ref{int.1mm}), (for all $n$). The primary goal of this paper is to find sharp results about the asymptotic behavior of $\wt X_{\al}=\{\wt X_{\al,j},j\in\mN\}$ as $j\to\ff$. The way we proceed is find finite excessive functions   $f$ for  $  X$  for which the asymptotic behavior of $\wt X_{\al}$ is the same as the asymptotic behavior of $  X_{1/2}$. Obtaining the asymptotic behavior of $  X_{1/2}$ is relatively simple because we are just dealing with Gaussian sequences.  To be more explicit, we  find finite excessive functions   $f$ such that
 \be
\limsup_{j\to \ff}\frac{\wt  X  _{\al, j}}{\phi _{j} }=\limsup_{j\to \ff}\frac{   X  _{1/2, j}}{\phi _{j} }\qquad 
a.s.\label{gena.kas}. 
\ee
The specific sequence of positive numbers $\phi=\{\phi_{j}\}$ is generally easily determined because $  X_{1/2}$ is a sequence of Gaussian squares.

\medskip	 
We get two classes of results. 
The first are general  limit theorems for permanental processes that  hold when  their kernels  $U$  and $\wt U$ satisfy certain general conditions. These are Theorems \ref{theo-mchain}--\ref{theo-1.10} given in  in Section \ref{sec-genres}. In Section \ref{sec-applic},  in    Theorems \ref{theo-3.1}--\ref{theo-lev} we apply these results to  the potential densities of specific families of Markov chains.  We consider birth and death processes, with and without emigration, and potentials that are the covariances of first and higher order autoregressive Gaussian sequences.

\subsection{General results}\label{sec-genres}

\medskip  	For any matrix $K$ let   $ K (l,n)$ denote the $n\times n$ matrix obtained by restricting the matrix  $ K $ to
 $  \{ l+1,\ldots,l+n\}\times \{ l+1,\ldots,l+n\}$. In the next theorem we consider $ U (l,n)^{-1}$. The reader should note that $  (U (l,n) )^{-1}$ is not  generally the same as the matrix $ U ^{-1}(l,n)$.
 
 	\medskip	For any invertible matrix $M$ we often denote $M^{-1}_{j,k}$ by $M^{j,k}$.

 \bt\label{theo-mchain}   Let $ X  $, $ U $,  $f$ and   $\wt  X  _{\al}$  be as   in Theorem \ref{theo-borelN}   and let  $\eta $ be a   Gaussian sequence with covariance $ U $.  Then  
\begin{equation}
\sum_{  k=1}^{n}( U (l,n))^{j,k}f_{k+l}\ge 0,\qquad j=1,\ldots n .\end{equation} 
  Suppose, in addition that,
\begin{equation}
\sum_{ j,k=1}^{n}( U (l,n))^ {j,k}f_{k+l}=o_{l}\(1\), \mbox{ uniformly in }n,\label{condreal}\ee 
and   there exists a sequence $\phi=\{\phi_{j}\}$  such that,    
 \be
 \limsup_{j\to \ff}\frac{  \eta_{j}}{( 2\,\phi_{j})^{1/2} }=1 \qquad 
a.s.\label{121.4a},
\ee
and   
\be  f_{j}=o(\phi _{j})\label{1.5mm}.
\ee Then  
 \be
\limsup_{j\to \ff}\frac{\wt  X  _{\al, j}}{\phi _{j} }=1\qquad 
a.s.\label{gena.ka}
\ee
 for all $\al\geq 1/2$. (Also, trivially, the upper bound holds for all $\al>0$.)
 \et
 
 In most of our applications of   this theorem we use results in \cite[Section 7]{MRejp} to show that the lower bound in  (\ref{gena.ka}) actually holds for all $\al>0$.

  The primary ingredient in   Theorem \ref{theo-mchain}  is the symmetric   potential density $    U=\{   U(j,k),j,k\in\mN\}$. We see in (\ref{condreal}) that   $( U(l,n))^{-1}$   must exist  for all $l$ and $n$. It follows from \cite[Theorem 13.1.2]{book} that this is the case.

  Theorem \ref{theo-mchain} is proved in Section \ref{sec-thm1.2}.

  \medskip	The next theorem  gives limit theorems for permanental sequences $\wt X_{\al}$ when the row sums of $ U$ in (\ref{1.10})
are uniformly bounded.    It has a  simpler more direct proof  than   Theorem \ref{theo-mchain} and doesn't require that we obtain the complicated estimate (\ref{condreal}).

    \bt\label{theo-1.8mm}      Let $ X, U, f$ and $\wt X_{\al}$  be as in Theorem \ref{theo-borelN}. 
If  
%(which  holds in particular when $f=  U h$ and $h\in \ell_{1}^{+}$),  **
    \begin{equation}
  \inf_{j} U_{j,j} >0, \quad \sup_{j} \sum_{k=1}^{\ff}   U_{j,k}<\ff,\quad\mbox{and}\quad f\in c^{+}_{0},\label{1.39}
   \end{equation}   
then
 \be
 \limsup_{n\to \ff}\frac{\wt X_{\al, n } }{   U_{n,n} \log n }=1\qquad 
a.s.\label{121.4mb}
\ee
\et

  Note that it follows from (\ref{1.39}) that $\sup_{n}  U_{n,n}<\ff$.
 
\medskip   	The proof of  Theorem \ref{theo-mchain} uses a result that compares the permanental sequence $\wt X_{\al}$ with the   Gaussian sequence $\eta $,  determined by the  covariance matrix $ U $. Therefore  $U$ must be symmetric.  The proof of    Theorem \ref{theo-1.8mm} does not involve Gaussian processes and so we don't need  $U$ to be  symmetric for that reason.    The requirement that $U$ must be symmetric is used because of Theorem \ref{theo-borelN}.  Theorem 6.1,  \cite{MRejp} is similar to Theorem \ref{theo-borelN}   but  does not require that $U$ is symmetric if $f$ is a left potential with respect to $U$, i.e., for all $k\in\mN,$  
\be
 f_{k}= \sum_{j=1}^{\ff}h_{j}  U_{j,k},\qquad\mbox{for some   $h\in\ell^{+}_{1}$ .}\label{1.38nnw}\ee
See \cite[(6.1)]{MRejp}.

\medskip	Using  \cite[Theorem 6.1]{MRejp} enables us to obtain limit theorems for permanental sequences with potentials of the form of (\ref{1.10}) in which $U$ is the potential of 
 Markov chains that  are not necessarily symmetric.

    \bt\label{theo-1.8gen}  Let   ${  
  X}\!=\!
(\Om,  \FF_{t}, X_t,\th_{t},P^x
)$ be a      transient Borel right process with state space $ \mN$ and  strictly positive potential  density  $ U$. 
Assume that  
%(which  holds in particular when $f=  U h$ and $h\in \ell_{1}^{+}$),  **
    \begin{equation}
  \inf_{j} U_{j,j} >0, \quad \sup_{j} \sum_{k=1}^{\ff}  U_{j,k}<\ff, \quad \mbox{and}\quad \sup_{k} \sum_{j=1}^{\ff}  U_{j,k}<\ff.\label{1.39gen}
   \end{equation}  
Let $f\in\ell_{\ff}$ be such that 
\be
 f_{k}= \sum_{j=1}^{\ff}h_{j}  U_{j,k},\qquad\mbox{for some $h\in \ell_{1}^{+}$},\label{1.38nn}\ee
  and let  $  \wt  U=\{\wt  U_{j,k}, j,k\in\mN\}$ where,
 \be
 \wt  U_{j,k}=  U_{j,k} +f_{ k},\qquad j,k\in\mN.\label{1.10ffs}
   \end{equation}
Then   for any   $\al>0$, $\wt  U$  
is the kernel of an $\al$-permanental sequence  $\wt X_{\al} $ and 
 \be
 \limsup_{n\to \ff}\frac{\wt X_{\al, n } }{   U_{n,n} \log n }=1\qquad 
a.s.\label{121.4mbgen}
\ee
\et

  Note that (\ref{1.10ffs})  looks the same as (\ref{1.10}) but here $U$ is not necessarily symmetric.
  Consequently, (\ref{121.4mbgen}) is of interest even for $f=0$.  (See Example \ref{ex-8.1}.)
  
\medskip	  Theorems \ref{theo-1.8mm} and \ref{theo-1.8gen} are  proved in Section \ref{sec-8}.

 \medskip	
   \medskip Let  $M$  be an $ \mN  \times  \mN  $ matrix and consider the operator norm on $\ell_{\ff}\to\ell_{\ff}$,  
  \begin{equation}
\|M\|=\sup_{\|x\|_{\ff}\leq 1}  \|Mx\|_{\ff}=\sup_{j} \sum_{k} |M_{j,k}|.\label{unif.1x1}
  \end{equation}
   We  say   that a Markov chain $X$  is uniform when  its $Q$ matrix has the property that $\|Q\|<\ff$. 
    Since all the row sums of $Q$ are negative,  
      \begin{equation}
\sup_{j}  |Q_{j,j}|\le \|Q\|\leq 2\sup_{j}  |Q_{j,j}|.\label{unif.4x1}
  \end{equation}
  (For information on uniform Markov chains, see \cite[Chapter 5]{F}.)

 \medskip	 The next theorem allows us to replace the hypotheses of Theorem  \ref{theo-1.8mm} with conditions on the $Q$ matrix of $X.$    Note that we call   $Q$  a   $(2k+1)-$diagonal matrix if $Q_{i,j}=0$ for all $|j-i|>k$.
 
\begin{theorem} \label{theo-1.10} Let $X$, $U$,  $f$ and   $\wt X_{\al}$  be as defined in Theorem \ref{theo-borelN} and assume furthermore that $X$ is a uniform Markov chain.  Then, if the row sums of the  Q-matrix of $X$ are bounded away from $0$, and $f\in c^{+}_{0}$,    
 \be
 \limsup_{n\to \ff}\frac{\wt X_{\al, n } }{  U_{n,n} \log n }=1\qquad 
a.s.\label{121.4mbj}
\ee

   Furthermore, when  the Q-matrix is   a $(2k+1)-$diagonal   matrix for some $k\geq 1$, $f\in c^{+}_{0}$ and   $f= U  h$ for $   h  \in c_{0}^{+}$ are equivalent.  
 \end{theorem}
 
 Theorem \ref{theo-1.10} is proved in Section \ref{sec-uniform}.

 \subsection{Applications}\label{sec-applic}

\medskip 
  
   The remaining  theorems in this section, Theorems \ref{theo-3.1}--\ref{theo-lev},
are applications of  the basic   Theorems \ref{theo-mchain}--\ref{theo-1.10}.  The basic theorems give general results for the quadruple $( X,\wt  X_{\al}, U,\wt U)$.   Our applications are examples based on specific choices of $ U$. We use different symbols for  the quadruple $( X,\wt  X_{\al}, U,\wt U)$ in the different examples. 

\medskip	 The simplest examples  of   symmetric transient Markov chains  are birth and death processes without emigration or explosion. We describe them by their   $Q$ matrix.

\medskip	  Let   $\mathbf s=\{s_{j},\,j\geq 1\} $ be a    strictly  increasing sequence  with   $s_{1}>0$ and   $\lim_{j\to\ff}s_{j}=\ff$, and let $ Y=\{ Y_{t},t\in R^{+}\}$ be a  continuous time birth and death process   on   $\mN  $ with $Q$ matrix    $  Q(\mathbf s)$ where,
\be -  Q({\mathbf s})=  \left (
\begin{array}{ cccccc cc}  
{a_{1}}+{a_{2}}&-{a_{2}}&0&\dots &0 &0 &\dots\\
-{a_{2}}&{a_{2}}+{a_{3}}&-{a_{3}}&\dots &0&0 &\dots\\
\vdots&\vdots&\vdots&\ddots&\vdots &\vdots &\ddots\\
0&0&0&\dots &{a_{j-1}}+{a_{j}}&-{a_{j}  }&\dots  \\
0&0&0&\dots &-{a_{j}}&{a_{j}}+{a_{j+1}}&\dots \\
\vdots&\vdots&\vdots&\ddots&\vdots &\vdots &\ddots\end{array}\right )\label{harr2.211},
  \ee
and  
  \begin{equation}
a_{1}=\frac{1}{s_{1}},\quad\mbox{and}\quad a_{j}=\frac{1}{s_{j}-s_{j-1}},\quad j\ge 1.\label{3.8bbx}
\end{equation}

Since all the row sums are equal to 0, except for the first row sum, we see that $Y$ is  a birth  and death process without emigration. (Except at the first stage. However, the first row can not also have a  zero sum because if it did,   $ Q({\mathbf s})([0,n])$ would not be invertible for any $n$.)

Since 
\begin{equation}
s_{j}=\sum_{k=1}^{j}  {1 \over a_{k}}, 
\end{equation}
the class of $Q$ matrices in (\ref{harr2.211}) include all symmetric birth and death processes for which  
\begin{equation}
\sum_{k=1}^{\ff}  \frac{1 }{a_{k}}=\ff.\label{nonexp.0x}
\end{equation}
This implies that $ Y $ does not explode, that is, it  does not run through all  $\mN  $ in finite time. See \cite[Theorem 5.1]{TK}.

We show in Theorem \ref{theo-3.1mm} that $ Y$ has potential densities $V$=\{ V$_{j,k},j,k\in \mN\}$ where,
\be
V_{j,k}= s_{j}\wedge s_{k}.\label{ypot.1q}
\ee
  
  The next theorem is an application of Theorem \ref{theo-mchain} to the quadruple $(Y,\wt Y_{\al},V,\wt V)$. This is an example of  $(X,\wt X_{\al},U,\wt U)$     in which $U=V$, in     (\ref{ypot.1q}).

\medskip	Let $s_{i}\uparrow \ff$ and define
\begin{equation}
\KK_{\mathbf s}(j) =\log \(\sum_{i=1}^{j-1}      1\wedge  \log (s_{i+1}/s_{i}) \).\label{ff}
\end{equation}
This function is introduced in \cite{Kovalb} to obtain limit theorems for certain Gaussian    sequences and is critical in our applications of Theorem \ref{theo-mchain}.

 	\begin{theorem}\label{theo-3.1}  Let $V$ be as given in (\ref{ypot.1q}).
  Let $f= Vh$, where $h\in\ell_{1}^{+}$, (which implies that $f_{j}=g(s_{j}),  j\ge 1$,
 where $g$ is an increasing  strictly concave function)  and let  $\wt Y_{\al}=\{\wt Y_{\al,j},j\in \mN\}$  be an $\al$-permanental sequence with kernel 
   $
 \wt V=\{ \wt V_{j,k};j,k\in\mN\}$, where
   \begin{equation}
 \wt   V _{j,k}=V_{j,k}+f_{k},\qquad j,k\in \mN. 
   \end{equation}
Then  
        \begin{equation}
 \limsup_{j\to \ff} {\wt Y_{\al,j} \over   s_{j}\KK_{\mathbf s}(j)}=1,  \qquad 
a.s.,\label{18.4mkw}\qquad \forall\, \al >0. \end{equation}
\end{theorem}

(We use the expression `$g$ is an increasing function' to include the case in which $g$ is non-decreasing.  We say that $g$ is a  strictly concave function when $\lim_{x\to\ff}g(x)/x=0.$ )

\medskip	Properties of   $   \KK_{\mathbf s}(j)$  are given in Lemma \ref{lem-2.7} and the examples   following it. Using them  we get  the following corollary of Theorem \ref{theo-3.1}.
  
   \begin{corollary}\label{cor-1.1} In Theorem \ref{theo-3.1}, 
   \begin{itemize}
 \item[(i)]
  if     $\limsup_{j\to\ff}s_{j}/s_{j-1}<\ff$, then
\begin{equation}
\limsup_{j\to \ff} { \,\wt Y_{\al,j} \over   s_{j}\log \log s_{j}}=1,  \qquad 
a.s.\label{1.24mm},   \qquad \forall\, \al>0.
\end{equation} 
 \item[(ii)] If  $\liminf_{j\to\ff}s_{j}/s_{j-1}>1$, then
\begin{equation}
\limsup_{j\to \ff} { \wt Y_{\al,j} \over     s_{j}\log {j }}=1,  \qquad 
a.s.\label{1.25mm}, \qquad \forall\, \al>0.
\end{equation} 
 \end{itemize}

 \end{corollary}

   We show in Section \ref{sec-BM} that     the potentials $f= Vh$, where $h\in\ell_{1}^{+}$, satisfy (\ref{condreal}) and (\ref{1.5mm}). This   allows us to apply 	Theorem \ref{theo-mchain}.  In Section \ref{sec-BM} we also give a Riesz representation theorem for functions  that are excessive for  $X$.

\medskip	 

In Section  \ref{sec-gen} we modify $Q(\mathbf s)$,   to obtain    $Q$ matrices for a large class of birth and death processes with emigration.
Let $\BB=diag(b_{1},b_{2},\ldots)$, i.e., $\BB$ is a diagonal matrix with diagonal elements $(b_{1},b_{2},\ldots  )$). Define  \begin{equation}
 -\wt Q(\mathbf s)=  -\wt Q(\mathbf s, \mb b)= \BB(-Q(\mathbf s))\BB.\label{4.1mm}
   \end{equation}
We show that when $b_{j}=g(s_{j})$, $j\ge 1$, where $g(x)$ is an  increasing strictly  concave function, $\wt Q(\mathbf s)$ is the $Q$ matrix of a birth and death  process with emigration.

\medskip	Let $V$ be as given in (\ref{ypot.1q}) and let $ W=\{ W_{j,k};j,k\in\mN\}$, where,
\begin{equation}
 W_{j,k}=b^{-1}_{j}V_{j,k}b^{-1}_{k},\qquad j,k\ge 1.\label{1.16nn}
   \end{equation}
We show in Lemma \ref{lem-specmar} that   $W$ is the potential    density  of a Markov chain  $Z$ with $Q$ matrix  $\wt Q(\mathbf s, \mb b)$.      
   
\medskip	 The next theorem generalizes Theorem \ref{theo-3.1} and  Corollary  \ref{cor-1.1} .

\begin{theorem}\label{theo-4.2mm}
  Let  $W$ be as given in (\ref{1.16nn}) and let $f= W  h$, where $   h  \in\ell_{1}^{+}$.   Let  $\wt Z_{\al}=\{\wt Z_{\al,j},j\in \mN\}$ be an $\al$-permanental sequence
  with   kernel
   $ \wt W=\{ \wt W_{j,k};j,k\in\mN\}$, where
   \begin{equation}
\wt W_{j,k}=W_{j,k}+f_{k},\qquad  j,k\in\mN. \label{4.38mm}
   \end{equation}
  Then 
\begin{equation}
 \limsup_{j\to \ff} { \,\wt Z_{\al,j} \over   W_{j,j}\KK_{\mathbf s}(j)}=1,  \qquad 
a.s.\label{18.4mkff},   \qquad \forall\, \al>0.
\end{equation} 
If  $\limsup_{j\to\ff}s_{j}/s_{j-1}<\ff$, then
\begin{equation}
\limsup_{j\to \ff} { \,\wt Z_{\al,j} \over   W_{j,j}\log \log s_{j}}=1,  \qquad 
a.s.\label{1.24mmx},   \qquad \forall\, \al>0.
\end{equation} 
If  $\liminf_{j\to\ff}s_{j}/s_{j-1}>1$, then
\begin{equation}
\limsup_{j\to \ff} { \wt Z_{\al,j} \over     W_{j,j}\log {j }}=1,  \qquad 
a.s.\label{1.25mmx}, \qquad \forall\, \al>0.
\end{equation}

\end{theorem}

In Lemma \ref{lem-4.4mm2} we show under the hypotheses of Theorem \ref{theo-4.2mm},   $\lim_{j\to\ff} f_{j}/W_{j,j} =0$.  

\medskip	
Clearly, when $\BB$ is the identity matrix, Theorem \ref{theo-4.2mm} gives Theorem \ref{theo-3.1}, It is useful to state Theorem \ref{theo-3.1} separately because it is instrumental in the proof of Theorem \ref{theo-4.2mm}.

\medskip	 An interesting class of examples of potentials of the form of (\ref{1.16nn})  is when  $b_{j}=s^{1/2}_{j}$ for all $j\in\mN$.
 We denote this potential   density  by $\WW=\{\WW_{j,k},j,k\in \mN\}$, and see that 
\begin{equation}
  \WW_{j,k} =\frac{b_{j}}{b_{k}},\qquad j\le k\label{3.28ja}.
   \end{equation}
This expression is much more interesting if we  set $b_{j}=e^{v_{j}}$. Then we get 
 \be
 \WW_{j,k}=e^{-|v_{k}-v_{j}|},\qquad \forall\, j,k\in\mN.\label{5.6mmq}
 \ee
 Let $\xi=\{\xi(x),x\in R^{+}\}$ be a mean zero Gaussian process with covariance $\exp(-|x-y|)$ and note that $\xi$ is an Ornstein-Uhlenbeck process and  $\WW$ is the covariance of 
 $ \{\xi( v_{j}),j\in\mN\}$.
  In {Theorem} \ref{theo-4.2mmfq} we show what results  Theorem \ref{theo-4.2mm} gives when $W$ is written as in (\ref{5.6mmq}).

\medskip	In Section  \ref{sec-6} we take the  potential   densities, described abstractly in (\ref{1.16nn}), to be the covariance of a first order auto regressive Gaussian sequence.
Let $\{g_{j}, j\in \mN\}$ be    a sequence of  independent identically distributed standard normal random variables and $\{x_{n}\}$ 
an increasing  sequence    with $0< x_{n}\le 1$.   We consider first order autoregressive Gaussian sequences $\wh \xi=\{\wh \xi_{n},n\in\mN\}$,  defined by,
 \begin{equation}
\wh\xi_{1}=g_{1},\qquad\wh\xi_{n}=x_{n-1}\wh\xi_{n-1}+g_{n},\qquad n\ge 2.\label{auto.1oo}
 \end{equation}
The covariance of  $\wh\xi$ is  $\UU=\{\UU_{j,k},j,k\in \mN\}$, where,   \begin{equation}
     \UU_{j,k}=\sum_{i=1}^{j}\(\prod_{l=i}^{j-1}x_{l}\prod_{l=i}^{k-1}x_{l}\),\qquad j\le k.\label{comp.10kk}
      \end{equation}  
and $\{x_{j}\}$ is an increasing sequence, with  $0<x_{j}\le 1$. 
This has the form of (\ref{1.16nn}) with  
 \begin{equation}
   b_{j}=\prod_{l=1} ^{j-1}x_{l}^{-1},\qquad\mbox{and}\qquad s_{j}=\UU_{j,j}b_{j}^{2}=\sum_{i=1}^{j}b_{j}^{2},
   \end{equation}
  and consequently, as we  show,  is the potential  density  of a Markov chain    which we denote by   $\XX$.  
  In addition we show in Lemma \ref{lem-5.2mm} that  $\lim_{j\to\ff}\UU_{j,j}$ exists  and is strictly greater than 1.  
 	In this case Theorem \ref{theo-4.2mm} gives:

 		\begin{theorem}
 \label{theo-4.2mmf}   
  Let  $  \UU$ be as given in (\ref{comp.10kk}) and let $f$ be a finite  excessive function for $\XX$.  Let  $\wt \XX_{\al}=\{\wt \XX_{\al,j},j\in \mN\}$ be an $\al$-permanental sequence
 with kernel        
      $
  \wt    \UU=\{ \wt    \UU_{j,k};j,k\in\mN\}$, where
   \begin{equation}
 \wt    \UU_{j,k}=\UU_{j,k}+f_{k},\qquad  j,k\in\mN. \label{4.38nnq}
   \end{equation}

 \begin{itemize}

 \item [(i)]    If  $\lim_{j\to\ff}\UU_{j,j}=\ff$, or equivalently,  $\lim_{j\to \ff}x_{j}= 1 $, and $f= \UU  h$, where $   h  \in\ell_{1}^{+}$,  then
 \begin{equation}
   \limsup_{j\to\ff}\frac{\wt \XX_{\al,j}}{   \UU_{j,j} \log \log (\UU_{j,j}b_{j}^{ 2}) }= 1\quad a.s.\label{5.44mm}
   \end{equation}
   
  In particular, if    
$\,\UU_{j,j}$ is a regularly varying function with index $0<\bb<1$,   then
 \begin{equation}
   \limsup_{j\to\ff}\frac{\wt \XX_{\al,j}}{   \UU_{j,j}  \log j }= 1-\bb\quad a.s.\label{5.44vv}
   \end{equation}

 \item[(ii)] If  $\lim_{j\to\ff}\UU_{j,j}=1/(1-\de^{2})$, for some $0<\de<1$, or equivalently,  $\lim_{j\to \ff}x_{j}=  \de<1 $, and $f\in c^{+}_{0}$,     then
  \begin{equation}
   \limsup_{j\to\ff}\frac{\wt \XX_{\al,j}}{   \log j }= \frac{1}{1-\de^{2}}\quad a.s.\label{5.45mm}
   \end{equation}
    \end{itemize}
 Furthermore, when  $\lim_{j\to\ff}\UU_{j,j}=1/(1-\de^{2})$, for some $0<\de<1$, $f\in c^{+}_{0}$ and   $f= \UU  h$ for $   h  \in c_{0}^{+}$ are equivalent.
\end{theorem}

    The statement in (\ref{5.44mm})   and even the one in (\ref{5.44vv}) do not   seem too useful because there are too  many unknowns. Ultimately everything depends on the sequence $\{x_{j}\}$. We give some examples.
 They are arranged in order of decreasing values of $ x_{j}$, (for large $j$).

  \begin{example} \label{ex-1.1}{\rm 

 \begin{itemize}

 \item[(i)] If  $j(1-x ^{2}_{j})\to
  0$ as $j\to\ff$,
 \begin{equation}
   \limsup_{j\to\ff}\frac{\wt \XX_{\al,j}}{   j \log \log j }= 1\qquad a.s.\label{5.44z}, \qquad \forall\, \al>0.
   \end{equation}
This includes the case where $\prod_{j=0}^{\ff}x _{j} >0 $.
 \item[(ii)] If  $ j(1-x^{2} _{j})\sim c $ as $j\to\ff$   for some $c>0$,
  \begin{equation}
   \limsup_{j\to\ff}\frac{\wt \XX_{\al,j}}{   j \log \log j }= \frac{1}{1+c}, \qquad a.s.\label{5.44za}, \qquad \forall\, \al>0.
   \end{equation}
   
  \item[(iii)] If  $ j^{\bb}(1-x^{2} _{j})\sim 1 $ as $j\to\ff$,       for   $0<\bb<1$,
 \begin{equation}
   \limsup_{j\to\ff}\frac{\wt \XX_{\al,j}}{  j^{\bb}   \log j }= 1-\bb,\quad a.s.\label{5.44mh}, \qquad \forall\, \al>0.
   \end{equation}
    \end{itemize}
%    \item[(i)]
% If    $\lim_{j\to \ff}x_{j}=  \de<1 $,
%\begin{equation}
% \limsup_{j\to \ff} {  \wh Y_{\al,j} \over    \log j}={1 \over   1-\de^{2} },  \qquad 
%a.s., \qquad \forall\, \al>0.
%\label{18.4z}
%\end{equation} 
 }\end{example}
   
   In Section \ref{sec-7} we take the symmetric potential    densities $U$ in (\ref{1.10}) to be the covariance of a   $k$-th  order autoregressive Gaussian sequence,   $k\ge 2$.  Let $\{g_{j}, j\in \mN\}$ be  a sequence of  independent identically distributed standard normal random variables  and $\{p_{i}\}_{i=1} ^{k} $  a decreasing sequence of probabilities  with  $\sum_{l=1}^{k}p_{l}\le 1$.   We define the Gaussian sequence  $\xi=\{\xi_{n},n\in\mN \}$  by,
\begin{equation}
\xi_{1}=g_{1}, \hspace{.2 in}\mbox{and}\quad\xi_{n}=\sum_{l=1}^{k}p_{l}\xi_{n-l}+g_{n},\hspace{.2 in}n\geq 2,\label{ark.1o}
\end{equation}
where $\xi_{i}=0$ for  all $i\le 0$.   Let   
\be
\VV=\{\VV_{m,n}, m,n\in \mN\},\label{W}
\ee
be the covariance of $\xi$.

 We show that with certain  additional conditions, $\VV$ is the potential  density of a continuous Markov chain $\YY$ on $\mN$ with a  $Q$ matrix that  is a symmetric T\"oeplitz matrix  which is completely determined  by $\{Q_{n,m}\}_{m\ge n}$,  i.e.,  
\bea
  && Q_{n,n}=-\(1+\sum_{l=1}^{k}p_{l}^{2}\),\quad  Q_{n,j}=\bb_{j}>0,\quad j\in[n+1,n+k],\nn\\
  &&Q_{n,j}=0,\quad j> n+k,\label{1.49mm}
   \eea
   where $\bb_{j}$ are functions of $\{p_{i}\}_{i=1} ^{k} $.
   In addition, the row sums of the $n$-th row of $Q$, for $n\ge k+1$, is equal to $(1-\sum_{l=1}^{k}p_{l})^{2}$. 
   
\medskip	 We can consider these Markov chains as population models  which, when at stage $n\ge k+1$,  increase or decrease by 1 to $k$ members,  and so are generalizations of birth and death processes.  When $\sum_{l=1}^{k}p_{l}=1$, there is no emigration once the population size reaches $k$. When $\sum_{l=1}^{k}p_{l}<1$, there is   emigration at each stage.

\begin{theorem}\label{theo-4.2mmfx}  Let  $\VV$ be as defined in (\ref{W}) with the additional property that $p_{i}\downarrow$, and let $f$ be a finite  excessive function for $\YY$.  Let  $\wt \YY_{\al}=\{\wt \YY_{\al,j} ,j\in \mN\}$ be an $\al$-permanental sequence
 with kernel  
    $
\wt \VV=\{\wt \VV_{j,k};j,k\in\mN\}$ where,
   \begin{equation}
\wt \VV_{j,k}=\VV_{j,k}+f_{k}. \label{7.85}
   \end{equation}
 
  \begin{itemize}
 
\item[(i )] If    $\,\sum_{j=1}^{k}p_{j}<1$ and   $f\in c^{+}_{0}$,     then  
\begin{equation}
 \limsup_{j\to \ff} { \wt \YY_{\al,j}\over    \log j}=c^{*},\qquad\forall\,\al>0, \qquad 
a.s.\label{7.87}
\end{equation} 
 for some constant
\begin{equation}
   1+p_{1}^{2}\le c^{*}\le   \frac{1}{1-(\sum_{l=1}^{k}p_{l})^{2}}.\label{1.55mm}
   \end{equation} 
  The precise value of $c^{*}$ is given in (\ref{7.4nn}). 
 \item[(ii)]  
 If   $\,\sum_{j=1}^{k}p_{j}=1$ and  in addition $f=\VV h$ where $h\in \ell_{1}^{+}$,  then
 \begin{equation}
   \limsup_{j\to\ff}\frac{\wt \YY_{\al,j}}{   j \log \log j }= \frac{1}{\(\sum_{l=1}^{k}lp_{l}\)^{2}}\qquad a.s.,\qquad\forall\,\al\ge 1/2.\label{7.86}
   \end{equation}

   \end{itemize} 
   Furthermore, when  $\,\sum_{j=1}^{k}p_{j}<1$, $f\in c^{+}_{0}$ and   $f= \VV  h$ for $   h  \in c_{0}^{+}$ are equivalent.
   \end{theorem} 
   
       The limits in (\ref{7.86}) and (\ref{7.87}) may also hold for certain sequences $\{p_{i} \}$ that are not decreasing. See Remark \ref{rem-dec}.

\medskip   We show in   Lemma \ref{lem-4.3ja} that  when $\sum_{j=1}^{k}p_{j}=1$, the condition   $f=\VV h$ where $h\in \ell_{1}^{+}$,  holds for all concave increasing functions $f$ satisfying $f _{j}=o(j)$ as $j\to\ff$.  Furthermore it is trivial that the upper bound in (\ref{7.86}) holds for all $\al>0$. But we need additional conditions on the potentials $f$ to show that the lower bound holds for all $\al>0$.

\medskip	
In the next theorem we show that    (\ref{7.86}) holds for all $\al>0$, when the potentials $f$ are such that $f_{j}=o(j^{1/2})$ as $j\to\ff$. We don't think that this restriction is required but we need it to use the techniques that we have at our disposal.

\begin{theorem} \label{theo-rest} Under the hypotheses of Theorem \ref{theo-4.2mmfx} assume in addition that $f_{j}=o(j^{1/2})$ as $j\to\ff$. Then (\ref{7.86}) holds for all $\al>0$.
 \end{theorem}

When $\sum_{j=1}^{k}p_{j}<1$, the condition   $f=\VV h$ where $h\in \ell_{1}^{+}$, implies that $f\in \ell_{1}$. 
    In Remark \ref{rem-6.2} we give an explicit formula for $c^{*}$ in terms of the roots of the polynomial
   \begin{equation}
   P(x)=1-\(\sum_{l=1}^{k}p_{l}x^{l}\).
   \end{equation}

 \medskip  We   use ideas from the proofs of   Theorems  \ref{theo-1.8mm} and  \ref{theo-1.8gen} to get limit theorems for permanental sequences with kernels that are related to the potentials of   L\'evy processes that are not necessarily symmetric.

 \begin{theorem} \label{theo-lev}   Let $X$ be a L\'evy process on $\mathbb Z $ that is  killed  at the end of an independent exponential time, with potential density   $U=\{U_{j,k};j,k\in \mathbb Z \}$. 
 
Let  $f=\{f_{k}, k\in \mathbb Z\}$ be a finite excessive function for $X$,    and let  $  \wt  U=\{\wt  U_{j,k}, j,k\in\mathbb Z\}$ where,
 \be
 \wt  U_{j,k}=  U_{j,k} +f_{ -k},\qquad j,k\in\mathbb Z.\label{lev.0}
   \end{equation}
Then   for any   $\al>0$, $\wt  U$  
is the kernel of an $\al$-permanental sequence  $\wt X_{\al} $, 
and  if  $\lim_{k\to\ff}f_{-k}=0$, then  
 \be
 \limsup_{n\to \ff}\frac{\wt   X_{\al, n } }{   \log n }=U_{0,0}.\qquad 
a.s.\label{lev.1},
\ee
%where $\{p_{t}(j,k);j,k\in\mN\}$ are the transition probabilities of $X$ and $p_{t}(0)=p(j,j)$ for all $j %\in\mN.$

     Furthermore, if   $g= U h$, for a positive sequence $   h   $, then  $  g  \in c_{0}^{+} (\mathbb Z)$ if and only if  $   h  \in c_{0}^{+} (\mathbb Z)$. 
 \end{theorem}

Note that  when $X$ is not symmetric, (\ref{lev.1}) is of interest even for $f=0$.

\medskip  Theorems \ref{theo-borelN}--\ref{theo-1.10}, are results for the broad classes of permanental processes described by quadruple $( X,\wt  X_{\al}, U,\wt U)$.  Theorem  \ref{theo-borelN} is   given in   \cite[Theorem 1.11]{MRejp}.   
 Theorem \ref{theo-mchain} is proved in Section \ref{sec-thm1.2}.  Theorem \ref{theo-1.8mm} and \ref{theo-1.8gen} are  proved in Section \ref{sec-8} and Theorem \ref{theo-1.10} is proved in Section \ref{sec-uniform}.  Theorems \ref{theo-3.1}--\ref{theo-lev}
are applications of  Theorems \ref{theo-borelN}--\ref{theo-1.10} in which the matrices  $U$ are  the potential densities of specific families of Markov chains.
 We use different symbols for  $( X,\wt  X_{\al}, U,\wt U)$ in the different examples.

  In Section  \ref{sec-BM} we take $U=V=\{   V_{j,k},j,k\in\mN\}$ where,
  \begin{equation}
   V_{j,k}=s_{j}\wedge s_{k},\qquad \mbox{for }  s_{j}\uparrow \ff,  \label{1.59mm}
   \end{equation}
   and give the proof of Theorem \ref{theo-3.1}.
 
   In Section \ref{sec-gen} we take $U=W=\{   W_{j,k},j,k\in\mN\}$ where,
 \be  
 W_{j,k}= \frac{s_{j}\wedge s_{k}}{b_{j}b_{k}}\label{1.60mm}
   \end{equation}
and $  b=\{b_{j}\}$ is a finite potential for  the Markov process determined by   $V$.  Theorem \ref{theo-4.2mm} is proved in this section. We consider the specific example given in (\ref{3.28ja}) and (\ref{5.6mmq}) in which,
$U=\WW=\{   \WW_{j,k},j,k\in\mN\}$ where,
 \be
\WW_{j,k}=e^{-|v_{k}-v_{j}|},\qquad \forall\, j,k\in\mN,\label{5.6mmqmm}
 \ee
For a sequence $v_{j}\uparrow\ff$.  Theorem \ref{theo-4.2mmfq} gives limit theorems for permanental processes based $\WW$.

\medskip	 
In Section  \ref{sec-6} we take   $U=\VV$ to be the  covariance of a first order autoregressive Gaussian sequence. In this case $\UU_{j,k}$ is also an example of (\ref{1.60mm}) in which  
 \begin{equation}
   b_{j}=\prod_{l=1} ^{j-1}x_{l}^{-1},\qquad\mbox{and}\qquad s_{j}=\sum_{i=1}^{j} {b^{2}_{i}}.
   \end{equation}
We  give the proof of Theorem \ref{theo-4.2mmf} in this section.

\medskip	     In Sections \ref{sec-BM}--\ref{sec-6} the potentials are all examples of (\ref{1.60mm}). The Markov chains with these potentials only move between their nearest neighbors.    In Section \ref{sec-7} we take the symmetric potential $U$ in (\ref{1.10}) to be the covariance of a   $k$-th  order autoregressive Gaussian sequences for    $k\ge 2$, and  denote it by $\VV$.      Markov chains   with these potentials move amongst  their $k$ nearest neighbors. We can not find the potentials of these chains precisely but we can estimate the potentials sufficiently well to give a proof of Theorems \ref{theo-4.2mmfx} and Theorem \ref{theo-rest}.

  We thank Pat Fitzsimmons and  Kevin O'Bryant for several helpful conversations.

\section{Birth and death processes}\label{sec-BM}

Let   $\mathbf s=\{s_{j},\,j\geq 1\} $ be a    strictly  increasing sequence  with $s_{j}>0$ and   $\lim_{j\to\ff}s_{j}=\ff$, and let $ \ov Y=\{ \ov Y_{t},t\in R^{+}\}$ be the  continuous time birth and death process   on $\mN$,  without emigration, with $Q$ matrix  $\ov Q(\mathbf s)$ where,
\be -\ov Q({\mathbf s})= \frac{1}{2}\left (
\begin{array}{ cccccc cc}  
{a_{1}}+{a_{2}}&-{a_{2}}&0&\dots &0 &0 &\dots\\
-{a_{2}}&{a_{2}}+{a_{3}}&-{a_{3}}&\dots &0&0 &\dots\\
\vdots&\vdots&\vdots&\ddots&\vdots &\vdots &\ddots\\
0&0&0&\dots &{a_{j-1}}+{a_{j}}&-{a_{j}  }&\dots  \\
0&0&0&\dots &-{a_{j}}&{a_{j}}+{a_{j+1}}&\dots \\
\vdots&\vdots&\vdots&\ddots&\vdots &\vdots &\ddots\end{array}\right )\label{harr2.211q},
  \ee
and  
  \begin{equation}
a_{1}=\frac{1}{s_{1}},\quad\mbox{and}\quad a_{j}=\frac{1}{s_{j}-s_{j-1}},\quad j\ge 1.\label{3.8bb}
\end{equation}
Since 
\begin{equation}
s_{j}=\sum_{k=1}^{j}  {1 \over a_{k}}, 
\end{equation}
the class of $Q$ matrices in (\ref{harr2.211q}) include all symmetric birth and death processes for which  
\begin{equation}
\sum_{k=1}^{\ff}  {1 / a_{k}}=\ff.\label{nonexp.0}
\end{equation}
This implies that $\{ \ov Y_{t},t\in R^{+}\}$ does not explode, that is it, does not run through all  $\mN$ in finite time; see \cite[Theorem 5.1]{TK}.  

\begin{theorem} \label{theo-3.1mm}
The continuous time birth and death  process $ \ov Y$ has potential densities,  
\be
\ov V_{j,k}=2 \(s_{j}\wedge s_{k}\),\quad j,k\in\mN.\label{ypot.1}
\ee
\end{theorem}

\Proof     It is easy to see that  
$\ov V \,\ov Q(\mathbf s)=\ov Q(\mathbf s)\ov V=-I $ in the sense of matrix multiplication.  However,       generally, this is not sufficient  to show that $ \ov Y$ has potential densities $\ov V_{j,k}$,  (unless $\sup_{j}  a_{j}<\ff$, see Lemma \ref{lem-invpot}). We see in Lemma \ref{lem-3.3ll} that there are functions $f$ with $\ov Q(\mathbf s)f=0$.  

\medskip		Let $\ov B=\{\ov B_{t},t\in R^{+}\}$ be Brownian motion killed the first time it hits $0$.  $\ov B$ has potential densities  
\begin{equation}
U_{\ov B}(x, y)=2 \(x\wedge y\),\qquad x,y>0.\label{3.23j}
\end{equation}
We use $\ov B$ to prove (\ref{ypot.1}). To do this we first make the connection between $ \ov Y$ and $\ov B$.

Using (\ref{harr2.211q}) and the relationship between the Q matrix and the   jump matrix   of the Markov chain,  (see  \cite[Section 2.6]{Norris}), we have that  for all $n\geq 2$,    
\bea
P_{ \ov Y}\(n,n+1\)&=&{a_{n+1} \over a_{n }+a_{n+1}}\label{harr2.a}\\
&=&{s_{n}-s_{n-1} \over s_{n+1}-s_{n-1}} =P^{s_{n}}_{\ov B}\(T_{s_{n+1}}<T_{s_{n-1}}\),\nn
\eea
where we use   \cite[Chapter II, Proposition 3.8]{RY} for the last equality. (As usual, $T_{x}$ is the first hitting time of $x$.) 

Similarly,  
\bea
P_{ \ov Y}\(n,n-1\)&=&{a_{n} \over a_{n }+a_{n+1}} \label{harr2.b}\\
&=&{s_{n+1}-s_{n} \over s_{n+1}-s_{n-1}}=P^{s_{n}}_{\ov B}\(T_{s_{n-1}}<T_{s_{n+1}}\).\nn
\eea
In the same manner we   have,  
\begin{equation}
P_{ \ov Y}\(1,2\)={a_{2} \over a_{1 }+a_{2}}=  P^{s_{ 1}}_{\ov B}\(T_{s_{ 2}}<T_{\De}\),\label{harr2.a1}
\end{equation}
where $\De$ is the cemetery state,  and  

\begin{equation}
P_{ \ov Y}\(1,\De\)={a_{1} \over a_{1 }+a_{2}}=P^{s_{ 1}}_{\ov B}\(T_{\De}<T_{s_{ 2}}\).\label{harr2.b1}
\end{equation}

Now, let  $L_{t}^{x}$ denote the  local time of Brownian motion.  It follows from  \cite[Chapter VI, (2.8)]{RY}, that for all $n\geq 1$, 
\bea 
E_{\ov B}^{s_{n}} \(L^{s_{n}}_{T_{s_{n-1}}\wedge T_{s_{n+1}}}\)&=&2{\(s_{n+1}-s_{n}\)\(s_{n}-s_{n-1}\) \over s_{n+1}-s_{n-1}}\label{harr2.c} \\
&=&{2 \over a_{n}+a_{n+1}}.\nn
\eea 
We see from \cite[Section 2.6]{Norris} and the $Q$ matrix  in  (\ref{harr2.211q}) that
the holding time of  $ \ov Y$ at $n$ is an exponential  random variable with parameter $(a_{n}+a_{n+1})/2 $ that is  independent of everything else.  	  This holding time has expectation ${2 / (a_{n}+a_{n+1}})$.

\medskip To obtain  (\ref{ypot.1}) we  show that the behavior of $ \ov Y$  and  $\ov B$ are similar
in the following sense: Begin  $ \ov Y$ at $j$ and $\ov B$ at $s_{j}$. The  next visit of $ \ov Y$ to an integer will be  to either $j+1$, with probability    (\ref{harr2.a}),    or to  $j-1$ with  probability  (\ref{harr2.b}). These are the same probabilities that the    next visit of $ \ov B$ is to $s_{j+1}$    or   $s_{j- 1}$.    During the time interval that $ \ov Y$   and $\ov B$ make this transition,   it follows from the last paragraph that the expected value of the increase  in $L_{t}^{s_{j}}$   is  the expected amount of time that $ \ov Y$
spends at $j$.   We repeat this analysis  until the processes  move to $\De$,   at which time they die. It follows from this   that,
\begin{equation}
\ov V_{j,k}=E^{j}\(\int_{0}^{\ff} 1_{k}\( \ov Y_{t}\)  \,dt \)=E_{\ov B}^{s_{j}} \(  L^{s_{k}}_{\ff}\)=U_{\ov B}(s_{j},s_{k}), \label{harr2.d}
\end{equation}
which, by (\ref{3.23j}), gives (\ref{ypot.1}).\qed

To simplify the    notation we consider the continuous time Markov chain   
\be
Y=\{Y_{t},t\in R^{+}\}=\{\ov Y_{2t},t\in R^{+}\},\label{3.12nn}
\ee which has potential densities   given by the matrix $V=\{V_{j,k};j,k\in\mN\}$ with, 
\begin{equation}
V_{j,k}=  s_j\wedge s _{k}  = E\(B_{s_{j}}B_{s_{k}}\),\label{3.28kja}
\end{equation} 
where $\{B_{t},t\in R^{+}\}$ is standard Brownian motion, and    $Q$ matrix, 
\begin{equation}
Q(\mathbf s)=2\ov Q(\mathbf s).\label{3.14mm}
\end{equation}

\medskip	One of our goals is to study permanental processes with kernels of the form (\ref{1.10}). To that end we now 
describe the   finite potentials and excessive functions  of $Y$.

\begin{theorem}\label{theo-3.1f} A potential $ f=Vh$ is finite if and only if $ h\in \ell_{1}^{+}.$
When this is the case the following  equivalent conditions hold:   
\begin{itemize}
     \item[(i)]  
\begin{equation}
\frac{f_{n}-f_{n-1}}{s_{n}-s_{n-1}}\downarrow 0\label{3.38mp} ,
\ee
where we take $f_{0}=s_{0}=0$.
   \item[(ii)]  the function   $g(s_{n})=f_{n} $ is   concave on  $\{   0\}\cup \{s_{j},\,j\geq 1\} $  and  
\be   \quad \frac{f_{n}}{s_{n}}\downarrow 0\label{3.38mm} .
\end{equation}
 
    \end{itemize}
\et

\Proof  We point out
in the second paragraph following Theorem \ref{theo-borelN} that   $ f=Vh$ is finite   when  $ h\in \ell_{1}^{+}.$   The reverse implication follows from the fact that \begin{equation}
f_{1}=\sum_{k=1}^{\ff} V_{ 1,k}h_{k}=s_{1 }\sum_{k=1}^{\ff}h_{k},\label{2.18mm}
\end{equation} 
where we use   (\ref{3.28kja}).  

In general we have 
\be 
f_{n}  = \sum_{k=1}^{n }s_{k}h_{k}+s_{n }\sum_{k=n+1}^{\ff}h_{k}, \label{3.35mm}
\ee
and
\be
f_{n+1}  =  \sum_{k=1}^{n }s_{k}h_{k}+s_{n+1}\sum_{k=n+1}^{\ff}h_{k}.\label{3.35mmj}
\ee
Therefore,   
\begin{equation}
\frac{f_{n+1}-f_{n}}{s_{n+1}-s_{n}}= \sum_{k=n+1}^{\ff}h_{k}.\label{3.37mm}
\end{equation}
This   and (\ref{2.18mm}) gives $(i)$. It also shows that     $g(s_{n})=f_{n} $ is  a concave function on    $\{   0\}\cup\{s_{j},\,j\geq 1\} $.

Note that if we divide  (\ref{3.35mm}) by $s_{n}$ and (\ref{3.35mmj}) by $s_{n+1}$,  and use the fact that $s_{j}$ is strictly increasing, we have
\begin{equation}
\frac{f_{n+1} }{s_{n+1} }<   \frac{f_{n} }{s_{n} }.\label{3.37mmj}
\end{equation}
This shows that if $f=Vh$ then
$
{f_{n} }/{s_{n} }\downarrow$.

\medskip	To see that $(i)$ implies $(ii)$  set
\begin{equation}
\de_{j}= \frac{f_{j}-f_{j-1}}{s_{j}-s_{j-1}}\ge 0.
\end{equation}
We write
\bea
\frac{f_{n}}{s_{n}}&=&\frac{f_{p}+(f_{n}-f_{p})}{s_{p}+(s_{n}-s_{p})}= \frac{f_{p}+\sum_{j=p+1}^{n} (f_{j}-f_{j-1})}{ s_{p}+(s_{n}-s_{p})}\\
&=&\frac{f_{p}+\sum_{j=p+1}^{n}\de_{j}(s_{j}-s_{j-1})}{ s_{p}+(s_{n}-s_{p})} \le \frac{f_{p}+\de_{p }\sum_{j=p+1}^{n} (s_{j}-s_{j-1})}{s_{p}+(s_{n}-s_{p})}\nn\\
&= &\nn\frac{f_{p}+\de_{p } (s_{n}-s_{p})}{s_{p}+(s_{n}-s_{p})}.
\eea
Consequently,
\begin{equation}
\limsup_{n\to\ff}    \frac{f_{n}}{s_{n}}\le \de_{p }.
\end{equation}
Since this holds for all $p$ we see that  $(i)$ implies $(ii)$.
 
  That  $(ii)$ implies $(i)$ is an elementary property of concave functions.
\qed

We now describe the finite  excessive functions for     $ Y$. These are the finite   functions $f$ for which $- Q ({\bf s})f  \ge 0$.  

%moved next display since not used in next proof. } We note for future reference that 
%\bea
%-\( Q({\bf s})f\)_{1}&=&\(a_{1}+a_{2}\)f_{1}  -a_{2}f_{2} \label{suphar.1}\\
%&=&\(\displaystyle\frac1{s_{1}}+\displaystyle\frac1{s_{2} -s_{1}}\)f_{1}  -\displaystyle\frac{f_{2}}{s_{2} -s_{1}}  \nn\\
%&=&  { f_{1} \over  s_{1}}     -{f_{2}-f_{1} \over s_{2}-s_{1}}. \nn
%\eea 

 	\begin{lemma} \label{lem-3.1mm}The function $f$ is a finite excessive function  for $ Y$ if and only if,  
\begin{equation}
\frac{f_{n}-f_{n-1}}{s_{n}-s_{n-1}}\downarrow \de\ge 0, \label{3.45mm}
\end{equation}
   where we take $f_{0}=s_{0}=0$. 
\end{lemma}

\Proof 
  For all $m\ge 1$,   
\bea
-\( Q({\bf s})f\)_{m}&=&-  a_{m} {f_{m-1} } + \( {a_{m}}+ {a_{m+1} }\)f_{m}  -  a_{m+1} {f_{m+1} } \label{suphar.2}\\
&=& a_{m}({f_{m}-f_{m-1}}) - a_{m+1} ( f_{m+1}-f_{m})  \nn\\
&=&  {f_{m}-f_{m-1} \over s_{m}-s_{m-1}}-  {f_{m+1}-f_{m} \over s_{m+1}-s_{m}} \nn.
\eea

Since $-  Q ({\bf s})f  \ge 0$,  this shows that  ${(f_{m}-f_{m-1})}/{(s_{m}-s_{m-1})}$
is decreasing   and consequently has a positive limit which we denote by $\de$.  \qed

We know that  unless $\de=0$, $f$ is not a potential.

\medskip	We sum up these results in the following lemma:  

\begin{lemma}\label{lem-3.3ll}   Let $f$ is a finite excessive function for  $Y$  and set   $g(0)=0$ and  
\begin{equation}
g(s_{j})=  f_{j},\qquad j\in \mN,\label{3.38jj}
\end{equation}
then  $g$ is a concave function   on    $ \{0\}\cup  \{s_{j},\,j\geq 1\} $. 

If in addition 
the function $f$ is a finite  potential for   $Y$ then $g(s_{j})=o(s_{j})$, as $j\to\ff$.

The function  $f=\{f_{j}\}$, where 
\begin{equation}
f_{j}=\de s_{j},\qquad \forall\, j\in \mN,\label{3.39jj}
\end{equation}
is an excessive function for  $Y$,   (in fact $  Q ({\bf s})f \equiv 0$),  but it is not a potential for  $Y$.
\end{lemma}

\Proof    The first   statement  follows   because the terms in      (\ref{suphar.2})    are positive.

The second   statement  follows from
Theorem  \ref{theo-3.1f}, $(ii)$.  

  Obviously  $  Q ({\bf s})f \equiv 0$ so $f$ in (\ref{3.39jj}) is an excessive function for $V$. It follows from the second   statement that it is not a potential.  
 \qed

\bl\label{lem-hh}Let $f$ be a finite  excessive function for   $Y$ such that,  
\begin{equation}
\frac{f_{n}-f_{n-1}}{s_{n}-s_{n-1}}\downarrow 0,\label{3.38mph}
\ee
  where we take $f_{0}=s_{0}=0$. 
Then $f=Vh$ where  $h=-  Q ({\bf s})f \in \ell_{1}^{+}.$
\el

\Proof Since  $f$ is finite and excessive,  
\be
h_{k} =(-  Q ({\bf s})f)_{k}\geq 0,\qquad \forall\,k\in\mN.
\ee
By  (\ref{suphar.2}) and (\ref{3.38mph}) we see that   
\be
\|h\|_{1}=\frac{ f_{1} }{s_{1}} .\label{3.37nn}
\ee

It remains to show that $f=Vh$. 
Using   (\ref{suphar.2}), and setting $f_{0}=s_{0}=0$, we have   that for $n\ge 2$,
\begin{eqnarray}
\sum_{k=1 }^{n-1}s_{k}h_{k}&=&\sum_{k=1 }^{n-1}s_{k}\({f_{k}-f_{k-1} \over s_{k}-s_{k-1}}-  {f_{k+1}-f_{k} \over s_{k+1}-s_{k}}\)
\\
\qquad&=&f_{1} + \sum_{k=2 }^{n-1}(s_{k}-s_{k-1})\({f_{k}-f_{k-1} \over s_{k}-s_{k-1}}\)-s_{n-1}\(  {f_{n}-f_{n-1} \over s_{n}-s_{n-1}}\)\nonumber\\
\qquad &=&f_{n-1}  -s_{n-1}\(  {f_{n}-f_{n-1} \over s_{n}-s_{n-1}}\).\nonumber
\end{eqnarray}
 Furthermore. using  (\ref{suphar.2}) and (\ref{3.38mph}) we see that,
\begin{equation}
\sum_{k=n }^{\ff}h_{k}=  \frac{f_{n }-f_{n-1}}{s_{n }-s_{n-1}}.\label{3.37mm34}
\end{equation}
Consequently,  that for $n\ge 2$,
\begin{eqnarray}
(Vh)_{n}&=&\sum_{k=1 }^{n-1}s_{k}h_{k}+s_{n}  \sum_{k=n }^{\ff}h_{k}
\\
&=& f_{n-1}  -s_{n-1}\(  {f_{n}-f_{n-1} \over s_{n}-s_{n-1}}\)+s_{n}   \frac{f_{n }-f_{n-1}}{s_{n }-s_{n-1}}=f_{n}.\nonumber
\end{eqnarray}
In addition, by (\ref{3.37nn}), 
\begin{equation}
  (Vh)_{1}= s_{1}\sum_{k=1 }^{\ff}h_{k}=f_{1}.
  \end{equation}
\qed

  The next corollary sums up the results of    Theorem \ref{theo-3.1f} and the following lemmas.

\begin{corollary}\label{cor-2.1}
Let $f\geq 0$ be
a finite function. 
Then $f=Vh$ where  $h  \in \ell_{1}^{+}$ if and only if
 \begin{equation}
\frac{f_{n}-f_{n-1}}{s_{n}-s_{n-1}}\downarrow 0,\label{3.38mphn}
\ee
  where we take $f_{0}=s_{0}=0$. 
\end{corollary}

\Proof It follows from Theorem \ref{theo-3.1f} that if $f=Vh$ for some   $h  \in \ell_{1}^{+}$ then $f$ is finite and (\ref{3.38mphn}) holds.   It follows from Lemma \ref{lem-3.1mm} that (\ref{3.38mphn}) implies that $f$ is a finite excessive function for $Y$. Therefore, using Lemma \ref{lem-hh} we see that 
$f=Vh$ for some  $h  \in \ell_{1}^{+}$.\qed

\medskip	We have the following  Riesz decomposition theorem for functions which are excessive for $  Y$.  

\bt \label{theo-3.2}Let $f$ be
a finite   excessive function for $  Y$. Then, necessarily,  $f$  satisfies (\ref{3.45mm}) for some $\de\ge 0$, and 
\begin{equation}
f_{n}=  \wt f_{n}+\de s_{n}\label{3.36jj},\quad\mbox{for all}\quad n\in \mN,
\end{equation} 
where $ \{\wt f_{n},n\in\mN\}$ is
a potential  for   $Y$. 
\et

	\Proof
	Let $f$ be an excessive function for $Y$ and define 
\begin{equation}
\wt f_{n}=f_{n}-\de s_{n} 
\end{equation}
for $\de\geq 0$ as defined in (\ref{3.45mm}).
This implies that $\wt f$ is a finite   excessive function for   $Y$, and 
\begin{equation}
\frac{\wt f_{n+1}-\wt f_{n}}{s_{n+1}-s_{n}}=\frac{f_{n+1}-f_{n}}{s_{n+1}-s_{n}}-\de,\label{3.48pp}
\end{equation}
which together with (\ref{3.45mm}) gives
\be
\lim_{n\to\ff}\frac{\wt f_{n}-\wt f_{n-1}}{s_{n}-s_{n-1}}\downarrow 0\label{3.46mm }.
\ee
By Lemma \ref{lem-hh} we see that $\wt f$ is a potential for   $Y$.  \qed

We now consider the asymptotic properties of permanental processes with kernels that are not symmetric but are modifications of symmetric potentials.
  Let    $\wt  Y_{\al}=\{\wt Y_{\al, n}, n\in\mN\}$ be a   permanental processes with kernel $\wt V=\{\wt  V_{j,k},j,k\in\mN\}$ where,
\begin{equation}
\wt  V_{j,k}=s_{j}\wedge s_{k}+f_{k}.
\end{equation}
and   $f=\{f_{k}, k\in \mN\}$ is    a finite potential for $Y$. 

Since we  use Theorem \ref{theo-mchain} to find the asymptotic behavior of   $\wt  Y_{\al}$ we need only deal with finite sections of kernels.

\begin{lemma} \label{lem-3.1}Let $V(1,n)$ be an  $n\times n$ matrix with elements
\begin{equation}
V(1,n)_{j,k}=s_{j}\wedge s_{k}\qquad j,k=1,\ldots,n ,
\end{equation}
in which  $s_{j}$ is a strictly increasing sequence. Then
\be V(1,n)^{-1}= \left (
\begin{array}{ cccccc cc}  
{ a_{1}}+ { a_{2}}&- { a_{2}}&0&\dots &0 &0  \\
- { a_{2}}& { a_{2}}+ { a_{3}}&- { a_{3}}&\dots &0&0  \\
\vdots&\vdots&\vdots&\ddots&\vdots &\vdots  \\
0&0&0&\dots & { a_{n-1}}+ { a_{n}}&- { a_{n}}   \\
0&0&0&\dots &- { a_{n}}& { a_{n}} \end{array}\right ),\label{harr2.211w}
  \ee
  where $\{a_{j}\}$ is given in (\ref{3.8bb}).
  \el

\Proof    It is easy to verify that this is the inverse of $V(1,n)$.\qed

Note that the first n-1 rows of $V(1,n)^{-1}$ are the same as the first n-1 rows of the matrix in (\ref{harr2.211}).

\begin{lemma} \label{lem-3.2}Let $V(l,n)$ be an  $n\times n$ matrix with elements
\begin{equation}
V(l,n)_{j,k}=s_{j}\wedge s_{k}\qquad j,k=l+1,\ldots,l+n ,
\end{equation}
in which  $s_{j}$ is a strictly  increasing sequence.  Then    
  \[
\hspace{-4.3in}V(l,n)^{-1}=
\]
\be   \left (
\begin{array}{ cccccc c }  
1/s_{ l+1  }+ {a _{l+2}}&- {a _{l+2}}&0&\dots &0 &0  \\
- {a _{l+2}}& {a _{l+2}}+ {a _{l+3}}&- {a _{l+3}}&\dots &0&0 \\
\vdots&\vdots&\vdots&\ddots&\vdots &\vdots  \\
0&0&0&\dots & {a _{l+n-1}}+ {a _{l+n}}&- {a _{l+n}}    \\
0&0&0&\dots &- {a _{l+n}}& {a _{l+n} } \end{array}\right ).\label{harr2.65mm}
  \ee

  \el

\Proof This follows immediately from Lemma \ref{lem-3.1} by relabeling the $a_{\cd}$ and    taking  $a_{1}=1/s_{l+1} $.    An alternate proof is simply to verify  that (\ref{harr2.65mm}) is the inverse of $V(l,n)$.\qed

In the next lemma we give the estimate that enables us to apply Theorem \ref{theo-mchain}.    Recall that for any invertible matrix $M$ we often denote $M^{-1}_{j,k}$ by $M^{j,k}$. 

\begin{lemma} \label{lem-3.5jj}
Let $f$ be a potential for $Y$.  Then  
\begin{equation}
\sum_{j,k= 1}^{ n}V(l,n)^{k,j}f_{l+ j}=  o_{l}(1),\label{3.12ff}
\end{equation}
uniformly in $n$.
\end{lemma}

\Proof   
Note that  
\begin{equation}
\sum_{j,k= 1}^{ n}V(l,n)^{k,j}f_{l+ j}= \sum_{ j= 1}^{ n}f_{l+j}   \sum_{k= 1}^{ n} V(l,n)^{k,j}=\frac{f_{l+1}}{ s_{l+1}},\label{3.12mm}
\end{equation}
where we use the fact that all the column sums of $V(l,n)^{-1}$ are equal to zero except for the first one. Therefore, (\ref{3.12ff}) follows from (\ref{3.38mm}).\qed

\noindent{\bf Proof of  Theorem \ref{theo-3.1} }  We first use Theorem \ref{theo-mchain}.   Therefore, we need    to obtain the denominator in (\ref{121.4a}) for the Gaussian sequence $\xi=\left\{\xi_{j},j\in \mN\right\}$ where, 
\be
E(\xi_{j}\xi_{k})=s_{j}\wedge s_{k},\qquad j,k\in \mN. \label{3.14k}
\ee
We use Koval's Theorem,   \cite[page 1]{Kovalb}.  This involves the  function 
\begin{equation}
\KK_{s_{i},M}(j)=\log \(\sum_{i=1}^{j-1}      M\wedge  \log \(s_{i+1}/s_{i}\) \),\label{ffy}
\end{equation}
for any number $M>0$. (Note that in the notation introduced in (\ref{ff}),    $\KK_{s_{i}}(j)=\KK_{s_{i},1}(j)$.) Since for   any  $M>0$,
\begin{equation}
 \lim_{j\to\ff} \frac{\KK_{s_{i} }(j)}{\KK_{s_{i},M}(j)} =1,\ \label{3.73y}
   \end{equation}
we use $\KK_{s_{i} }(j)$    to avoid ambiguity.

\medskip	 Koval's Theorem states that  
\begin{equation}
\varlimsup_{j\to \ff} \frac{\xi_{j} }{  (2s_{j}\KK_{s_{i} }(j))^{1/2}}=1,   \qquad 
a.s.\label{18.3}
\end{equation}

\medskip  	Note that for  any  $M>0$,
\begin{equation}
   \lim_{j\to\ff}\KK_{s_{i},M}(j)=\ff.\label{3.73jb}
   \end{equation}
This is obvious when $\limsup_{i\to\ff}\log \(s_{i+1}/s_{i}\) >M$ because there would be an infinite number of the terms $M$ in the sum. If $\limsup_{i\to\ff}\log \(s_{i+1}/s_{i}\) \le M$, then  replacing $M$ by $2M$, we can find an $i_{0}$ such that,
\begin{equation}
   \sum_{i=1}^{j-1}      2M\wedge  \log \(s_{i+1}/s_{i}\) > \sum_{i=i_{0}}^{j-1}         \log \(s_{i+1}/s_{i}\) =\log s_{j }-\log s_{i_{0}}.
   \end{equation}
By (\ref{3.73y}), $ \lim_{j\to\ff}  \KK_{s_{i},2M }(j)/\KK_{s_{i},M}(j) =1$ so we get (\ref{3.73jb}).

\medskip	By Theorem \ref{theo-borelN}, $ \wt   V $ is the kernel of $\al$-permanental processes for all $\al>0$. In addition we see by Lemma \ref{lem-3.5jj} that (\ref{condreal}) is satisfied, and by (\ref{3.38mm})   and (\ref{3.73jb}) that (\ref{1.5mm}) is satisfied. Consequently, we can use Theorem \ref{theo-mchain} to get (\ref{18.4mkw}) for all   $1/2\le \al<\ff$. Since $\wt Y_{\al}$ is infinitely divisible and positive, it is obvious that the upper bound in (\ref{18.4mkw})  holds for all $\al>0$. 

\medskip

  We now show that the lower bounds in (\ref{18.4mkw})  holds for all $\al>0$. To   show this  it suffices to find a subsequence $\{s_{p_{j}}\}$ of $\{s_{ {j}}\}$ such that  
 \begin{equation}
\limsup_{j\to \ff} { \,\wt Y_{\al, p_{j}} \over   s_{p_{j}} \KK_{s_{i}}(p_{j})}\ge 1,  \qquad 
a.s.\label{m2.72mm},   \qquad \forall\, \al>0.
\end{equation} 
We choose $\{s_{p_{j}}\}$ recursively as follows:
\begin{equation}
   s_{p_{j+1}}=\min\left\{k:\frac{s_{k}}{s_{p_{j}}}\ge \th\right\}\label{2.57mm}
   \end{equation} 
 where $s_{p_{1}}=1$ and $\th\ge e$. Clearly  
 \begin{equation}
   \frac{s_{p_{j+1}}}{s_{p_{j}}}\ge \th,\qquad\mbox{and}\qquad \frac{ s_{p_{j+1}-1}}{s_{p_{j}}}<\th .\label{2.74mm}
   \end{equation}
   Consequently,  
 \bea
   \sum_{k= p_{l}}^{ p_{l+1}-1}1\wedge \log\frac{s_{k+1}}{s_{k}}&\le& \(   \sum_{k=p_{l}}^{ p_{l+1}-2}  \log\frac{s_{k+1}}{s_{k}}\)+1 \\&\le&\nn\log \th+1<2\log \th.
   \eea
 It follows from these relationships that,   \bea
    \KK_{s_{i}}(p_{j})&=&\log\(\sum_{i=1}^{p_{j}-1}1\wedge \log \frac{s_{i+1}}{s_{i}}\)\\
    &=&\nn\log\(\sum_{l=1}^{j}\sum_{i=p_{l}}^{p_{l+1}-1}1\wedge \log \frac{s_{i+1}}{s_{i}}\)\\
    &\le&\nn \log j+\log  \log\th^{2}. 
   \eea
   This shows that 
  \begin{equation}
\limsup_{j\to \ff} { \,\wt Y_{\al, p_{j}}  \over   s_{p_{j}} \KK_{s_{i}}(p_{j})}\ge  \limsup_{j\to \ff} { \,\wt Y_{\al, p_{j}}  \over   s_{p_{j}}\log j}. \label{2.77mm}\end{equation} 
   Therefore, to obtain  (\ref{m2.72mm}) it suffices to show that 
\begin{equation}
   \limsup_{j\to \ff} { \,\wt Y_{\al, p_{j}}  \over   s_{p_{j}}\log j}\ge 1.\label{2.62mm}
   \end{equation}
 To do this we first extend  and relabel
  $\wt Y_{\al}^{({p})} =\{\wt Y_{\al, p_{j}} ,j\in  \mN\}$ to the  permanental process  $\wh Y_{\al}^{({p})}=  \{\wh Y_{\al, j},j\in \{0\}\cup\mN\}$   with kernel,
\begin{equation}
K_{j,k}=s_{p_{j}}\wedge s_{p_{k}}+f(s_{p_{k}}),\qquad j,k\in\mN,\label{3.56nny}
\end{equation} 
\vspace{-.25in}
\begin{equation}
K_{0,0}=1, \qquad K_{j,0}=1, \qquad j\in\mN, \qquad \mbox{and }\qquad  K_{0,k}=f(s_{p_{k}}),\qquad k\in\mN.      \nn
\end{equation}
 It is clear that   $\wh Y_{\al}^{({p})}\stl \wt Y_{\al}^{({p})} $ on $\mN$, so that to obtain  (\ref{2.62mm}) it suffices to show that, 
\begin{equation}
   \limsup_{j\to \ff} { \, \wh Y_{\al, j} \over   s_{p_{j}}\log j}\ge 1.\label{2.64mm}
   \end{equation}
(Note that by definition, to show that $\wh Y_{\al}^{({p})}$ is a permanental process it suffices to show that for all $\{i_{1},\ldots,i_{n}\}\in  \{0\}\cup\mN$,   $\{K_{i_{j},i_{k}}\}_{j,k=0}^{n}$ is the kernel of a permanental process.   It follows as in  in  (\ref{10.6mm})--(\ref{19.40}) that $\{K_{i_{j},i_{k}}\}_{j,k=0}^{n}$ is an inverse M-matrix. Hence by  \cite[Lemma 4.2]{EK}  it is the kernel of a permanental process.) 

\medskip		 Let $V^{({p})}=\{V^{({p})}_{j,k},j,k\in\mN\}$ where,
\begin{equation}
V^{(p)}_{j,k}=s_{p_{j}}\wedge s_{p_{k}} .\label{2.65mm}
\end{equation}

Let  $K(0,n+1)$ denote the matrix  $\{K_{j,k}\}_{j,k=0}^{n}$.   It follows from   (\ref{19.40})   that  for $j\geq 1 $ the reciprocal of  the diagonal element of the $j$-th row of $(K(0,n+1))^{-1}$, i.e., $1/K(0,n+1)^{j,j}$, satisfies,
\begin{equation}
    1/K(0,n+1)^{j,j}=1/ V^{(p)}(1,n)^{j,j},\qquad 1\le j\le n.
 \end{equation}

It follows from Lemma \ref{lem-3.1} with $s_{j}$  replaced by $s_{p_{j}}$,  and   the second equality in (\ref{harr2.c}),
 that for   $1\le j< n$, 
\begin{eqnarray}
1/V^{(p)}(1,n)^{j,j}&=&{\(s_{p_{j+1}}-s_{p_{j }}\)\(s_{p_{j }}-s_{p_{j-1}}\) \over s_{p_{j+1}}-s_{p_{j-1}}}   \label{3.58nn}\\
&=&s_{p_{j }}\({\(1-s_{p_{j }}/s_{p_{j+1 }}\)\(1-s_{p_{j -1}}/s_{p_{j }}\) \over 1-s_{p_{j -1}}/s_{p_{j +1}}}\)   \nonumber\\
&\ge & s_{p_{j }}  \(1-s_{p_{j }}/s_{p_{j+1 }}\)\(1-s_{p_{j -1}}/s_{p_{j } }\) .     
\nn
\end{eqnarray}
Using (\ref{2.74mm}) we see that    for     $1\le j< n$,
\begin{equation}
   1/K(0,n+1)^{j,j}\ge s_{p_{j}}(1-1/\th)^{2}.  \label{3.58ynn}
   \end{equation}
Since this  holds for all $n$, it follows from
\cite[Lemma 7.3]{MRejp} that,  
\begin{equation}
\varlimsup_{j\to \ff} {\wh Y_{\al, j}\over   s_{p_{j}}\log j}\geq \(1-(1/\th) \)^{2} \qquad 
a.s.,\label{18.newe}
\end{equation}
and since we can take $\th$ arbitrarily large we get (\ref{2.64mm}).\qed

    We   continue to study the behavior of the function $\KK_{s_{i}}(j).$

 \bl \label{lem-2.7}     
\begin{equation}
   \KK_{s_{i}}(j)\le \log\log s_{j}\wedge \log j.\label{3.69y}
   \end{equation}
   Furthermore, if 
   \begin{equation}
   \liminf_{i\to\ff}\frac{s_{i}}{s_{i-1}}>1,\qquad\mbox{then}\qquad \label{3.70y}\lim_{j\to\ff}   \frac{  \KK_{s_{i}}(j)}{\log j}=1
   \end{equation}
   and if  \begin{equation}
  \limsup_{i\to\ff}\frac{s_{i}}{s_{i-1}}<\ff \qquad\mbox{then}\qquad  \lim_{j\to\ff}   \frac{  \KK_{s_{i}}(j)}{\log \log s_{j}}=1.\label{2.64}
   \end{equation}
  \el

\Proof The statement in (\ref{3.69y}) is trivial. To continue, consider $\KK_{s_{i},M}(j)$ in (\ref{ffy}). 
If $   \liminf_{i\to\ff}\frac{s_{i}}{s_{i-1}}>1$ holds there exist numbers   $m_{0}>0$ and   $i_{0}$ such that, 
\begin{equation}
 \inf_{i\ge i_{0}}   \frac{s_{i}}{s_{i-1}}\ge e^{m_{0}},
   \end{equation}
which implies that,
\begin{equation}
 \inf_{i\ge i_{0}}\log    ({s_{i}}/{s_{i-1}})>  m_{0}>0.
   \end{equation}
Therefore,
\begin{equation}
\lim_{j\to\ff}   \frac{  \KK_{s_{i},  m_{0}}(j)}{\log j}=1,
   \end{equation}
which, by (\ref{3.73y}) gives (\ref{3.70y}).

  To get (\ref{2.64}) we simply take $M=1+\log D$ in (\ref{ffy}), where $  D= \limsup_{i\to\ff} {s_{i}}/ s_{i-1}.$
\qed

When $   \liminf_{i\to\ff} {s_{i}}/{s_{i-1}} =1$, we can't simplify $ \KK_{s_{i}}(j)$ without imposing additional conditions. It can oscillate between $\log j$ and $\log\log s_{j}$ when $\log\log s_{j}<\log j$. (Of course it is possible that  $\log\log s_{j}>\log j$ for some $j$,  or even for most $j$, but because of (\ref{3.69y}) we needn't be concerned with these cases.)

\medskip	We can be more precise when,  
 \be \lim_{j\to\ff}\frac{s_{j+1}}{s_{j}}=1.
  \ee
 We can write,  
\begin{equation}
 s_{j}=\exp\({\sum_{k=1}^{j}\ep_{k}}\),\qquad\mbox{where $\ep_{k}>0$,  $ \lim_{k\to\ff}\ep_{k}= 0$},
 \end{equation}   
and the sum diverges.   
Since $\ep_{k}\to 0$ we have $\log\log s_{j}<\log j $  for all $j$ sufficiently large, but we may still have\begin{equation}
\lim_{j\to\ff}  \frac{\log\log s_{j}}{\log j}=1.\label{3.64mmj}
 \end{equation} 
This is the case if $\ep_{k}=1/\log k$, which implies that   $s_{j}\sim \exp({ j/\log j})$, as $j\to\ff$, and the right-hand side of (\ref{3.70y}) still holds.

\medskip  	We give some more examples.
 
\begin{example} {\rm 
\begin{itemize}\item[(i)]
 If   $\ep_{k}=k^{ \al-1}$, for $0<\al<1$, we have   $s_{j}\sim \exp({j^{\al}/\al}$), as $j\to\ff$, and,
\begin{equation}
\lim_{j\to\ff}  \frac{\log\log s_{j}}{\log j}=\al.\label{3.64mmjx}
 \end{equation}  
Consequently,  by (\ref{2.64}),
\begin{equation}
\lim_{j\to\ff}   \frac{  \KK_{s_{i}}(j)}{\log \log s_{j}} =\lim_{j\to\ff}   \frac{  \KK_{s_{i}}(j)}{\al \log  {j}}=1.\label{3.71z}
   \end{equation}
 
 \item[(ii)] If $\ep_{k}=k^{-1}$, we have   $s_{j}\sim   j$ as $j\to\ff$, and,
 
  \begin{equation}
\lim_{j\to\ff}   \frac{  \KK_{s_{i}}(j)}{\log \log s_{j}} =\lim_{j\to\ff}   \frac{  \KK_{s_{i}}(j)}{\log  \log  {j}}=1.\label{3.71zx}
   \end{equation}
 
\item[(iii)]  If $\ep_{k}=1/(k\log k) $, we have   $s_{j}\sim  \log j$ as $j\to\ff$, and,
 
  \begin{equation}
\lim_{j\to\ff}   \frac{  \KK_{s_{i}}(j)}{\log \log s_{j}} =\lim_{j\to\ff}   \frac{  \KK_{s_{i}}(j)}{\log  \log \log  {j}}=1.\label{3.71zr}
   \end{equation}
\end{itemize}

 }\end{example}

 	  	\section{Birth and death processes with emigration } \label{sec-gen}

 A continuous time birth and  death process with emigration is  a Markov chain with a tridiagonal $Q$ matrix. When all the row sums of the $Q$ matrix, except for the first row sum, are equal to zero,  it is called, simply, a birth and  death process. In this section we generalize the Q matrix  $Q(\mathbf s)$ defined in (\ref{3.14mm}) to get a large class of $Q$ matrices of continuous time birth and  death process with emigration.

 \medskip	  For any    sequence $b=(b_{1},b_{2},\ldots)$ define $D_{b}=\mbox{diag }(b_{1},b_{2},\ldots)$. We have the following obvious but important lemma:
   
   \begin{lemma} \label{lem-4.1}
Let  $Q$ denote  the Q-matrix of a    Markov chain $Y$  on $\mN$.  If   $b$ is an excessive function for $Y $, then 
\begin{equation}
   D_{b}QD_{b}
   \end{equation}  is also a    Q-matrix. 
    \end{lemma} 

   \Proof This follows immediately since 
$b$ is positive and $Qb\leq 0$.\qed

We apply Lemma \ref{lem-4.1} to   $Q(\mathbf s)$ defined in (\ref{3.14mm}). We point  out in the paragraph containing (\ref{3.14mm}) that  $Q(\mathbf s)$ is the $Q$-matrix of a 
 continuous time Markov chain   $Y$
 with  potential densities   $V=\{V_{j,k},j,k\in \mN\}$ where,
\begin{equation}
V_{j,k}=  s_j\wedge s _{k} ,\qquad j,k\in \mN,\label{3.28kjam}
\end{equation} 
and  $\mathbf s=\{s_{j},\,j\geq 1\} $ is  a    strictly  increasing sequence  with $s_{j}>0$ and   $\lim_{j\to\ff}s_{j}=\ff$. The next lemma is a significant generalization of this observation.

\bl\label{lem-specmar} Let   $Z=\{ Z_{t},t\in R^{+}\}$ be a   continuous symmetric transient Markov chain    on $\mN$ with $Q$ matrix  $D_{b}  Q(\mathbf s)D_{b}$, where $b$ is a finite potential for  the Markov chain  $Y$ defined in  (\ref{3.12nn}). Then   $W=\{W_{j,k},j,k\in \mN\}$ where
   \begin{equation}
W_{j,k}={1 \over b_{j}}V_{j,k}{1 \over b_{k}}\label{gend.2},    \end{equation}
   is the potential density for  $Z$.  
\el

\begin{remark} \label{rem-4.1}{\rm   
In Lemma \ref{lem-specmar} we take $b$ to be a finite potential for  $Y$. It follows from Theorem  \ref{theo-3.1f} that the function   $g(s_{j})=b_{j}$ is an increasing   concave function of $\{s_{j}\}$   and  $s_{j}/b_{j}\uparrow\ff$.

Consider  $\{f_{j}\}$, the finite potentials of $W$. We have
\begin{equation}
   f_{j}=\sum_{k=1}^{\ff}W_{j,k}h_{k}=\sum_{k=1}^{\ff} {1 \over b_{j}}V_{j,k}{h_{k} \over b_{k}}.
   \end{equation}
   Consequently
\begin{equation}
  b_{j} f_{j}=\sum_{k=1}^{\ff}W_{j,k}h_{k}=\sum_{k=1}^{\ff} V_{j,k}{h_{k} \over b_{k}}.\label{3.5mm}
   \end{equation}
Therefore,  $\{b_{j}f_{j}\}$, is a finite potential for $Y$. As noted in the first paragraph of this remark this implies that $g(s_{j})=b_{j}$ is an increasing   concave function of $\{s_{j}\}$. Therefore we can write $f$ as
  \begin{equation}
  f_{j}=\frac{g(s_{j})}{h(s_{j})},\qquad \forall j\in \mN,\label{4.18jj}
  \end{equation}
  where $g$ and $h$ are positive  strictly concave functions.   
 }\end{remark}

\noindent{\bf Proof of Lemma \ref{lem-specmar} }  It is easy to see that   $Q(\mathbf s)V=-I$ in the sense of multiplication of infinite matrices. Consequently, since $   D_{b}W=VD_{b}^{-1}$, it follows that we also have  
\begin{equation}
D_{b}  Q(\mathbf s)D_{b} W=D_{b}  Q(\mathbf s)VD_{b}^{-1}  =-I.
   \end{equation}
   Let $\ov W$ be the potential density for $Z$.   Using  Lemma \ref{lem-9.1mm} we see that $D_{b}  Q(\mathbf s)D_{b}\ov  W\newline =-I$.  Consequently,  \begin{equation}
Q(\mathbf s)D_{b}\(W-\ov  W\)=0. \label{specmar.1}
\end{equation}

Consider the equation $Q(\mathbf s)g=0$.  Using  (\ref{suphar.2}) we see that we must have 
\begin{equation}
   \frac{g_{j}-g_{j-1}}{s_{j}-s_{j-1}}=c_{0}\qquad \forall j\ge 1,
   \end{equation}
for some fixed   constant $c_{0}$ where we set $g_{0}=s_{0}=0$. Therefore, all solutions of $ Q(\mathbf s)g=0$ are of the form $g=c_{0}(s_{1}, s_{2}, \ldots)$.  

Consider   the components of (\ref{specmar.1}). We see that for all $k\in\mN$,
\begin{equation}
   \sum_{j}(Q(\mathbf s))_{l,j} (D_{b}  (W-\ov  W) )_{j,k}=0.
   \end{equation}
Therefore, using the observations in the preceding paragraph, we have that for each $k\in \mN$,
\begin{equation}
W_{j,k}-\ov  W_{j,k}=c_{k}s_{j}/b_{j}, \hspace{.2 in} \forall j\geq 1,\label{specmar.2}
\end{equation}
for   some constant $c_{k}$.

We now show  that  $c_{k}=0$ for all $k$. Let  $ P^{j}$ denote probabilities for $  Z$. We have  
\begin{equation}
\ov  W_{j,k}=  P^{j}\(T_{k}<\ff\)\ov  W_{k,k}.\label{specmar.3}
\end{equation}
Using this and (\ref{specmar.2}) we see that,  
 \begin{equation}
{s_{k} \over b_{j}b_{k}}=W_{j,k}=  P^{j}\(T_{k}<\ff\)\ov  W_{k,k}   +c_{k}s_{j}/b_{j}, \hspace{.2 in} \forall j\geq k.\label{specmar.5}
\end{equation} 
Since $b_{j} $ is increasing and $s_{j}/b_{j}\uparrow\ff$, this is only possible if   $c_{k}=0$.
\qed

 \begin{remark} {\rm 
This Lemma also applies if $b$ is a general finite excessive function for $Y$.  That is, by Theorem 
\ref{theo-3.2}, if we add $\de s_{j}$ to the present $b_{j}$. In that case, the left-hand side of (\ref{specmar.5}) goes to zero as $j\to\ff$,   and the last term in (\ref{specmar.5}) converges to $c_{k}/\de$, which again shows that $c_{k}=0$.  
}\end{remark}

 \begin{remark} {\rm 
We also note that if $Y=\{ Y_{t},t\in R^{+}\}$ is a process with   Q-matrix $Q(\mathbf s)$, then 
$Z$ can be obtained from $Y$  by first doing a $b$-transform   and then a time change by the inverse of  the  continuous additive functional  $\int_0^t b(Y_s)^{-2}\,ds$.   This gives an alternate proof of   Lemma \ref{lem-specmar}.
}\end{remark}

Our goal in this section is to use Theorem \ref{theo-mchain} to prove Theorem \ref{theo-4.2mm}. We use the next   two lemmas to obtain (\ref{condreal}).

 \begin{lemma}\label{lem-4.4mm1} Let $f=W  h$, where $  h\in\ell_{1}^{+}$. Then  
\begin{equation}
\sum_{ j,k=1}^{n}(W(l,n))^{j,k}f_{k+l}= \frac{f_{l+1}}{W_{l+1,l+1}}+o_{l}\(1\), \mbox{ uniformly in }n.\label{condrealq}
\end{equation} 

\end{lemma}

\Proof  
For any $l,k\in\mN $ we set $f^{(l)}_{k}=f_{l+k}$.  Similarly, we set $a^{(l)}_{k}=a_{l+k}$ and $b^{(l)}_{k}=b_{l+k}$. We have  
\be 
     \sum_{j,k= 1}^{ n}W(l,n)^{ j,k} f^{(l)}_{k}=  \sum_{ k= 1}^{ n}f^{(l)}_{k}b^{(l)}_{ {k}}\sum_{j=1}^{n}(V(l,n)) ^{ k,j} b^{(l)}_{j}.   \label{4.29}
  \ee 
  For any sequence $\{c_{k}\}$ we use the standard notation $\De c_{k}=c_{k+1}-c_{k}.$ 

Using  (\ref{harr2.65mm}) we see that, 
  \begin{equation}
     \sum_{j=1}^{n}(V(l,n)) ^{ 1,j} b^{(l)}_{j}=\frac{b^{(l)}_{1}}{s ^{(l)}_{1}}- \De b^{(l)}_{1}a^{(l)}_{2 },
  \end{equation}
and for $1<k<n$,
\bea
  \sum_{j=1}^{n}(V(l,n)) ^{ k,j} b^{(l)}_{j} &=&-b^{(l)}_{k-1}a^{(l)}_{k }+b^{(l)}_{k}(a^{(l)}_{k }+a^{(l)}_{k+1})-b^{(l)}_{k+1}a^{(l)}_{k+1}\\
   &=& \De b^{(l)}_{k-1}a^{(l)}_{k }-\De b^{(l)}_{k}a^{(l)}_{k+1 }\nn,
    \eea  
    and  
    \begin{equation}
     \sum_{j=1}^{n}(V(l,n)) ^{ n,j} b^{(l)}_{j}=\De b^{(l)}_{n-1}a^{(l)}_{n }.\label{4.31s}
  \end{equation} 
  It follows from (\ref{4.29})--(\ref{4.31s}) that,  
  \be 
     \sum_{j,k= 1}^{ n}W(l,n)^{ j,k} f^{(l)}_{k}=\(f^{(l)}_{1}  b^{(l)}_{1}\) \frac{b^{(l)}_{1}}{s_{1} ^{(l)}}+\sum_{k=2}^{n}a^{(l)}_{k}\De b^{(l)}_{k-1}\De  \(f^{(l)}_{k-1}  b^{(l)}_{k-1}\).   \label{4.29p}
  \ee

  Set   $\wt f  =D_{b}f$.  Then, since $f=Wh$, for some $h\in\ell_{1}^{+}$, we see that 
   \begin{equation}
   \wt f_{k}=b_{k}\sum_{j=1}^{\ff}W_{k,j}h_{j}=\sum_{j=1}^{\ff}V_{k,j}\frac{h_{j}}{b_{j}},\qquad\forall\, k\in\mN.
   \end{equation}
   This shows that $\wt f=V(D^{-1}_{b}h)$. Therefore, by (\ref{3.37mm}), we see that
   for all $n\ge 1$,
\begin{equation}
\frac{\wt f_{n+1}-\wt f_{n}}{s_{n+1}-s_{n}}= \sum_{k=n+1}^{\ff}\frac{h_{k}}{b_{k}}.\label{4.20mm}
\end{equation}

We now use  (\ref{4.29p})  and  (\ref{4.20mm}) and the fact that $\wt f  =D_{b}f$ to get,  
   \bea
     \sum_{j,k= 1}^{ n}W(l,n)^{ j,k}  f^{(l)}_{k}&= &  f^{(l)}_{1}  \frac{\(b^{(l)}_{1}\)^{2}}{s_{1} ^{(l)}}+\sum_{j=2}^{n}a^{(l)}_{j}\De b^{(l)}_{j-1}\De \wt f^{(l)}_{j-1} \label{4.33mm}\\
     &=&\nn   f^{(l)}_{1}  \frac{\(b^{(l)}_{1}\)^{2}}{s_{1} ^{(l)}}+ \sum_{j=2}^{n}\frac{\De b^{(l)}_{j-1}\De   \wt f^{(l)}_{j-1}}{s^{(l)}_{j}-s^{(l)}_{j-1}}\\
     &=&\nn  \frac{f_{l+1}}{W_{l+1,l+1}}+\sum_{j=2}^{n} {\De b^{(l)}_{j-1} } \sum_{k=j}^{\ff} {h^{(l)}_{k} \over b^{(l)}_{k}}.      \eea
 
Since 
\bea 
 \lefteqn{    \sum_{j=2}^{n} {\De b^{(l)}_{j-1} } \sum_{k=j}^{\ff} {h^{(l)}_{k} \over b^{(l)}_{k}}}\\
   &&\nn\quad=-b^{(l)}_{1}\sum_{k=2}^{\ff} {h^{(l)}_{k} \over b^{(l)}_{k}}+  \sum_{j=2}^{n-1}b^{(l)}_{j}\(\sum_{k=j}^{\ff} {h^{(l)}_{k} \over b^{(l)}_{k}}-\sum_{k=j+1}^{\ff} {h^{(l)}_{k} \over b^{(l)}_{k}}\)+b^{(l)}_{n}\sum_{k=n}^{\ff} {h^{(l)}_{k} \over b^{(l)}_{k}}\\
   &&\nn\quad\le   \sum_{j=1}^{n-1}h^{(l)}_{j} +b^{(l)}_{n}\sum_{k=n}^{\ff} {h^{(l)}_{k} \over b^{(l)}_{k}}\le \sum_{j=1}^{n-1}h^{(l)}_{j} + \sum_{k=n}^{\ff} {h^{(l)}_{k} }=\sum_{j=1}^{\ff}h^{(l)}_{j},
 \eea 
   we get (\ref{condrealq}).\qed
 
Using the next lemma with    Lemma \ref{lem-4.4mm1} we get (\ref{condreal}).

\begin{lemma}\label{lem-4.4mm2}   Let  $f=Wh$, $h\in \ell_{1}^{+}$. Then
\begin{equation}
  \lim_{j\to\ff}\frac{f_{j}}{W_{j,j}}=0.\label{4.41mm}
  \end{equation}
  \end{lemma}
 
\Proof   
For $p<j+2$,  we have, 
   \bea
f_{j} &=&\sum_{k=1}^{\ff} W_{j,k}h_{k}=\frac{1}{ b _{j}}\sum_{k=1}^{\ff} V_{j,k}\frac{h_{k}}{ b_{k}}\label{3.23mm}\\
 &=&\nn\frac{1}{ b _{j}}\sum_{k=1}^{p} \frac{s_{k}h_{k}}{ b_{k}}+\frac{1}{ b _{j}}\sum_{k=p+1}^{j-1} \frac{s_{k}h_{k}}{ b_{k}}+\frac{s_{j}}{ b _{j}}\sum_{k=j}^{\ff} \frac{ h_{k}}{ b_{k}}\\
 &\le &\nn\frac{1}{ b _{j}}\frac{s_{p}}{ b_{p} }\sum_{k=1}^{p}   h_{k}  +\frac{s_{j}}{ b^{2} _{j}}\sum_{k=p+1}^{\ff}  { h_{k}}=\frac{1}{ b _{j}}\frac{s_{p}}{ b_{p} }\sum_{k=1}^{p}   h_{k}  +W_{j,j}\sum_{k=p+1}^{\ff}  { h_{k}}, 
        \eea
  where for the last line we note that by Remark \ref{rem-4.1},   $b_{k}$ is increasing and $s_{k}/b_{k}\uparrow$.  Therefore, 
   \begin{equation}
   \frac{f_{j}}{W_{j,j}}\le \frac{s_{p}}{b_{p}} \frac{b_{j}}{s_{j}}\|h\|_{1}+\sum_{k=p+1}^{\ff}  { h_{k}}. 
   \end{equation}    
   Using Remark \ref{rem-4.1} again we see that for all $p>0$,
   \begin{equation}
     \lim_{j\to\ff}\frac{f_{j}}{W_{j,j}}\le \sum_{k=p+1}^{\ff}  { h_{k}}.
   \end{equation}
    This gives   (\ref{4.41mm}).\qed

 	  \noindent{\bf Proof of Theorem  \ref{theo-4.2mm}  }   Let $\xi=\{\xi_{j},j\in \mN\}$ be a Gaussian sequence with covariance $W$. It follows from Koval's Theorem that  
  \begin{equation}
  \limsup_{j\to\ff}\frac{ \xi_{j}}{ (2W_{j,j} \KK_{s_{i}}(j))^{1/2}}= 1\quad a.s.\label{3.16kuj}
  \end{equation}
   Therefore, for $\al\ge 1/2$, (\ref{18.4mkff})  follows from Theorem  \ref{theo-mchain}.   Note that Lemmas \ref{lem-4.4mm1} and \ref{lem-4.4mm2} give (\ref{condreal}). In addition  Lemma  \ref{lem-4.4mm2} and  (\ref{3.73jb})  shows that (\ref{1.5mm}) holds. Also, as we have pointed out,    the upper bound in  (\ref{18.4mkff}) actually holds for all $\al>0$.   

\medskip	 We now show that the lower bounds in (\ref{18.4mkff})  holds for all $\al>0$. To this  it suffices to find a subsequence $\{s_{p_{j}}\}$ of $\{s_{ {j}}\}$ such that  
 \begin{equation}
\limsup_{j\to \ff} { \,\wt Z_{\al, p_{j}} \over W_{p_{j},p_{j}} \KK_{s_{i}}(p_{j})}\ge 1,  \qquad 
a.s.\label{2.72nn},   \qquad \forall\, \al>0.
\end{equation} 
If we choose $\{s_{p_{j}}\}$ as in (\ref{2.57mm}), this   follows if we show that,
\begin{equation}
\limsup_{j\to \ff} { \,\wt Z_{\al, p_{j}}\over W_{p_{j},p_{j}} \log j }\ge 1,  \qquad 
a.s.\label{2.72nnc},   \qquad \forall\, \al>0.
\end{equation}

	Consider the permanental process     $\wt Z_{\al}^{({p})} =\{ \wt Z_{\al, p_{j}},j\in  \mN\}$. 
As in the proof of Theorem \ref{theo-3.1}
 we   extend   and relabel 
  $\wt Z_{\al}^{({p})} $ to  get a permanental process    $\wh Z_{\al}^{({p})}=  \{\wh Z_{\al,j  },j\in \{0\}\cup\mN\}$   with kernel $\ov K=\{\ov K_{j,k};j,k\in \{0\}\cup\mN\}$ where,
\begin{equation}
\ov K_{j,k}= \frac{s_{p_{j}}\wedge s_{p_{k}} }{ b_{p_{j}} b_{p_{k}}}+f(s_{p_{k}}),\qquad j,k\in\mN,\label{3.56nn}
\end{equation}
\begin{equation}
\ov K_{0,0}=1, \qquad \ov K_{j,0}=1, \qquad j\in\mN, \qquad \mbox{and }\qquad  \ov K_{0,k}=f(s_{p_{k}}),\qquad k\in\mN.      \nn
\end{equation}

 It is clear that   $\wh Z_{\al}^{({p})}\stl \wt Z_{\al}^{({p})} $ on $\mN$, so that to obtain (\ref{2.72nnc}) it suffices to show that, 
\begin{equation}
   \limsup_{j\to \ff} { \,\wh Z_{\al,j  } \over   W_{p_{j},p_{j}}\log j}\ge 1.\label{2.64mmx}
   \end{equation}

    Let $\ov W^{({p})}=\{\ov W^{({p})}_{j,k},j,k\in\mN\}$ where,  \begin{equation}
\ov W^{({p})}_{j,k}=\frac{s_{p_{j}}\wedge s_{p_{k}} }{ b_{p_{j}} b_{p_{k}}}  ,\label{3.56xx}
\end{equation}
and let  $\ov K(0,n+1)$ denote the matrix $\{\ov K_{j,k}\}_{j,k=0}^{n}$. As in the proof of  Theorem \ref{theo-3.1}, it follows from   (\ref{19.40})   that  for 
$j\geq 1 $,
\begin{equation}
    1/\ov K(0,n+1)^{j,j}=1/ \ov W^{(p)}(1,n)^{j,j},\qquad 1\le j\le n.
 \ee
 It is easy to see that
 \begin{equation}
   \ov W^{({p})}(1,n)^{j,j}=b_{j}^{2}   V^{({p})}(1,n)^{j,j},\label{3.32nn}
   \end{equation}
 where $V^{({p})}$ is given in (\ref{2.65mm}). Therefore, analogous to (\ref{3.58nn}) and (\ref{3.58ynn}) we see that,
 \begin{equation}
   1/\ov K(0,n+1)^{j,j}\ge \frac{s_{p_{j}}}{b_{j}^{2}}(1-1/\th)^{2}=W_{p_{j},p_{j}}(1-1/\th)^{2}
.  
   \end{equation}
 As in the proof of Theorem \ref{theo-3.1} this implies (\ref{2.64mmx}).\qed

 \begin{remark} {\rm Let $\wt Y_{\al}$ be as in Theorem \ref{theo-3.1}. The kernel of $\wt Y_{\al}$ is 
  $
\wt V=\{ \wt V_{j,k};j,k\in\mN\}$, where
   \begin{equation}
  \wt V_{j,k}=V_{j,k}+f_{k},\qquad j,k\ge 1.
   \end{equation}
   
 It follows from (\ref{int.1}) that $  Z'_{\al}:=D_{b}^{-2}\wt Y_{\al}=(b^{-2}_{1}Y_{\al,1 }, b^{-2}_{2}Y_{\al,2},  \ldots)$ has kernel   $
 W'=\{ W'_{j,k};j,k\in\mN\}$ where
   \begin{equation}
 W'_{j,k}=\frac{V_{j,k}}{b_{j}b_{k}}+\frac{f_{k}}{b_{j}b_{k}}= {W_{j,k}} +\frac{f_{k}}{b_{j}b_{k}},\qquad j,k\in\mN.
   \end{equation}
This is because for all $n\times n$ matrices $K$ and $D_{b'}^{-1}$,
\begin{equation}
     |I+ K D_{b'}^{-2}S| =   |I+ D_{b'}^{-1}KD_{b'}^{-1}S|. 
   \end{equation}
 It follows from Theorem \ref{theo-3.1} that 
   \begin{equation}
  \limsup_{j\to\ff}\frac{ b^{2}_{j}  Z'_{\al}(j) }{  V_{j,j} \KK_{s_{i}}(j) }= 1\quad a.s.\label{4.sp},
  \end{equation}
 or, equivalently
      \begin{equation}
  \limsup_{j\to\ff}\frac{ Z'_{\al}(j) }{  W_{j,j} \KK_{s_{i}}(j) }= 1\quad a.s.\label{4.spx}
  \end{equation}
 as in (\ref{18.4mkff}).  
 
 This is easy, but $ Z'_{\al}$ has kernel $W'$ whereas $\wt Z_{a}$ in Theorem \ref{theo-4.2mm} has kernel $\wt W$, in (\ref{4.38mm}). In Theorem \ref{theo-borelN} we set out to consider symmetric kernels perturbed by an excessive function $f$ as in (\ref{1.10}). This is what we do in Theorem \ref{theo-4.2mm}.
 
}\end{remark}

\medskip	In Lemma \ref{lem-4.4mm2} we use the explicit representation of $W$. It is interesting to note that (\ref{4.41mm}) holds in great generality when the diagonals of the matrix go to infinity.

\bl\label{lem-littleo} Let $f=Wh$, $h\in\ell_{1}^{+}$, for some infinite matrix  $W$ such that  
\begin{equation}
 W_{k,j}\leq   W_{k,k},\qquad \forall k\in\mN  .
   \end{equation} 
   Then 
 \begin{equation}
  f_{k}\le  W_{k,k}\,\,\|h\|_{1}, \mbox{ and if $\,\,\lim_{k\to\ff} W_{k,k}=\ff$,}\quad  f_{k}=o( W_{k,k}).\label{aa}
  \end{equation}
 \el

 \Proof  Let $\ep>0$. we have 
 \bea
f_{k}&=&\sum_{\{j: W_{j,j}\le \ep W_{k,k}\}}W_{k,j}h_{j}+\sum_{\{j:W_{j,j}> \ep W_{k,k}\}}W_{k,j}h_{j} \\
  &\le& \ep W_{k,k}\sum_{\{j: W_{j,j}\le \ep W_{k,k}\}}h_{j} +W_{k,k}\sum_{\{j:W_{j,j}> \ep W_{k,k}\}}h_{j} \nn\\
   &\le& \ep W_{k,k}\,\,\|h\|_{1} +W_{k,k}\sum_{\{j:W_{j,j}> \ep W_{k,k}\}}h_{j}. \nn
  \eea
  Therefore,
  \begin{equation}
   \lim_{k\to\ff}\frac{f_{k}}{W_{k,k}}\le  \ep  \,\,\|h\|_{1}+   \lim_{k\to\ff}\sum_{\{j:W_{j,j}> \ep W_{k,k}\}}h_{j}. 
   \end{equation}
  If $\lim_{k\to\ff} W_{k,k}=\ff$ the last sum goes to 0. This gives the second statement in (\ref{aa}). The first statement is obvious. \qed

  Let $\WW=\{\WW_{j,k},j,k\in \mN\}$ where
  \be
 \WW_{j,k}=e^{-|v_{k}-v_{j}|},\qquad \forall\, j,k\in\mN,\label{5.6mmqo}
 \ee
  as defined in (\ref{5.6mmq}). The next theorem applies Theorem \ref{theo-4.2mm} when $W$ is written in this way.  
   
Set  
 \begin{equation}
\ov\KK_{\mathbf v}(j)=\log \(\sum_{i=1}^{j-1}      1\wedge 2 (v_{i+1}-v_{i}) \).\label{ffx}
\end{equation}

	\begin{theorem}
 \label{theo-4.2mmfq}    Let  $\WW$  be the potential density of a continuous time Markov chain $\ZZ$ as given in (\ref{5.6mmqo}) and let $f$ be a finite excessive function for $\ZZ$.  Let  $\wt  \ZZ_{\al}=\{\wt 
   \ZZ_{\al,j},j\in \mN\}$ be the $\al$-permanental process
 with kernel 
   $
   \wt \WW=\{ \wt \WW_{j,k};j,k\in\mN\}$ where,  
   \begin{equation}
   \wt \WW_{j,k}=\WW_{j,k}+f_{k},\qquad  j,k\in\mN. \label{4.38mma}
   \end{equation}

 \begin{itemize}

 \item[(i)]  If  $f= \WW  h$, where $   h  \in\ell_{1}^{+}$, then  
   \begin{equation}
 \limsup_{j\to \ff} { \wt \ZZ_{\al,j} \over   \ov\KK_{\mathbf v}((j)  }=1,  \qquad 
a.s.\label{18.4q}, \qquad \forall\, \al>0,
\end{equation}  
 \item[(ii)]  If  $f= \WW  h$, where $   h  \in\ell_{1}^{+}$ and   $\limsup_{j\to\ff} (v_{j}-v_{j-1})<\ff$,  then  
\begin{equation}
 \limsup_{j\to \ff} { \wt \ZZ_{\al,j} \over    \log v_{j}}=1,  \qquad 
a.s.\label{18.4mkwb}, \qquad \forall\, \al>0.
\end{equation} 

 \item[(iii)] If $\liminf_{j\to\ff}(v_{j}-v_{j-1})>0$ and $f\in c^{+}_{0}$,   then  
\begin{equation}
 \limsup_{j\to \ff} { \wt \ZZ_{\al,j} \over      \log {j }}=1,  \qquad 
a.s.\label{18.4mks}, \qquad \forall\, \al>0.
\end{equation} 
\end{itemize}
Furthermore when   $\liminf_{j\to\ff}(v_{j}-v_{j-1})>0$, the conditions $f\in c^{+}_{0}$, and    $f= \wt \WW  h$,  $   h  \in c^{+}_{0}$, are equivalent. 
    \end{theorem}

 	 \noindent{\bf \Proof of Theorem \ref{theo-4.2mmfq} }
 When    $f= \WW  h$, where $   h  \in\ell_{1}^{+}$, this  is simply an application of  Theorem \ref{theo-4.2mm} with $s_{j}$ replaced by $e^{2v_{j}}$ and $b_{j}=s_{j}^{1/2}$.
 
  That (\ref{18.4mks}) extends to $\al$-permanental processes $\wt  \ZZ_{\al} $ with kernels $\wt \WW$   in which $f$  an excessive function for $\ZZ$   with the property that $f\in c^{+}_{0}$ follows from Theorem \ref{theo-1.8mm} 
 since  $\liminf_{j\to\ff}(v_{j}-v_{j-1})>0$ implies $\|\WW \|<\ff$ and $\WW _{ j,j}=1$ for all $j\in\mN.$ 
 
 The fact that $f\in c^{+}_{0}$ if and only if  $f= \WW   h$, where $   h  \in c_{0}^{+}$ follows from Lemma \ref{lem-8.2nn} once we show that (\ref{8.14}) holds.  The condition $\liminf_{j\to\ff}(v_{j}-v_{j-1})>0$ implies that there exists a $j_{0}$ such that $(v_{j}-v_{j-1})\ge \de>0$ for all $j\ge j_{0}+1$. Therefore, for $j\ge 2j_{0}$,
   \begin{equation}
 \sum_{k=1}^{j/2}\WW _{j,k}\le \sum_{k=1}^{j/2}e^{-(v_{j} -v_{j/2})}\le  \frac{j}{2}\de^{j/2} .\label{8.14aa}
   \end{equation}
This shows that (\ref{8.14}) holds. \qed

 Clearly, the complete statement involving (\ref{18.4mks}) follows from Theorem \ref{theo-1.8mm} and Lemma \ref{lem-8.2nn}. One doesn't need the much more complicated Theorem \ref{theo-mchain}.

 	\begin{example}  \label{ex-4.1a}{\rm Consider the special case of  Theorem \ref{theo-4.2mmfq} $(iii)$ in which  $v_{j}=\la j$, $\la>0$, and  $r=e^{-\la}$ so that, 
   \begin{equation}
\wh\WW_{j,k} = e^{- \la|k-j|}= r^{ |k-j|},\qquad j,k\in\mN.\label{4.4nn}
   \end{equation}
By Lemma \ref{lem-specmar},     $\wh\WW$ is the potential density for the   continuous symmetric transient Markov chain    on $\mN$ with $Q$ matrix  $D_{b}  Q(\mathbf s)D_{b}$, where   $b_{j}=s_{j}^{1/2}=r^{-(j-1)}$.

  This example also follows from Lemma \ref{lem-invpot}. 
We claim that  $\wh\WW$ is the potential density of  a   continuous symmetric transient Markov chain    on $\mN$ with $Q$ matrix
\be -Q  ={1 \over 1-r^{2}}\left (
\begin{array}{ cccccc cc}  
 1 &-r&0&0&\dots   \\
-r&1+r^{2}&-r&0&\dots   \\
0&-r&1+r^{2}&-r &\dots  \\\vdots&\vdots&\vdots&\vdots&\ddots   \end{array}\right ).\label{harr2.211x}
      \ee
To  see this  write out 
\be \wh\WW= \left (
\begin{array}{ cccccc cc}  
 1 &r&r^{2}&r^{3}&\dots   \\
 r&1 & r&r^{2}&\dots   \\
r^{2}& r&1 & r &\dots  \\\vdots&\vdots&\vdots&\vdots&\ddots   \end{array}\right ).\label{harr2.21dd}
      \ee
      It is easily seen that $\wh\WW Q=-I$.    }\end{example}

 Consider the birth and death process studied in Section \ref{sec-BM} which is defined in terms of a  strictly  increasing sequence     $\mathbf s=\{s_{j},\,j\geq 1\} $   with   $s_{1}>0$ and  $\lim_{j\to\ff}s_{j}=\ff$.  This process has  potential densities 
\be
V_{j,k}=s_{j}\wedge s_{k},\qquad j,k\in\mN. 
\ee
We shift the sequence $\mathbf s$
  by a constant $\De>-s_{1}$ and  obtain a new birth and death process defined by the  sequence $\mathbf s'=\{s'_{j}=s_{j}+\De,\,j\in\mN\} $ with    potential density 
  \be
V'_{j,k}=s'_{j}\wedge s'_{k},\qquad j,k\in\mN. 
\ee

 	\begin{lemma} 
 Let  $f=Vh$ for some $h\in \ell_{1}^{ +}$, then  $f=V'h'$ for some $h'\in \ell_{1}^{ +}$ if
\begin{equation}
\frac{f_{1}}{s_{1}+\De}\geq \frac{f_{2}-f_{1}}{s_{2}-s_{1}}.\label{34.b}
\end{equation}   
 \end{lemma}
 
 \Proof Since $s_{j}-s_{j-1}=s'_{j}-s'_{j-1}$, $j\ge 2$,
 \begin{equation}
    \frac{f_{2}-f_{1}}{s'_{j}-s'_{j-1}}= \frac{f_{2}-f_{1}}{s_{j}-s_{j-1}}\qquad j\ge 2.\label{3.59}
    \end{equation}
By Corollary \ref{cor-2.1} the right-hand of (\ref{3.59}) and consequently  the left-hand of (\ref{3.59}) is decreasing. Therefore by  Corollary \ref{cor-2.1} again
for $f=V'h'$ for some $h'\in \ell_{1}^{ +}$ we only need in addition that,
 \begin{equation}
\frac{f_{1}}{s'_{1} }\geq \frac{f_{2}-f_{1}}{s'_{2}-s'_{1}}.\label{34.bq}
\end{equation}  
This is (\ref{34.b}).\qed

\medskip	 The next lemma generalizes Lemma 	 \ref{lem-specmar},

\bl\label{lem-specmarx} If $b=Vh$ for some $h\in \ell_{1}^{ +}$ and  $b_{2}>b_{1}$ and, 
\begin{equation}
   \De\le \frac{b_{1}s_{2}-b_{2}s_{1}}{b_{2}-b_{1}} ,\label{3.63mm}
   \end{equation}
   then   $W'=\{W'_{j,k},j,k\in \mN\}$ where,
    \begin{equation}
 W'_{j,k}={1 \over b_{j}}V'_{j,k}{1 \over b_{k}}=W_{j,k}+{\De \over b_{j} b_{k}},\label{wshift.1}
 \end{equation}
 is the potential density of a Markov chain.
\el 

\Proof  Since $b=Vh$ for some $h\in \ell_{1}^{ +}$ then it follows that if  (\ref{34.b}) holds  with $f$  replaced by $b$, that is if we have,
 \begin{equation}
\frac{b_{1}}{s_{1}+\De}\geq \frac{b_{2}-b_{1}}{s_{2}-s_{1}},\label{34.bb}
\end{equation} 
  then $b=V'h'$ for some $h'\in \ell_{1}^{ +}$ in which case  (\ref{wshift.1}) follows from Lemma \ref{lem-specmar}. 
  The condition in (\ref{3.63mm}) is simply a rearrangement of (\ref{34.bb}).\qed

 \medskip	 We see from  (\ref{harr2.211q}) and (\ref{3.8bb}) that $Q(s')$ differs from $Q(s)$ only in the $(1,1)$ entries which are,\label{page40}  
\begin{equation}
Q(s)_{1,1}=-\frac{1}{s_{1}}-Q(s)_{1,2}\hspace{.2 in}\mbox{   and }\hspace{.2 in}Q(s')_{1,1}=-\frac{1}{s_{1}+\De}-Q(s)_{1,2}.\label{qshift.1}
\end{equation}
  Set
  \begin{equation}
\QQ(b,\mathbf s) = D_{b} Q(\mathbf s)D_{b}.\label{6.1nnw}
   \end{equation}
 This is  the $Q$ matrix for $W$. Consequently, $ \QQ(b,\mathbf s')$ is  the $Q$ matrix for $W'$. Since $s_{j}-s_{j-1}=s'_{j}-s'_{j-1}$ for all $j\ge 2$ and $b$ is unchanged we see that 
     $ \QQ(b,\mathbf s')$ differs from $ \QQ(b,\mathbf s)$  only in the $(1,1)$ entry. Using    (\ref{qshift.1}) and the fact that $\QQ(b,\mathbf s)_{1,2}=b_{1}b_{2}Q(s)_{1,2}$ we have, 
 \begin{equation}
\QQ(b,\mathbf s)_{1,1}=-\frac{b^{2}_{1}}{s_{1}}-\frac{b _{1}}{b_{2}}\QQ(b,\mathbf s)_{1,2}\hspace{.1 in}\mbox{and}\hspace{.1 in}\QQ(b,\mathbf s')_{1,1}=-\frac{b^{2}_{1}}{s_{1}+\De}-\frac{b _{1}}{b_{2}}\QQ(b,\mathbf s)_{1,2}.\label{qshift.2}
\end{equation}
Since $W'$ is a potential we know that $ \QQ(b,\mathbf s')$ is a $Q$ matrix. Therefore the row sum of its first row must be less than or equal to 0. That is we must have,
\begin{equation}
   \QQ(b,\mathbf s')_{1,1}\le - \QQ(b,\mathbf s')_{1,2}=-\QQ(b,\mathbf s)_{1,2}.
   \end{equation}
   Using (\ref{qshift.2}) and the fact that $\QQ(b,\mathbf s)_{1,2}=b_{1}b_{2}/(s_{2}-s_{1})$ we see that this inequality is the same as (\ref{34.bb}).

 \begin{example} \label{ex-3.2}{\rm 
 Here are some concrete  examples of the  relationship between  $\QQ(b,\mathbf s)$ and $W$ and $\QQ(b,\mathbf s')$ and $W'$. We take for $\QQ(b,\mathbf s)$ and $W$ the matrices in (\ref{harr2.211x}) and (\ref{harr2.21dd}). 
  In this case we have $b_{j}=s_{j}^{1/2}=r^{-(j-1)}$, $j\ge 1$,  and $W_{j,k}=r^{|k-j|}$, $j,k\ge 1$.  Therefore,    $b_{1}=s_{1}=1,  b_{2}=r^{-1}$, $ s_{2}=r^{-2}$ and, 
   \begin{equation}
         \QQ(b,\mathbf s')_{1,1}=-\frac{1+r^{2}\De}{(1+\De)(1-r^{2})}.
        \end{equation}
        Using the fact that $\De>-s_{1}$ and (\ref{3.63mm}) we see that we must have,
    \begin{equation}
-1<\De\leq   \frac{1  }{r},\quad\mbox{or equivalently,} \quad  \QQ(b,\mathbf s')_{1,1}\le - \frac{r}{1-r^{2}}.\label{34.bbb}
  \end{equation}

   \begin{itemize}
 \item [(i)] For $p>0$ set,
 \begin{equation}
   \De=-r^{p}\quad\mbox{or equivalently} \quad    \QQ(b,\mathbf s')_{1,1} =-\frac{1-r^{ p+2}}{(1-r^{p})(1-r^{2})} .
   \end{equation}  
Then by (\ref{wshift.1}),
\begin{equation}
   W'_{j,k} =r^{|k-j|}-r^{j+k+p-2}.\label{qshift.2x}
   \end{equation}
   \item [(ii)] For $ p\ge -1$ set,
 \begin{equation}
   \De= r^{p}\quad\mbox{or equivalently} \quad    \QQ(b,\mathbf s')_{1,1} =-\frac{1+r^{ p+2}}{(1+r^{p})(1-r^{2})}.
   \end{equation}  
Then by (\ref{wshift.1}),
\begin{equation}
   W'_{j,k} =r^{|k-j|}+r^{j+k+p-2}.\label{qshift.2xx}
   \end{equation}
 	\item [(iii)]  More generally for $\bb\geq r-r^2$     take,
  \begin{equation}
     \De=\frac{1-\bb-r^{2}}{\bb}\quad\mbox{or equivalently} \quad       \QQ(b,\mathbf s')_{1,1} =-\frac{\bb+r^{2}}{1-r^{2}}.
     \end{equation}
 Then by  (\ref{wshift.1}), 
\bea
W'_{j,k} &=& r^{|k-j|}+ \frac{1-\bb-r^{2}}{\bb}r^{j+k-2}\label{7.4maa}
\\&=&r^{|k-j|}-\frac{r^{j+k}}{\bb}+ \frac{1-\bb }{\bb}r^{j+k-2}\nn.
\eea

 \end{itemize}

\begin{remark} \label{rem-3.5}  {Let  $f=W'h$  for some  $h\in\ell_{1}^{+}$.  Let  $\wt Z'_{\al}=\{\wt Z'_{\al,j},j\in \mN\}$  be an $\al$-permanental sequence with kernel 
   $
 \wt W'=\{ \wt W'_{j,k};j,k\in\mN\}$, where
   \begin{equation}
 \wt   W' _{j,k}=W'_{j,k}+f_{k},\qquad j,k\in \mN. 
   \end{equation}
Then (\ref{18.4mkff})--(\ref{1.25mmx}) hold with $\wt Z _{\al}$ replaced by $\wt Z'_{\al}$,  and $W_{j,j}$ replaced by $W' _{j,j}$. 

This is easy to see. It follows from Theorem \ref{theo-4.2mm} itself that (\ref{18.4mkff})--(\ref{1.25mmx}) hold with $\wt Z _{\al}$ replaced by $\wt Z'_{\al}$ and $W_{j,j}$ replaced by $W '_{j,j}$ and $\mathbf s$ replaced by $\mathbf s'$.  Since $\lim_{j\to\ff}s_{j}=\ff$ implies that $s'_{j}\asymp s_{j}$,  $\KK_{\mathbf s}(j)\asymp \KK_{\mathbf s'}(j)$   we see that (\ref{18.4mkff})--(\ref{1.25mmx}) hold with $\wt Z _{\al}$ replaced by $\wt Z'_{\al}$,  and $W_{j,j}$ replaced by $W ' _{j,j}$ as stated. 
  }\end{remark}
  }\end{example}

  	   \section{Birth and death processes with emigration  related to first order Gaussian  autoregressive sequences}\label{sec-6}

 Let $g_{1}, g_{2}, \ldots $ be  a sequence of  independent identically distributed standard normal random variables and   $\{x_{n}\}$  
a   sequence of   positive numbers.    A first order autoregressive Gaussian sequence $\wh \xi=\{\wh \xi_{n}\}$ is  defined by,
 \begin{equation}
\wh\xi_{1}=g_{1},\qquad\wh\xi_{n}=x_{n-1}\wh\xi_{n-1}+g_{n},\qquad n\ge 2.\label{auto.1}
\end{equation}
It is easy to see that  
      \begin{equation}
     \wh\xi_{n}=\sum_{i=1}^{n}\(\prod_{l=i}^{n-1}x_{l}\)g_{i}\label{comp.10},
      \end{equation}
in which we take the  empty product $\prod_{l=n}^{n-1}x_{l}=1$.

\medskip	   We consider these processes with the added assumptions  that  $0< x_{n}< 1$, and $x_{n}\uparrow.$

 \medskip	Let  $\UU=\{\UU_{j,k};j,k\in\mN\}$ be  the covariance matrix for  $  \wh\xi  $.
 It follows that for $j\leq k$, 
      \begin{equation}
     \UU_{j,k}=\sum_{i=1}^{j}\(\prod_{l=i}^{j-1}x_{l}\prod_{l=i}^{k-1}x_{l}\)=\sum_{i=1}^{j}\(\prod_{l=i}^{j-1}x^{2}_{l}\prod_{l=j}^{k-1}x_{l}\)=\UU_{j,j}\prod_{l=j}^{k-1}x_{l}.\label{comp.10s}
      \end{equation}  
For $j>k$ we use the fact that $\UU$ is symmetric.

\medskip	  We now show that  $\UU$  can be written in the form of (\ref{1.16nn}).

\begin{lemma} \label{lem-7.1nn}
\begin{equation}
   \UU_{j,k} =\frac{s_{j}\wedge s_{k}}{b_{j}b_{k}},\qquad j,k\in\mN,\label{7.4mm}
   \end{equation}
where,
\begin{equation}
   b_{j}=\prod_{l=1} ^{j-1}x_{l}^{-1},\qquad \mbox{and}\qquad   s_{j} = \sum_{i=1}^{j}b^{2}_{i},\qquad \forall  j\in\mN.\label{arg.105a}
   \end{equation}
      Furthermore, $\{s_{j}\}$ is   a   strictly increasing  convex function of $j$.   (In particular, $\lim_{j\to\ff}s_{j}=\ff$.)  
 
\el

\Proof   By (\ref{comp.10s}) we have, 
\bea
     \UU_{j,j}& =&\sum_{i=1}^{j} \prod_{l=i}^{j-1}x^{2}_{l}   =\sum_{i=1}^{j}  \prod_{k=1}^{i-1}x^{ -2}_{k}  \prod_{l=1}^{j-1}x^{ 2}_{l}\label{pom.1o}\\
     &=&\frac{\sum_{i=1}^{j}b_{i}^{2}}{b_{j}^{2}} =\frac{s_{j}}{b_{j}^{2}}.\nn
   \eea
Using this and  (\ref{comp.10s}) again, we see that
for $j\le k$,  
\bea
 \UU_{j,k} &=& \UU _{j,j}\prod_{l=j} ^{k-1} x_{l}  = {\UU _{j,j} \over    \prod_{l=j} ^{k-1}x_{l}^{-1}} \\&=&  {\UU _{j,j}\prod_{l=1} ^{j-1}x_{l}^{-1}\over   \prod_{l=1} ^{k-1}x_{l}^{-1}} 
 = {\UU _{j,j}\prod_{l=1} ^{j-1}x_{l}^{-2}\over  \prod_{l=1} ^{j-1}x_{l}^{-1} \prod_{l=1} ^{k-1}x_{l}^{-1}}=\frac{\UU _{j,j}b^{2}_{j}}{b_{j}b_{k}}=\frac{s_{j}}{b_{j}b_{k}} \nn,
\eea
which is (\ref{7.4mm}).
 
  Since   $b_{j}> 1$,  for all   $j>1$, we see that   $s_{j+1}-s_{j}> 1$, for all $j\in\mN$, so that $s_{j}\uparrow\ff$. 
  
  We can say more than this. By   (\ref{arg.105a}),
  \begin{equation}
   (s_{j+1}-s_{j})- (s_{j}-s_{j-1})=b_{j+1}^{2}-b^{2}_{j}> 0,\label{5.8nn}
   \end{equation} 
   which shows that $\{s_{j}\}$ is   an increasing  convex function of $j$. \qed

 \begin{lemma}\label{lem-gconc}  Let $b=\{b_{j}\}$, $j\in\mN$, be as given in (\ref{arg.105a}).   Then,\begin{equation}
     b=\UU h,\qquad \mbox{for some  $h  \in \ell_{1}^{+}$}, 
    \end{equation}
  and $\UU$ is the potential density of a continuous symmetric transient Markov chain on $\mN$.

    Furthermore, the function  $g(s_{j})=b_{j}$, $j\in\mN$, $g(0)=0$, is an increasing   concave function of $\{s_{j}\}$   and  $s_{j}/b_{j}\uparrow\ff$. 
\end{lemma}  
 
  \Proof For $j\geq 2 $ we have,
  \begin{equation}
 {b _{j}-b _{j-1} \over s_{j}-s_{j-1}}= {b _{j}-b _{j-1} \over b^{2}_{j} }= 
 \frac{1}{b_{j}}\(1-\frac{b_{j-1}}{b_{j }}\)=\frac{ 1-x _{j-1} }{b_{j}},
 \label{pom.2}
  \end{equation}
  which is decreasing in $j$. In particular
  \begin{equation}
  {b _{2}-b _{1} \over s_{2}-s_{1}} =\frac{ 1-x _{1} }{b_{2}}=x _{1}(1-x _{1})<1\label{pom.2z}.
  \end{equation}
Therefore,  since $b _{1}=s _{1}=1$, 
 this  shows that    $g(s_{j})=b_{j} $ is    a concave function on  $\{   0\}\cup \{s_{j}, j\geq 1\}$. It follows from Corollary \ref{cor-2.1} that $b=\UU h$   for some  $h  \in \ell_{1}^{+}$.
  Using this and  Lemma \ref{lem-specmar} it follows  that $\UU$ is the potential density of a continuous symmetric transient Markov chain on $\mN$.
  
It is easy to see that 
  \begin{equation}
 \frac{ s_{j} }{b_{j}}\uparrow\ff. \label{pom.3}
  \end{equation}
We use the fact that $b_{j}\uparrow$, which implies that $\lim_{j\to\ff}b _{j} $   exists. If the limit is finite,  (\ref{pom.3})   is trivial because $s_{j}\to\ff$.

 If   $\lim_{j\to\ff}b _{j}=\ff$,   (\ref{pom.3})  follows because,  
 \begin{equation}
   \frac{s_{j}}{b_{j}}=\UU_{j,j}b_{j}> b_{j},
   \end{equation}
   since $\UU_{j,j}>1$, (see (\ref{5.12mm})).  \qed

 	 Although it is not obvious the next lemma shows that $\UU_{j,j}$ is strictly increasing.  
  
 	\begin{lemma} \label{lem-5.2mm}  The terms  $\,\UU_{j,j}$ are strictly increasing. Consequently,   $\lim_{j\to\ff}\UU_{j,j}$ exists. 
Furthermore, for all $j\ge 2$, 
\begin{equation}
 \UU_{j,j}\ge 1+x^{2}_{j-1}.\label{5.12mm}
   \end{equation} 
 \end{lemma}

\Proof    By  (\ref{comp.10s}),  
\be 
   \UU_{j,j} =\sum_{i=1}^{j}\prod_{l=i}^{j-1}x_{l}^{2}= \sum_{p=0}^{j-1}\prod_{l=j-p}^{j-1}x_{l}^{2}\label{4.13f}.
   \ee 
Similarly, 
\be
  \UU_{j+1,j+1}=  \sum_{p=0}^{j}\prod_{l=j+1-p}^{j }x_{l}^{2}\nn= \prod_{l=1}^{j }x_{l}^{2}+\sum_{p=0}^{j-1}\prod_{l=j+1-p}^{j }x_{l}^{2}.
   \ee 
Since $x_{l}\uparrow$,    
\begin{equation}
 \prod_{l=j+1-p}^{j }x_{l}^{2}\geq  \prod_{l=j-p}^{j-1}x_{l}^{2}.
   \end{equation}
   Therefore,  
   \begin{equation}
        \UU_{j+1,j+1}\ge   \UU_{j ,j }+\prod_{l=1}^{j }x_{l}^{2}.
      \end{equation}
This shows that   $\{\UU_{j,j}\}$ is strictly increasing.
 
 To obtain (\ref{5.12mm}) we simply note that for $j\ge 2$,   by the first equality in (\ref{4.13f})
\begin{equation}
   \UU_{j,j}\ge \prod_{l=j}^{j-1}x_{l}^{2}+\prod_{l=j-1}^{j-1}x_{l}^{2}= 1+x_{j-1}^{2} . 
   \end{equation}
\qed

\vspace{-.25 in}
 \begin{remark}\label{rem-5.1} {\rm In Lemma \ref{lem-gconc} we saw that  $g(s_{j})=b_{j}$ is an increasing   concave function of $\{s_{j}\}$   and  $s_{j}/b_{j}\uparrow\ff$.   Since $s_{j}/b^{2}_{j}=\UU_{j,j}$,  Lemma \ref{lem-5.2mm} strengthens  this to        $s_{j}/b^{2}_{j} \uparrow$. Although it is possible that $\lim_{j\to\ff}s_{j}/b^{2}_{j}<\ff$.
  }\end{remark}

	 \begin{lemma} \label{lem-5.3mm}
 \be
\lim_{j\to\ff}x_{j}=\de<1\quad\mbox{if and only if}\quad  \lim_{j\to\ff}\UU_{j,j}=\frac{1}{1-\de^{2}}.\label{5.88mm} 
\ee
and
  \be
\lim_{j\to\ff}x_{j}=1\quad\mbox{if and only if}\quad  \lim_{j\to\ff}\UU_{j,j}=\ff,\label{5.87mm}
\ee
 \end{lemma} 
 
\Proof Suppose   $\sup_{j}\UU_{j,j}<\ff$. Then by Lemma \ref{lem-5.2mm},   $\lim_{j\to\ff}\UU_{j,j}=c$, for some    $c>1$.  
 Note that by (\ref{auto.1}), 
 \be
 \UU_{j+1,j+1}=x_{j}^{2}\UU_{j ,j }+1.\label{6.13nn}
 \ee
 It   follows from this  that 
\be
\lim_{j\to\ff}x_{j} =\(\frac{c-1}{c}\)^{1/2}.
\ee
Setting this last expression equal to $\de$ shows that if  $ \lim_{j\to\ff}\UU_{j,j}= {1}/{1-\de^{2}}$, for $0<\de<1$,  then $\lim_{j\to\ff}x_{j}=\de<1.$

Now suppose that $\lim_{j\to\ff}x_{j}=\de<1.$ Then
\begin{equation}
   \UU_{j,j}\le \sum_{i=1}^{j}\de^{2(i-1)}\le \frac{1}{1-\de^{2}}.
   \end{equation}
This show that $\lim_{j\to\ff}   \UU_{j,j}=d$, for some $d<\ff$. Taking the limit as $j\to\ff$ in (\ref{6.13nn}) we see that
\begin{equation}
   d=\de^{2}d+1.
   \end{equation}
Solving for $d$ we see that $ \lim_{j\to\ff}\UU_{j,j}= {1}/{(1-\de^{2})}$.

 \medskip	 The statement in (\ref{5.87mm}) is implied by (\ref{5.88mm}). 
\qed

\medskip	 \noindent {\bf  Proof of Theorem \ref{theo-4.2mmf}}   When   $f= \UU  h$, where $   h  \in\ell_{1}^{+}$ it follows immediately from Theorem \ref{theo-4.2mm} that 
  if  $\limsup_{j\to\ff}s_{j}/s_{j-1}<\ff$, then,
 \begin{equation}
   \limsup_{j\to\ff}\frac{ \wt \XX_{\al,j}}{   \UU_{j,j} \log \log (\UU_{j,j}b_{j}^{ 2}) }= 1\quad a.s.,\label{5.44kk}
   \end{equation}
and if
$\liminf_{j\to\ff}s_{j}/s_{j-1}>1$, then,
\begin{equation}
\limsup_{j\to \ff} {\wt \XX_{\al,j} \over     \UU_{j,j}\log {j }}=1,  \qquad 
a.s.\label{1.25mmxo}, \qquad \forall\, \al>0.
\end{equation} 

It follows from (\ref{6.13nn}) and   (\ref{pom.1o}) that,
 \begin{equation}
   \frac{s_{j+1}}{s_{j}}=\frac{\UU_{j+1,j+1}}{\UU_{j,j}x_{j}^{2}}=1+\frac{1}{x^{2}_{j }\UU_{j ,j }}.
   \end{equation}
   We know from Lemma \ref{lem-5.2mm} that $ \lim_{j\to\ff}\UU_{j,j}  $ exists.
 If $ \lim_{j\to\ff}\UU_{j,j} = \ff$,   $\lim_{j\to\ff}s_{j+1}/s_{j}=1$, which gives (\ref{5.44kk}) and (\ref{5.44mm}). 
 
 If  $\limsup_{j\to\ff} \UU_{j,j}  =c<\ff$, then
 \be
 \liminf_{j\to\ff}s_{j+1}/s_{j}\ge 1+\frac{1}{c},
 \ee
 which gives (\ref{1.25mmxo}) and, by (\ref{5.88mm}), also (\ref{5.45mm}). 
 
  The proof of (\ref{5.44vv}) is given in Lemma \ref{lem-5.7nn} below.
  
\medskip	  That (\ref{5.45mm}) extends to $\al$-permanental processes $\wt  \XX_{\al} $ with kernels  $\wt \UU$   in which $f$  an excessive function for $ \XX$   with the property that $f\in c^{+}_{0}$,  follows from Theorem \ref{theo-1.8mm}. We show that the conditions in (\ref{1.39}) are satisfied.
  Since 
\begin{equation}
     \UU_{j,k}=\sum_{i=1}^{j}\(\prod_{l=i}^{j-1}x_{l}\prod_{l=i}^{k-1}x_{l}\),\qquad j\le k,\label{comp.10kkx}
     \end{equation}  
and $\{x_{j}\}$ is an increasing sequence, we have 
\begin{equation}
  \UU_{j,k}\leq \sum_{i=1}^{j}\de^{(j-i)+(k-i)}= \de^{j+k} \sum_{i=1}^{j}\de^{- 2 i }\le \de^{k-j} \frac{1}{1-\de^{2}},\qquad j\le k.\label{comp.10kky}
\end{equation}
Therefore,
\begin{equation}
\sum_{k=j}^{\ff}   \UU_{j,k}\leq \frac{1}{(1-\de)(1-\de^{2})}\label{rree}.
\end{equation}
Since this also holds when $k\le j$ we see that 
  $\|\UU\|<\ff$. In addition it follows from (\ref{comp.10kkx}) that   $\UU_{j,j}\ge 1.$   

The fact that $f\in c^{+}_{0}$ if and only if  $f= \UU  h$, where $   h  \in c_{0}^{+}$ follows from Lemma \ref{lem-8.2nn} once we show that (\ref{8.14}) holds.   
This is easy to see since,
\begin{equation}
 \sum_{j=1}^{k/2}   \UU_{j,k}\leq  \frac{\de^{k/2}}{ (1-\de^{2})}  \sum_{j=1}^{k/2} \de^{k/2-j}\leq \frac{\de^{k/2}}{(1-\de)(1-\de^{2})}.
\end{equation} 
 \qed

 We use the next lemma to obtain Example \ref{ex-1.1} $(i)$ and $(ii)$.

 \begin{lemma} \label{lem-5.4nn} If $\lim_{j\to\ff} j(1-x^{2}_{j})= c$,   for some $c\ge 0$, then
 \begin{equation}
    \UU_{j,j}\sim \frac{1}{1+c}\, j\quad\mbox{as} \quad j\to\ff,\label{5.31nn}
   \end{equation}
   and
   \begin{equation}
   \log \log (\UU_{j,j}b_{j}^{ 2}) \sim\log\log j\qquad\mbox{as $j\to\ff.$}\label{5.19mm}
   \end{equation}
 \end{lemma}

  \Proof 	We have,  
   \begin{equation}
         \UU_{j,j}=\sum_{i=1}^{j}\(\prod_{l=i}^{j-1}x^{2}_{l} \). \label{5.26nn}
   \end{equation}
   For some  $\ep>0$ let $0<a<c<b$ be such that $|a-b|<\ep$. Since $\lim_{l\to\ff} l(1-x^{2}_{l})= c$ we can find a $j_{0}$ such that,
\begin{equation}
- \frac{b}{l}< \log x_{l}^{2}<-\frac{a}{l},\qquad\forall\, l\ge j_{0}.\label{5.27nn}
   \end{equation} 
   Consequently, 
 \bea 
   \sum_{i=j_{0}}^{j}\prod_{l=i }^{j-1}x^{ 2}_{l}&=& \sum_{i=j_{0}}^{j}\exp\( \sum_{l= i}^{j-1} \log x_{l}^{2}\)\le \sum_{i=j_{0}}^{j}\exp\(-\sum_{l=i}^{j-1} \frac{a}{l}\)\label{5.28nn}\\
   &\sim&   \frac{1}{(j-1)^{a}} \sum_{i=j_{0}}^{j}    i^{a}              \sim \frac{j }{1+a},\qquad\mbox{as $j\to\ff.$} \nn
\eea 
   Similarly,
   \begin{equation}
    \sum_{i=j_{0}}^{j}\prod_{l=i }^{j-1}x^{ 2}_{l}\ge \frac{j }{1+b},\qquad\mbox{as $j\to\ff.$} \label{5.29nn}
   \end{equation} 
 Note that the left-hand sides of   (\ref{5.28nn})   and (\ref{5.29nn}) differ from (\ref{5.26nn}) by some finite number. Therefore, since (\ref{5.27nn})   holds  for all $\ep>0$, we get (\ref{5.31nn}) when  $c>0$.

  When $c=0$ the left-hand side of (\ref{5.27nn})  holds for all $b>0$. Therefore,  
     \begin{equation}
   \liminf_{j\to\ff}\frac{  \UU_{j,j} }{j}\ge 1.
   \end{equation}
    However, since $\UU_{j,j} \le j$, we get (\ref{5.31nn}) when $c=0$. 
 
  To get (\ref{5.19mm}) we   note that by (\ref{5.31nn}),
  \begin{equation} 
  \log \log (\UU_{j,j}b_{j}^{ 2}) \sim \log\(\log \(\frac{j}{1+c}\)+ \log b^{2}_{j}\) ,\qquad\mbox{as $j\to\ff$}\label{5.44mc}. 
   \end{equation}
   Furthermore, by   (\ref{5.31nn}) 
 \bea  
\log  b_{j}^{2}&=&\log \prod_{l=1 }^{j-1}x^{ -2}_{l}= - \sum_{l= 1}^{j-1} \log x_{l}^{2}  \label{5.28nh}\\
    & =& C  - \sum_{l= j_{0}}^{j-1} \log x_{l}^{2}  \le C+   \sum_{l= j_{0}}^{j-1} \frac{b}{l}  \le C+ b\log  {j},   \nn
   \eea 
   where $C= - \sum_{l= 1}^{j_{0}-1} \log x_{l}^{2}$. Using (\ref{5.44mc}) and (\ref{5.28nh}) and the fact that $\log b_{j}\ge0$, we get  (\ref{5.19mm}).\qed

   \noindent{\bf Proof of Example \ref{ex-1.1} $(i) $ }   This follows immediately from Lemma \ref{lem-5.4nn}. We now show that this includes the case where   $\prod_{j=1}^{\ff}x _{j} >0 $. Set $x_{j}=1-\ep_{j}$. Therefore,  for some $C>0$,
   \begin{equation}
   \prod_{j=1}^{\ff}x _{j}\leq C\exp\(-\sum_{j=1}^{\ff}\ep_{j}\).
   \end{equation} 
   If $\prod_{j=1}^{\ff}x _{j} >0 $ then   $\sum_{j=1}^{\ff}\ep_{j}<\ff$.  
Since $\ep_{j}$ is decreasing it follows from  Lemma \ref{lem-spec} that  $\ep_{j}=o(1/j)$, as $j\to\ff.$ Therefore the condition that $j(1-x_{j} ^{2})\to 0$, as $j\to\ff$ includes the case when $\prod_{j=1}^{\ff}x _{j} >0 $ . \qed

   \begin{lemma} \label{lem-spec} If  $\ep_{j}\downarrow$ and $\sum_{j=1}^{\ff}\ep_{j}<\ff$, then 
   $\ep_{j}=o(1/j)$, as $j\to\ff.$
 \end{lemma}
 
 \Proof Suppose that $\limsup_{j\to\ff}j\ep_{j}\ge \de>0$. Then we can find a subsequence $\{j_{k}\}$ such that
 \begin{equation}
 \frac{j_{k}}{j_{k+1}}\le \frac{1}{2},\qquad\mbox{and}\qquad  \ep_{j_{k}}\ge \frac{\de}{j_{k}}.
   \end{equation}
  Therefore,
  \bea
   \sum_{j=1}^{\ff}\ep_{j}&=&\sum_{k=1}^{\ff}\sum_{l=j_{k}}^{j_{k+1}-1}\ep_{j}\ge \sum_{k=1}^{\ff} \ep_{j_{k+1}}(j_{k+1}-j_{k })\\
   &\ge&\de\sum_{k=1}^{\ff}  \(1-\frac{j_{k}}{j_{k+1} }\)=\ff.\nn
      \eea
  \qed
  
   \noindent{\bf Proof of Example \ref{ex-1.1} $(ii)$  }   This follows immediately from Lemma \ref{lem-5.4nn}. \qed
  
    The next lemma gives some useful information about $\UU_{j,j}$:
	\begin{lemma}\label{lem-5.6mm}  	\begin{equation}
 \UU_{j,j} \le \frac{1}{1-x_{j}^{2}}. \label{5.22k}
   \end{equation}
	Furthermore the following are equivalent:
 \begin{equation}
   \UU_{j,j}\sim  \frac{1}{1-x_{j}^{2}},\quad\mbox{as} \quad j\to\ff,\label{5.30mmq}
   \end{equation}
 and
  \begin{equation}
 \UU_{j+1,j+1}-\UU_{j ,j}\to 0,\quad\mbox{as $j\to\ff$}.  \label{5.30ws}
   \end{equation}
   \el
    
 \Proof  By (\ref{6.13nn}) and  Lemma \ref{lem-5.2mm},
 \begin{equation}
 (1-  x_{j}^{2})\UU_{j,j}=1-(\UU_{j+1,j+1}-\UU_{j,j})\le 1.\label{5.42nn}
   \end{equation}
All the statements in this lemma follow easily from this.\qed

Using Lemma \ref{lem-5.6mm} we make (\ref{5.44mm}) more specific:

 \begin{lemma}\label{lem-5.7nn} In Theorem \ref{theo-4.2mmf} assume  that $\UU_{j,j}$ is a regularly varying function with index $0<\bb<1$  then,
 \begin{equation}
   \limsup_{j\to\ff}\frac{ \wt \XX_{\al,j}}{   \UU_{j,j}   \log {j}} = 1-\bb,\qquad a.s.\,\qquad\forall\,\al>0\label{5.40nn}.
   \end{equation}
\end{lemma}

 	 \Proof   Considering  (\ref{5.44mm}) we need to show that
   \begin{equation}
   \log \log (\UU_{j,j}b_{j}^{ 2}) =   \log (\log  \UU_{j,j}+\log b_{j}^{ 2})\sim(1-\bb)\log j,\qquad\mbox{as $j\to\ff.$}\label{5.19kk}
   \end{equation}
   Furthermore, since $\UU_{j,j}\le j$, we need to show that, 
   \begin{equation}
   \log\log b_{^{j} }^{2}\sim(1-\bb)\log j,\qquad\mbox{as $j\to\ff.$}\label{5.40kk}
   \end{equation}

   To be specific let $\UU_{j,j}=g(j)=j^{\bb}L(j)$, where $L$ is a slowly varying function.  Then clearly, (\ref{5.30ws}) holds. Therefore, by (\ref{5.30mmq}), $(1-x_{j}^{2})\sim 1/g(j)$.
  It follows that for all $\ep>0$ and $0<a<1<b$ such that $|a-b|\leq \ep$ we can find an integer  $j_{0}$ such that
\begin{equation}
- \frac{b}{g(l)}< \log x_{l}^{2}<-\frac{a}{g(l)},\qquad\forall\, l\ge j_{0}.\label{5.27nn5}
   \end{equation} 
 Similar to (\ref{5.28nh}),  for all $j$ sufficiently large,
  \bea
 \log     b_{j}^{2}&=& - \sum_{l= 1}^{j-1} \log x_{l}^{2} \\
 &\leq  & - \sum_{l= 1}^{j_{0}-1} \log x_{l}^{2}+\sum_{l=j_{0}}^{j-1} \frac{b}{g(l)}  \sim        \frac{b j^{1-\bb}}{(1-\bb)L(j)}.  \nn    
            \eea  
   Likewise,       
    \be 
 \log     b_{j}^{2}\ge   - \sum_{l= 1}^{j_{0}-1} \log x_{l}^{2}+\sum_{l=j_{0}}^{j-1} \frac{a}{g(l)}  \sim        \frac{a j^{1-\bb}}{(1-\bb)L(j)}.  \nn            \ee           
  These two inequalities give (\ref{5.40kk}).\qed           
        
  In the proof of Lemma \ref{lem-5.7nn} we use the fact that when $\UU_{j,j}$ is a regularly varying function with index $0<\bb<1$, then $(1-x_{j}^{2})\sim 1/\UU_{j,j}$. What we do not show is that when $(1-x_{j}^{2})\sim h(j)$, for some regularly varying function $h(j)$ with index $-1<\bb'<0$ then $\UU_{j,j}\sim 1/h(j)$. We only consider this in the special case given in Example \ref{ex-1.1} (iii).

 \medskip	\noindent{\bf Proof of Example \ref{ex-1.1} (iii)  }
 We have,  
   \begin{equation}
         \UU_{j,j}=\sum_{i=1}^{j}\(\prod_{l=i}^{j-1}x^{2}_{l} \). \label{5.26nnj}
   \end{equation}
   For some  $\ep>0$ let $0<a<1<b$ be such that $|a-b|<\ep$. Since $\lim_{l\to\ff} l^{\bb}(1-x^{2}_{l})= 1$ we can find a $j_{0}$ such that,
\begin{equation}
- \frac{b}{l^{\bb}}< \log x_{l}^{2}<-\frac{a}{l^{\bb}},\qquad\forall\, l\ge j_{0}.\label{5.27nnj}
   \end{equation} 
   Consequently, 
 \bea 
   \sum_{i=j_{0}}^{j}\prod_{l=i }^{j-1}x^{ 2}_{l}&=& \sum_{i=j_{0}}^{j}\exp\( \sum_{l= i}^{j-1} \log x_{l}^{2}\)\le \sum_{i=j_{0}}^{j}\exp\(-\sum_{l=i}^{j-1} \frac{a}{l^{\bb}}\)\label{5.28nnj}\\
   &\sim&  \exp\(-\frac{aj^{1-\bb}}{1-\bb}\)           \int_{j_{0}} ^{j}\exp\(\frac{ax^{1-\bb}}{1-\bb}\)\,dx, \qquad\mbox{as $j\to\ff$, }
 \nn
\eea 
and   
\bea
      \int_{j_{0}} ^{j}\exp\(\frac{ax^{1-\bb}}{1-\bb}\)\,dx &=&\frac{1}{a} \int_{j_{0}}^{j}x^{\bb}d\(\exp\(\frac{ax^{1-\bb}}{1-\bb}\)\)\\
 &\sim&\nn \frac{j^{\bb}}a \exp\(\frac{aj^{1-\bb}}{1-\bb}\), \qquad\mbox{as $j\to\ff$, }
   \eea
where, for the last line we use integration by parts.     Therefore,
\begin{equation}
    \sum_{i=j_{0}}^{j}\prod_{l=i }^{j-1}x^{ 2}_{l}\le \frac{j^{\bb}}a,\qquad\mbox{as $j\to\ff$}.\label{5.46}
   \end{equation}
A similar argument shows that the left-hand side of (\ref{5.46}) is greater than or equal to $j^{\bb}/b$ as $j\to\ff$. Using these observations and following the proof of Lemma \ref{lem-5.4nn} we   see  that
\begin{equation}
   \UU_{j,j}\sim j^{\bb}, \qquad\mbox{as $j\to\ff$.}\label{5.30nn}
   \end{equation}
 Therefore (\ref{5.44mh}) follows from (\ref{5.44vv}).\qed

We now explicitly describe the $Q$ matrix corresponding to $\UU$ in Lemma \ref{lem-7.1nn}, which is 
\be
Q(b,\mathbf s ):= D_{b}  Q(\mathbf s)D_{b}.
\ee
  (See (\ref{harr2.211q}) and (\ref{3.14mm}).)

It follows from   (\ref{arg.105a}) that,   
\begin{equation}
  a_{j} ={1 \over s_{j}-s_{j-1}}=b^{-2}_{j},\qquad j\geq 2, \label{qmer.1}
   \end{equation}
  and 
     \begin{equation}
  a_{1} ={1 \over s_{1} } =b^{-2}_{1}=1. \label{qmer.2}
   \end{equation}
   Therefore (\ref{qmer.1}) holds for all $j\geq 1$. Consequently, for all $j\ge 1$,   
   \begin{equation}
  -Q(b,\mathbf s )_{j,j+1}=-b_{j} a_{j+1} b_{j+1}=-\frac{b_{j }}{b_{j+1}}=-x_{j},\label{qmer.3}
   \end{equation}
   and
      \begin{equation}
  -Q(b,\mathbf s )_{j,j}=b_{j}\( a_{j}+ a_{j+1}\) b_{j}= b_{j}^{2}\(\frac1{b_{j}^{2}}+\frac1{b_{j+1}^{2}}\)=1+x^{2}_{j}.\label{qmer.4}
   \end{equation}
Since $Q(b,\mathbf s )_{j,j+1}=Q(b,\mathbf s )_{j+1,j }$ we have,
       \be -Q(b,\mathbf s )=\left (
\begin{array}{ cccccc cc} 1 +x^{2}_{1}&-x_{1}&0&0&\dots&0&0&\dots  \\
-x_{1}&1+x^{2}_{2}&-x_{2}&0&\dots &0 &0 &\dots\\
0&-x_{2}&1+x^{2}_{3}&-x_{3}&\dots &0&0 &\dots\\
\vdots&\vdots&\vdots&\vdots&\ddots&\vdots &\vdots &\ddots\\
  0&0&0&0&\dots &1+x^{2}_{m}&-x_{m}&\dots  \\
  0&0&0&0&\dots &-x_{m}&1+x^{2}_{m+1}&\dots \\
\vdots&\vdots&\vdots&\vdots&\ddots&\vdots &\vdots &\ddots\end{array}\right )\label{qmer.5}
      \ee

 	 \begin{example} \label{ex-4.1}{\rm Let $x_{j}=r$. Then $b_{j}=r^{-(j-1)}$ 
and
\be -Q(\mathbf{s,b}) = \left (
\begin{array}{ cccccc cc}  
 1+r^{2} &-r&0&0&\dots   \\
-r&1+r^{2}&-r&0&\dots   \\
0&-r&1+r^{2}&-r &\dots  \\\vdots&\vdots&\vdots&\vdots&\ddots   \end{array}\right ),\label{har}
      \ee
 In addition
 \begin{equation}
   \UU_{j,j}=1+r^{2}+r^{4}+\cdots +r^{2(j-1)}=\frac{1-r^{2j}}{1-r^{2}}
   \end{equation}
   and for $j\le k$,
 \begin{equation}
   \UU_{j,k}=\UU_{j,j}r^{k-j}=\frac{r^{k-j}-r^{k+j}}{1-r^{2}}.\label{5.67}
   \end{equation}

   Consequently,
    \begin{equation}
   \UU_{j,k}=\UU_{j,j}r^{k-j}=\frac{r^{|k-j|}-r^{k+j}}{1-r^{2}},\qquad\forall j,k\in\mN.\label{5.63n}
   \end{equation} 
 }\end{example}
   Compare (\ref{qshift.2}) with $p=2$. 
(Note that 
\begin{equation}
   {e^{-\sqrt{2\de} |x-y|}-e^{-\sqrt{2\de} x}e^{-\sqrt{2\de} y} \over \sqrt{2\de}},\label{killed}
   \end{equation}
the $\de$-potential density for Brownian motion killed the first time it hits $0$.)

\medskip	 We show in Lemma \ref{lem-gconc} that the covariance of  the  first order Gaussian autoregressive sequence $\wh \xi$  in (\ref{auto.1}) is the potential  of a continuous time Markov chain $\UU=\{U_{j,k};j,k\in\mN\}$ where,
\be
  \UU_{j,k} =\frac{ s_{j}\wedge   s_{k}}{  b_{j} b_{k}}. 
  \ee
    At the end of Section \ref{sec-gen} we consider  the effect of a shift  $s_{j}\to s'_{j}=s_{j}+\De$ on such potentials. We now show that when apply such a shift to $\UU$ we still have   the covariance of a first order Gaussian autoregressive sequence.

 As in (\ref{auto.1}), let $g_{1}, g_{2}, \ldots $ be  a sequence of  independent identically distributed standard normal random variables,   $ 0<x_{n}<1, $   $x_{n}\uparrow$,  and take   $\wt\de\ne 0 $.     Consider the  Gaussian sequences $\wt \xi=\{\wt \xi_{n}\}$    defined by, \label{page-52}
 \begin{equation}
\wt \xi_{1}=\wt\de g_{1},\qquad\wt\xi_{n}=x_{n-1}\wt\xi_{n-1}+g_{n},\qquad n\ge 2.\label{auto.1qa}
 \end{equation}
This generalizes $\wh\xi$ in (\ref{auto.1}).
 
  \begin{theorem} \label{theo-4.2g} Let   $\{s_{j}\}$ and $\{b_{j}\}$   be as given  in  Lemma \ref{lem-7.1nn}.
   The covariance of the first order Gaussian auto regressive sequence $\wt \xi$ is
   $\UU'=\{U'_{j,k};j,k\in\mN\}$ where,
\be 
  \UU'_{j,k} =\frac{ s'_{j}\wedge   s'_{k}}{  b_{j} b_{k}} ,\label{34.7}
 \ee

 \begin{equation}
    s'_{j}=s_{j}+\De,\qquad \mbox{and}\qquad \De=\wt \de^{ 2}-1,\quad 0<\wt\de<\ff.
    \end{equation} 
    
    Furthermore,  $ \UU'$ is the potential density of a transient Markov chain if
 \begin{equation}
 0<\wt \de^{ 2}\leq \frac{1}{x_{1}(1-x_{1})}.\label{34.7d}
 \end{equation}
 \end{theorem}
 
  \Proof	It is easy to see that  
      \begin{equation}
     \wt\xi_{n}=\sum_{i=1}^{n}\(\prod_{l=i}^{n-1}y_{l}\)g_{i}\label{comp.10o},\qquad n\ge 2,
      \end{equation}
where $y_{1}=\wt   \de x_{1}  $, $y_{l}=x_{l}$, $l\ge 2$ and in which we take the  empty product $\prod_{l=n}^{n-1}y_{l}=1$. Therefore,
 for $j\leq k$,  
      \bea
     \UU'_{j,k}&=&\sum_{i=1}^{j}\(\prod_{l=i}^{j-1}y_{l}\prod_{l=i}^{k-1}y_{l}\)=\sum_{i=1}^{j}\(\prod_{l=i}^{j-1}y^{2}_{l}\prod_{l=j}^{k-1}y_{l}\)= \UU'_{j,j}\prod_{l=j}^{k-1}y_{l} \label{comp.10sxx}\\
	     &=&\nn \frac{ \UU'_{j,j}}{\prod_{l=j}^{k-1}y^{-1}_{l}}=\frac{ \UU'_{j,j}\prod_{l=1}^{j-1}y^{-1}_{l}}{\prod_{l=1}^{k-1}y^{-1}_{l}}=\frac{ \UU'_{j,j}b_{j}}{b_{k}}.
      \eea 
      Set $   s'_{j} =b_{j}^{2}  \UU'_{j,j}$, $j\in\mN$, so that (\ref{34.7}) holds.
Also note that,
\bea 
   \UU'_{j,j}&=&\sum_{i=1}^{j}\prod_{l=i}^{j-1}y^{2}_{l}=\prod_{l=1}^{j-1}y^{2}_{l}+ \sum_{i=2}^{j}\prod_{l=i}^{j-1}y^{2}_{l}\\
  &=&\wt   \de^{ 2}\prod_{l=1}^{j-1}x^{2}_{l}+ \sum_{i=2}^{j}\prod_{l=i}^{j-1}x^{2}_{l}.\nn 
 \eea
Therefore, since $b_{1}=1$,
\begin{equation}
  s'_{j} =b_{j}^{2}   \UU'_{j,j}=\wt   \de^{ 2}+ \sum_{i=2}^{j}\prod_{l=1}^{i-1}x^{-2}_{l}= \wt   \de^{ 2}+ \sum_{i=2}^{j}b_{i}^{2}= (\wt   \de^{ 2}-1)+ s_{j},\label{4.93}
   \end{equation}
   where we use (\ref{arg.105a}) for the last equation.
This gives   $s'_{j}=s_{j}+\De$ with $\De=\wt   \de^{ 2}-1$. 

  It follows from Lemma \ref{lem-specmarx} and (\ref{arg.105a}) that when    (\ref{34.7d}) holds, $\UU'$ is the potential density of a Markov chain.
\qed

  \begin{remark} \label{rem-4.2}{\rm Assume condition (\ref{34.7d}), so that $ \UU' $ is the potential density of a transient Markov chain which we denote by $\XX'$.  Let $f $ be a finite  excessive function for $\XX'$.  Let  $\wt \XX'_{\al}=\{\wt \XX'_{\al,j},j\in \mN\}$ be   an $\al$-permanental sequence
 with kernel        
      $
  \wt    \UU'=\{ \wt    \UU'_{j,k};j,k\in\mN\}$, where
   \begin{equation}
 \wt    \UU'_{j,k}=\UU'_{j,k}+f_{k},\qquad  j,k\in\mN. \label{4.38nn}
   \end{equation}
   Then using the same argument used in Remark \ref{rem-3.5} we see  that   if $f=\UU'h$ for some $h\in\ell_{1}$ then (\ref{5.44mm}) and (\ref{5.44vv}) hold with $\wt X _{\al}$ replaced by $\wt X'_{\al}$. Item  $(ii)$ in Theorem \ref{theo-4.2mmf} also  holds with $\wt X _{\al}$ replaced by $\wt X'_{\al}$. 
   
\medskip	Example \ref{ex-1.1} also holds with $\wt X _{\al}$ replaced by $\wt X'_{\al}$ since the computations depend on the relationship between $\UU$ and $\{b_{j}\}$ and $\{b_{j}\}$ is unchanged. 
   }\end{remark}

\begin{remark}\label{rem-5.2} {\rm Condition   (\ref{34.7d}) is necessary for $\UU'$ to be the potential density  of a Markov chain whereas (\ref{34.7}) holds for all $\wt\de\neq 0$.   This gives    examples of a critical point at which a covariance matrix ceases to be an inverse $M$-matrix. This has interesting implications in the study of Gaussian sequences with infinitely divisible squares.
 }\end{remark}

    \section{Markov chains with potentials that are the covariances of  higher order Gaussian  autoregressive sequences}	  \label{sec-7}

Consider a class of  $k$-th  order autoregressive Gaussian sequences,  for $k\ge 2$.  
Let  $g_{1}, g_{2},\ldots$ independent standard normal random variables and let  $ p_{i}>0$, $i=1,\ldots,k$, with  $\sum_{l=1}^{k}p_{l}\le 1$.    We define the Gaussian sequence  $\wt \xi=\{\wt \xi_{n},n\in\mN \}$  by,
\begin{equation}
\wt \xi_{1}=g_{1}, \hspace{.2 in}\mbox{and}\quad\wt \xi_{n}=\sum_{l=1}^{k}p_{l}\wt\xi_{n-l}+g_{n},\hspace{.2 in}n\geq 2,\label{ark.1}
\end{equation}
where $\wt \xi_{i}=0$ for  all $i\le 0$. Let $\VV=\{\VV_{j,k};j,k\in \mN\}$ denote the covariance of $\wt\xi$.

Our goal is to prove Theorem \ref{theo-4.2mmfx}. We begin by exhibiting some simple properties $\VV$. Set 
 \begin{equation}
   \phi_{n}=E(\wt \xi_{n}g_{1}),\label{7.2mm}
   \end{equation}
 and note that
  \begin{equation}
 \phi_{1}=1\hspace{.2 in}\mbox{ and  }\hspace{.2 in} \phi_{n}=\sum_{l=1}^{k}p_{l}\phi_{n-l},\hspace{.2 in}n\geq 2,\label{gf.1j}
 \end{equation}
 where $\phi_{n}=0$ for all $n\le 0$. Since  $\phi_{2}=p_{1}<1$ and $\sum_{l=1}^{k}p_{l}\le 1$, we see that,
 \begin{equation}
  \phi_{n}<1, \hspace{.2 in}\forall n\ge 2.\label{gf.1jj}
 \end{equation}
 
 \medskip	We now write $\{\wt \xi_{n}\}$ as a series with terms that are independent Gaussian random variables.
 
 	 \begin{lemma}  \label{lem-8.1mm}   
\begin{equation}
\wt\xi_{n}=\sum_{j=1}^{n} \phi_{n+1-j}g_{j}=\sum_{j }  \phi_{n+1-j}g_{j},\qquad j\in\mN,\label{ark.2}
 \end{equation}
 (since the terms in the last sum are all equal to 0 when $j\notin[1,n]$). Therefore,  
\begin{equation} 
\VV_{m,n}=E(\wt\xi_{m}\wt\xi_{n})=\sum_{j=1}^{m\wedge n}\phi_{m+1-j}\phi_{n+1-j}=\sum_{j=0}^{(m\wedge n)-1}\phi_{m -j}\phi_{n -j},\label{7.27mm}
   \end{equation}
   which implies, in particular,  that  
   \be
\VV_{m,n}\leq m\wedge n, \quad \VV_{1,1}=1,\quad\mbox{and}\quad \VV_{n,n}=  E(\wt\xi_{n}^{\,2})\uparrow.\label{6.6new}
   \ee
  \end{lemma}
 
 \Proof  We give a proof by induction.   Clearly (\ref{ark.2}) is true for $n=1$.  Then using (\ref{ark.1}) and induction we have
 \begin{eqnarray}
 \wt\xi_{n}&=&\sum_{l=1}^{k}p_{l}\wt\xi_{n-l}+g_{n}
 \label{ark.3}\\
 &=& \sum_{l=1}^{k}p_{l}   \sum_{j\leq n-1 }  \phi_{n-l+1-j}g_{j}  +g_{n}, \nonumber
 \end{eqnarray}
 where in the second equality we change nothing by allowing $j\leq n-1$ rather than $j\leq n-l$, 
 since $\phi_{n}=0$ for $n<1$. Interchanging  the order of summation this is equal to
 \begin{eqnarray}
 && \sum_{j\leq n-1 } \(  \sum_{l=1}^{k}p_{l}    \phi_{n-l+1-j} \) g_{j}  +g_{n}
 \\
 &&=  \sum_{j\leq n-1 }  \phi_{n+1-j}g_{j}  +g_{n}\nn
 \end{eqnarray}
where the last equality came from  (\ref{gf.1j}), since for    $j\leq n-1$
  we have $n+1-j\geq 2$. This 
 gives (\ref{ark.2}).
 
 The statement in (\ref{7.27mm}) follows immediately from (\ref{ark.2}); (\ref{6.6new}) is an immediate consequence of (\ref{7.27mm}),    (\ref{gf.1jj}) and (\ref{gf.1j}), since $\VV_{1,1}=\phi^{2}_{1}$.\qed
 
 	  We now introduce the matrix $A$ which,  with the additional condition that its   off diagonal elements are less that or equal to 0,      is the negative of  the $Q$ matrix for the continuous time   symmetric Markov chains on $\mN$    with potential  densities  $\VV=\{\VV_{j,k}$, $j,k\in\mN\}$.

\begin{lemma}\label{lem-7.5mm}    Let $A=\{A_{m,n};m,n\in\mN\}$ where,
  \begin{equation}
  A_{m,m}= 1+\sum_{i=1}^{  k}   p^{2}_{i},\qquad \forall\, m\in \mN\label{7.46mm},
  \end{equation}
 \be 
 A_{m,n} =-p_{|m-n|}+\sum_{l=1}^{k-|m-n|}   p_{l}\,\, p_{|m-n|+l},\qquad \mbox{for all $    1\le |m-n|\le k$}\label{7.47mm},
  \ee
  and 
  \be
  A_{m,n} =0,\qquad\forall\,   |m-n|>k.\label{7.47nn}
  \ee
   Then   
    \begin{equation}
    \VV A=A \VV=I\label{wa},
    \end{equation} in the sense of   multiplication of matrices. That is, for each $i,l\in\mN$,
   \begin{equation}
 \sum_{j} \VV_{i,j}A_{j,l}=\de_{i,l},  \label{matrixj.1}
   \end{equation}
and similarly $A\VV$.
 \end{lemma}

   Clearly $A_{m,n}$ depends only on $|m-n|$. Set
 \begin{equation}
 a_{|m-n|}=A_{m,n},\qquad  n,m\in\mN.  \label{ark.100}
 \end{equation}
Note  that $A$ is a symmetric T\"oeplitz matrix and that for $j\ge k+1$,  the $j$-th row of $A$ has the form
\begin{equation}
   0,\ldots, 0,a_{k},\ldots, a_{1},a_{0},a_{1},\ldots,a_{k},0,0,\ldots, \label{6.48mmo}
   \end{equation}
 where the initial sequence of zeros has $j-k$ terms.
 
 For $j\le k$ 
 the $(j+1)$-st row of $A$ has the form
\begin{equation}
 a_{j},\ldots, a_{1},a_{0},a_{1},\ldots,a_{k},0,0,\ldots .\label{6.49mm}
   \end{equation}

 Here is an explicit example.

\begin{example} \label{ex-6.2}{\rm   When $k=2$,
\begin{equation}
  \hspace{-4.8in} A=
   \end{equation}
   \vspace{-.25in}
\be   \left (\!\!
\begin{array}{ cccccc cc}  
 1+p_{1}^{2}+p^{2}_{2} &-p_{1}+p_{1}p_{2}&-p_{2}&0&0&0&\dots   \\
-p_{1}+p_{1}p_{2}& 1+p_{1}^{2}+p^{2}_{2} &-p_{1}+p_{1}p_{2}&-p_{2}&0&0&\dots   \\
-p_{2}&-p_{1}+p_{1}p_{2}& 1+p_{1}^{2}+p^{2}_{2}&-p_{1}+p_{1}p_{2}&-p_{2}&0&\dots  \\
0&-p_{2}&-p_{1}+p_{1}p_{2}& 1+p_{1}^{2}+p^{2}_{2}&-p_{1}+p_{1}p_{2}&-p_{2}& \dots  \\\vdots&\vdots&\vdots&\vdots&\vdots&\vdots&\ddots   \end{array}\!\!\right ) \label{hardd}
      \ee
We see that in this case $A$ is a symmetric T\"oeplitz matrix with five non-zero diagonals. The row sums for all rows after the second row  are equal to  $(1-(p_{1}+p_{2}))^{2}$. Note also that    $-A$ is a Q-matrix since, $p_{1}p_{2}\leq p_{1}$.  
}\end{example}
 
 \noindent{\bf Proof of Lemma \ref{lem-7.5mm} } We introduce two   infinite matrices,
     \be L =\left (
\begin{array}{ cccccc ccc} 1 &0&0 &\dots&0&0&\vdots &\vdots&\dots  \\
-p_{1}&1 &0& \dots &0 &0 &\vdots &\vdots&\dots\\
-p_{2}&-p_{1}&1  &\dots &0&0 &\vdots &\vdots&\dots\\
\vdots&\vdots&\vdots&\ddots&\vdots&\vdots &\vdots &\vdots&\dots\\
 -p_{k}&-p_{k-1}&-p_{k-2} &\dots &-p_{1}&1&0&0&\dots  \\
  0& -p_{k}&-p_{k-1}&\dots&-p_{2} &-p_{1}&1 &0&\dots \\
\vdots&\vdots&\vdots&\ddots&\vdots&\vdots &\vdots &\vdots&\ddots\end{array}\right ),\label{ark.4}
      \ee 
      and  
         \be \Phi =\left (
\begin{array}{ cccccc ccc} \phi_{1} &0&0 &\dots&0&0&0&0&\dots  \\
 \phi_{2}& \phi_{1} &0& \dots &0 &0 &0&0&\dots\\
 \phi_{3}& \phi_{2}& \phi_{1}&\dots &0&0 &0&0&\dots\\
\vdots&\vdots&\vdots&\ddots&\vdots&\vdots &\vdots &\vdots&\ddots\\
 \phi_{n}& \phi_{n-1}& \phi_{n-2}&\dots & \phi_{1}&0&0&0&\dots  \\
\vdots&\vdots&\vdots&\ddots&\vdots&\vdots &\vdots &\vdots&\ddots\end{array}\right ).\label{ark.5}
      \ee 
 where $\{\phi_{n}\}$ is given in (\ref{gf.1j}).

   It is  easy to see that, 
 \begin{equation}
L\Phi =\Phi L=I,\qquad\mbox{and}\qquad L^{T}\Phi^{T} =\Phi^{T} L^{T}=I.\label{ark.6}
 \end{equation}
  We also give  is an analytical  proof.
Set $p_{0}=-1$ and $p_{j}=0$, $j<0$, and write,
 \begin{equation}
 L_{i,j}=-p_{i-j},\qquad i,j\in \mN. \label{ark.6p}
 \end{equation}
and 
 \begin{equation}
 \Phi _{i,j}=\phi_{i+1-j},\qquad i,j\in \mN.  \label{ark.7}
 \end{equation}
 Consequently, \bea
  \(L\Phi \)_{m,n}&=&-\sum_{j}  p_{m -j} \phi_{j+1-n}\label{ark.11} \\
&=&\phi_{m+1-n}   -\sum_{n\leq j<m}  p_{m -j} \phi_{j+1-n}\nn.
 \eea
 When $n=m$ there are no non-zero terms in the final sum in (\ref{ark.11})  and since $ \phi_{1}=1$ we have   $\(L\Phi \)_{n,n}=1$. If $m< n$, all the terms in the last line of  (\ref{ark.11}) are equal to 0, so we have  $\(L\Phi \)_{m,n}=0$. When $m> n$, we set $l=m-j$ and write (\ref{ark.11}) as,
  \bea
 &&\(L\Phi \)_{m,n}=\phi_{m+1-n}   -\sum_{l=1}^{m-n}  p_{l} \phi_{m+1-n-l}=0,\label{ark.12}
 \eea
which follows from   (\ref{gf.1j}).  Thus we see  that $L\Phi=I$.  The second equality in(\ref{ark.6}) follows similarly.  The  last two equalities in (\ref{ark.6})  follow immediately.

 We now obtain (\ref{wa}).  Note that  it follows from (\ref{7.27mm}) that for all $m,n\in\mN,$   
 \begin{equation}
\( \Phi \Phi ^{T}\)_{m,n}=\sum_{j=1}^{m\wedge n}   \phi_{m+1-j} \phi_{n+1-j}=E\(\wt \xi_{m}\wt \xi_{n}\) =\VV_{m,n}      \label{ark.8}.
 \end{equation}
  We show below that $A =L^{T}L$.  
Therefore,  %formally, by (\ref{ark.6}),
%\begin{equation}
%WA= \Phi \Phi ^{T}L^{T}L=\Phi \(\Phi ^{T}L^{T}\)L=\Phi L=I.\label{matrixj.2}
%\end{equation}
    \be 
    \sum_{j} \VV_{i,j}A_{j,l}= \sum_{j} \(\sum_{m}\Phi_{i,m} \Phi_{m,j} ^{T}\)\(\sum_{n}L_{j,n}^{T}L_{n,l}\).
   \label{matrixj.3}
   \ee

  It is easy to see that (\ref{matrixj.1}) holds, once we  show that we can interchange the order of summation 
in (\ref{matrixj.3}).  This allows  us to write,
    \bea 
    \sum_{j} \VV_{i,j}A_{j,l}&= & \sum_{m}  \sum_{n}  \Phi_{i,m}  \(\sum_{j} \Phi_{m,j} ^{T} L_{j,n}^{T}\) L_{n,l}
   \label{matrixj.3y}\\
   &= & \sum_{m}  \sum_{n}  \Phi_{i,m}  \de_{m,n} L_{n,l}=\sum_{n}  \Phi_{i,n}  L_{n,l}=\de_{i,l},\nn
   \eea
where we use   (\ref{ark.6}) twice. 
  
To show that we can interchange the order of summation 
in (\ref{matrixj.3}) it suffices to show that for $i$ and $l$ fixed all the sums in (\ref{matrixj.3}) are only over a finite number of terms that are not equal to 0. 
Making use of the fact that many of the terms in $L$ and $\Phi$ are equal to 0, we see that, 
\begin{equation}
   \sum_{m}\Phi_{i,m} \Phi_{m,j} ^{T}=   \sum_{m=1}^{i}\Phi_{i,m} \Phi_{m,j} ^{T}
   \end{equation}
  and  
  \begin{equation}
\sum_{n}L_{j,n}^{T}L_{n,l}=\sum_{n=l}^{l+k+1}L_{j,n}^{T}L_{n,l}.
   \end{equation}
     Furthermore, for each  $n$,   $L_{j,n}^{T}=0$   when $j>n$.     This shows that the summation in (\ref{matrixj.3}) is only over a finite number of terms.

  We show in (\ref{matrixj.3y})  that $\VV A=I$. Since both  $\VV$ and $A$ are symmetric, we also have $A\VV=I$.
   
 	To show that  $A=L^{T}L$ we take the product $L^{T}L$ to see that    
	\begin{equation}
A_{m,m}=(L^{T}L)_{m,m}=\sum_{j}L^{2}_{j,m}=\sum_{j }    p^{2}_{j-m}  =1+\sum_{i=1}^{  k}   p^{2}_{i},\label{ark.14}
  \end{equation}
 and for $n< m$,    
 \begin{eqnarray}
   (L^{T}L)_{m,n}&=&\sum_{j}L _{j,m}L _{j,n}=\sum_{j}   p_{j-m} p_{j-n} \\
  &=&-p_{m-n}+\sum_{j>m}    p_{j-m} p_{j-n} =-p_{m-n}+\sum_{l\geq 1}   p_{l}\,\, p_{(m-n)+l}\nn\\
  &=&-p_{m-n}+\sum_{l=1} ^{k-(m-n)}    p_{l}\,\, p_{(m-n)+l},    \label{ark.13} \nonumber
  \end{eqnarray}
 where we make the substitution   $l=j-m$ at the next to last step and use the fact that $p_{(m-n)+l}=0$ when $(m-n)+l>k$.
 
  Since $L^{T}L$ is symmetric we get the same result when $n$ and $m$ are interchanged. It is clear that when $|m-n|>k$, $L^{T}L_{m,n} =0$. This shows that $A=L^{T}L$.\qed

 The next lemma gives some properties of the matrix $A$.  Note that we are interested in the case in which $\VV$ is the potential density of a Markov chain. For this to be the case the off diagonal elements on $A$ must be negative.
 
 	 \begin{lemma} \label{lem-6.5mm} Let $A$ be as given in Lemma \ref{lem-7.5mm} and assume that
$\sum_{i=1}^{k}p_{i}\leq 1$. Then 
 \be 
  \sum_{n\in\mN}A_{m,n}=\(1-\sum_{i=1}^{ k}   p_{i}\)^{2},\qquad m>k. 
 \label{ark.15} 
  \ee
  Furthermore, when   $p_{i}\downarrow$, 
    \begin{equation}
   A_{m,n}\le 0,\qquad \forall\, n,m\in\mN, \, n\ne m,\label{6.62mm}
   \end{equation}
and 
 \be 
  \sum_{n\in\mN}A_{m,n}> \(1-\sum_{i=1}^{ k}   p_{i}\)^{2}, \qquad 1\le m\le k.
 \label{ark.15v} 
  \ee
  Therefore, $-A$ is a $Q$-matrix with uniformly bounded entries.   

 \end{lemma}
 
 \Proof  To prove (\ref{ark.15}) we note that by Lemma \ref{lem-7.5mm}, for $m>k$,  
  \begin{eqnarray}
  \sum_{n\in\mN}A_{m,n}&=&a_{0}+2 \sum_{j=1}^{k}a_{j}=1+\sum_{i=1}^{   k}   p^{2}_{i}  +  2\sum_{i= 1}^{k }\(-p_{i}+\sum_{l=1}^{k-i}   p_{l}\,\, p_{i+l}\)\nn
  \\
&=&\(1-\sum_{i=1}^{ k}   p_{i}\)^{2}. 
 \label{ark.15s} 
  \end{eqnarray}

  For (\ref{6.62mm}) we use (\ref{7.47mm}) to see that for all  $    1\le |m-n|\le k$,
  \be 
a_{|m-n|}= A_{m,n} \le -p_{|m-n|}+p_{|m-n|+1}\sum_{l=1}^{k-|m-n|}   p_{l} \le -p_{|m-n|}+p_{|m-n|+1}\le 0\label{7.47mm4}.
  \ee
  To get (\ref{ark.15v})  we note that by  (\ref{6.49mm})   the row sums of the first $k$ rows of $A$ omit some of the terms $a_{i}$, $1\le i\le k$, which are  less than or equal to 0.  

The final statement in the lemma   follows from (\ref{6.62mm}), (\ref{ark.15v}) and (\ref{7.46mm}).
 \qed
 
 \vspace{-.3in}
 \begin{remark} \label{rem-dec}{\rm It is clear that $-A$ can be a $Q$-matrix with uniformly bounded entries,  even when $\{p_{i}\}$ are not  decreasing. We see from Example \ref{ex-6.2} that when $k=2$,   $-A$ is always a $Q$-matrix with uniformly bounded entries.  Nevertheless, to keep the statement of Theorem \ref{theo-4.2mmfx} from being too cumbersome, we include the hypothesis that $p_{i}\downarrow.$
 }\end{remark}

The next theorem  ties certain $k$-th order linear regressions to   Markov chains.

   \bt\label{lem-invA} 
Assume that $p_{i}\downarrow$.   Then  $\VV$ is the potential density of  a  Markov chain  on $\mN$ 
 with $Q$-matrix, $-A$.  
  \et
  
  The proof of this theorem depends on the following general result.

  \bl\label{lem-invpot} Let $Q$ be the   Q-matrix of a  transient Markov chain $X$ on $\mN$ and assume that   $Q$ is    a  $(2m+1)$-diagonal matrix,  with
   \be
   \sup_{j\in\mN}  |Q_{j,j}|<\ff.\label{6.70mm}
   \ee
    Let $V $ be a matrix   satisfying,  
  \be
  VQ=-I,\qquad\mbox{and}\qquad \sup_{i\in\mN}|V_{k,i}|<\ff, \qquad \forall\, k\in\mN. \label{6.112a}
    \ee
Then $V$ is the potential  density of $X$, and in particular has positive entries.

If $\sum_{i}|V_{k,i}|<\ff,$ $\forall\, k\in\mN$, then the same results hold without the requirement that $Q$ is a   $(2m+1)$-diagonal matrix.
\el
  \Proof    Let $U$ be the potential  density  of $X$. 
 By  (\ref{inv.34a}),   
  \begin{equation}
 -\de_{i,l}= \sum_{j}  Q_{i,j}U_{j,l}.\label{ord.1}
  \end{equation}
  Therefore,  
  \begin{equation}
 -V_{k,l}=\sum_{i}V_{k,i}\sum_{j} Q_{i,j}U_{j,l}=\sum_{i}\sum_{j} V_{k,i}Q_{i,j}U_{j,l}.\label{ord.2}
  \end{equation}
 We show immediately below  that we can interchange the order of summation. Consequently,  by (\ref{6.112a}), for all $k,l\in\mN$, 
   \begin{equation}
- V_{k,l}=\sum_{j}\sum_{i} V_{k,i}Q_{i,j}U_{j,l}= -\sum_{j} \de_{k,j}U_{j,l}=-U_{k,l}.\label{ord.3}
  \end{equation}
  This shows that  $V$ is the potential   density   of $X$.
 
 To be able to interchange the order of summation in (\ref{ord.2}), we only need to show that 
 for each fixed $k$ and $l$,
   \begin{equation}
 \sum_{i}\sum_{j} |V_{k,i}| |Q_{i,j}|U_{j,l}<\ff.\label{ord.4}
  \end{equation}
  
 We have $U_{j,l}\leq U_{l,l}$ for all $j$, and  for each $j$ there are at most $2m+1$ elements   $|Q_{i,j}|$ that are not equal to 0. Therefore,  
  \bea
 \sum_{i}\sum_{j} |V_{k,i}||Q_{i,j}|U_{j,l}&\leq&  U_{l,l} \sum_{i}\sum_{j} |V_{k,i}||Q_{i,j}|=U_{l,l} \sum_{j}\sum_{i} |V_{k,i}||Q_{i,j}| \nn\\
  &\leq&  (2m+1)U_{l,l}\sup_{i} |V_{k,i}|\sum_{j}  |Q_{i,j}|.  \label{ord.5} 
  \eea 
Finally, using (\ref{unif.4x1}) we have
      \be 
\sup_{i} |V_{k,i}| \sum_{j}  |Q_{i,j}|   \label{ord.6}  
\le \(\sup_{i }|V_{k,i}|\) \,2\sup_{j}  |Q_{j,j}|<\ff. 
  \ee 
  Thus we get (\ref{ord.4}).
  
   If $\sum_{i}|V_{k,i}|<\ff,$ $\forall\, k\in\mN$, then in place of (\ref{ord.5}) we have 
    \bea
 \sum_{i}\sum_{j} |V_{k,i}||Q_{i,j}|U_{j,l}&\leq&  U_{l,l} \sum_{i}\sum_{j} |V_{k,i}||Q_{i,j}|=U_{l,l}\sum_{i} |V_{k,i}| \sum_{j}|Q_{i,j}| \nn\\
  &\leq& U_{l,l} \,\( 2\sup_{j}  |Q_{j,j}|\)\, \sum_{i} |V_{k,i}|<\ff.  \label{ord.5m} 
  \eea 
  \qed

\noindent{\bf Proof of Theorem \ref{lem-invA} }  The proof follows immediately from Lemma \ref{lem-invpot} once we show that the hypotheses of the lemma are satisfied. The fact that $-A$ is a $Q$-matrix is given in Lemma \ref{lem-6.5mm}. 

The property  that   $Q$ is   a $(2k+1)$-diagonal matrix,  the condition in (\ref{6.70mm}) and  the first condition in (\ref{6.112a})  are given in Lemma \ref{lem-7.5mm}.   

 The second condition in (\ref{6.112a}) is given in   (\ref{6.6new}).  \qed

 We now turn to the proof of  Theorem \ref{theo-4.2mmfx}. In this case we need   sharp estimates of the covariance $\VV$.
  To this end we introduce a generating function for $\{\phi_{n}\}$.   Set
    \begin{equation}
 g(x)= \sum_{n=0}^{\ff}\phi_{n}x^{n}=\sum_{n=1}^{\ff}\phi_{n}x^{n}\label{gf.2a},
 \end{equation}
 since $\phi_{0}=0$.
It follows from (\ref{gf.1jj}) that this converges for all $|x|<1$.

 \begin{lemma}\label{lem-6.5}       Let
 \begin{equation}
   P(x)=1-\sum_{l=1}^{k}p_{l}x^{l}.\label{6.48mm}
   \end{equation}
 Then for all $|x|<1$,
   \begin{equation}
 g(x)={x \over P(x)}.\label{gf.3}
 \end{equation}
 \end{lemma}
 
 \Proof We have
 \bea
   \sum_{n=1}^{\ff}\phi_{n}x^{n}&=&x+ \sum_{n=2}^{\ff}\phi_{n}x^{n}=x+ \sum_{n=2}^{\ff}\sum_{l=1}^{k}p_{l}\phi_{n-l}\, x^{n}\label{6.173}\\
   &=& x+ \sum_{l=1}^{k}p_{l}x^{l}\sum_{n=2}^{\ff}\phi_{n-l}\, x^{n-l},\nn
   \eea
 where we use the fact that $\phi_{n}=0$ for $n<1$. In addition
 \bea
   \sum_{n=2}^{\ff}\phi_{n-l}\, x^{n-l}&=&\sum_{n=2}^{l}\phi_{n-l}\, x^{n-l}+\sum_{n=l+1}^{\ff}\phi_{n-l}\, x^{n-l}\label{6.174}\\
   &=& \sum_{n=l+1}^{\ff}\phi_{n-l}\, x^{n-l}=\sum_{k=1}^{\ff}\phi_{k}\, x^{k}\nn=g(x).
   \eea
It follows from (\ref{6.173}) and (\ref{6.174}) that,  
\begin{equation}
 g(x) =x+ g(x)\sum_{l=1}^{k}p_{l}x^{l},\label{gf.2}
 \end{equation}
 which gives (\ref{gf.3}).\qed

 \begin{lemma} \label{lem-6.3nn}
Let $q_{1},\ldots, q_{k}$ be the  roots of $P(x)$
 which may be complex.  Then,
 \begin{itemize}
\item[(i)]
  $ \sum_{l=1}^{k}p_{l}=1 \iff   q_{1}=1$     is   a simple root    and   $|q_{l}|>1$, $l=2,\ldots,k.$ 
\vspace{.05 in}  \item[(ii)]
$  \sum_{l=1}^{k}p_{l}<1\iff     |q_{l}|>1  , \quad  l=1, \ldots,k. $
   \end{itemize}
 \end{lemma}

 \Proof  Assume first that $\sum_{l=1}^{k}p_{l}=1$. Then it is obvious that $q_{1}=1$ is   a root. Furthermore since,  
 \begin{equation}
P'(1)= -\sum_{l=1}^{k}lp_{l}<0, \label{gf.4a}
 \end{equation}
it is not a multiple root.   
Also, note that  
  \begin{equation}
 1=\Big |\sum_{l=1}^{k}p_{l}x^{l}\Big |\leq  \sum_{l=1}^{k}p_{l}|x|^{l},\label{gf.5}
 \end{equation}
  with strict inequality when   $  |x|=1$ and $  x\ne 1$. Therefore,  $|q_{l}|>1$ for all  $l= 2,\ldots,k$. 

If $\sum_{l=1}^{k}p_{l}<1$ it is clear from (\ref{gf.5}) that $|q_{l}|>1$ for all  $l= 1,\ldots,k$. 
\qed 
 
 	We now give a formula for $\phi=\{\phi_{n}\}$. Define 
 \begin{equation}
  c_{1}=\frac{1}{\sum_{l=1}^{k}lp_{l}}. \label{csub1}
 \end{equation}

 	\begin{lemma}\label{lem-6.7} Let $P(x)$ be as given in (\ref{6.48mm}) and  assume that  it   has distinct roots  $q_{l}$ of degree $d_{l}$, $l=1,\ldots, k'$. 
\begin{itemize}
\item[(i)] If $ \sum_{l=1}^{k}p_{l}<1$, where  $k=\sum_{l=1}^{k'}d_{l} $, 
then all  $|q_{l}|>1$ and,
\begin{equation}
 \phi_{n}= \sum_{l=1}^{k'}\sum_{j=1}^{d_{l}}  B_{j}(q_{l}) \times {j-1+n \choose j-1}\( \frac{1 }{q_{l}}\)^{n}\label{get.10},
\end{equation}
where
\begin{equation}
   B_{j}(q_{l})={(-1)^{j} \over  q_{l}^{j}(  d_{l}-j)!}\lim_{x\to q_{l}}   D^{(d_{l}-j)}{x(x-q_{l})^{d_{l}} \over P(x)} . \label{6.57mm}
   \end{equation}
   Furthermore,
   \begin{equation}
   \|\phi\|_{1}=
\frac{1}{P(1)}<\ff \label{6.58nn}.\end{equation}

\item[(ii)]
If $ \sum_{l=1}^{k}p_{l}=1$ the roots of $P(x)$ can be arranged so that $q_{1} =d_{1}=1$ and   $|q_{l}|>1$, for $l=2,\ldots, k'$. In this case,
\begin{equation}
 \phi_{n}= c_{1} +\psi_{n},\label{6.59ww}
\end{equation}
where,
\begin{equation}
   \psi_{n}=\sum_{l=2}^{k'}\sum_{j=1}^{d_{l}}  B_{j}(q_{l}) \times {j-1+n \choose j-1}\( \frac{1 }{q_{l}}\)^{n}\label{get.10ii}.
   \end{equation}
Furthermore,  
 \begin{equation}
      \|\psi\|_{1}<\ff.
   \end{equation}
      \end{itemize}
 \end{lemma}

 	\Proof  Suppose  more generally that  $P(x)$ is a polynomial with $P(0)\neq 0$ and  distinct roots  $q_{l}$ of degree $d_{l}$, $l=1,\ldots, k'$.  Then we can write, 
\begin{equation}
{x \over P(x)}=\sum_{l=1}^{k'}\sum_{j=1}^{d_{l}}{a_{l,j} \over (x-q_{l})^{j}}\label{got.1}
\end{equation}
where
\begin{equation}
a_{l,j}={1 \over (d_{l}-j)!}  \lim_{x\to q_{l}}  D^{(d_{l}-j)}{x(x-q_{l})^{d_{l}} \over P(x)} .\label{got.2}
\end{equation}

For lack of a  suitable reference we provide  a simple proof. Let
\begin{equation}
f(x)={x \over P(x)}-\sum_{l=1}^{k'}\sum_{j=1}^{d_{l}}{a_{l,j} \over (x-q_{l})^{j}}.\label{got.1a}
\end{equation}
The function $f(x)$ is a rational function which can only have finite poles at $q_{l}$ of degrees $\leq d_{l}$, $l=1,\ldots, k'$. Consider 
\bea
&&(x-q_{l})^{d_{l}}f(x)\label{6.65nn}\\
&&\qquad
={x(x-q_{l})^{d_{l}} \over P(x)}-\sum_{l'=1,l'\ne l}^{k'}\sum_{j=1}^{d_{l'}}{a_{l',j}(x-q_{l})^{d_{l}} \over (x-q_{l'})^{j}}+\sum_{j=1}^{d_{l}}a_{l,j}(x-q_{l})^{d_{l}-j}. \nn
   \eea
  Considering the definition  of the $a_{l,j}$ in (\ref{got.2}), we see that,
  \begin{equation}
  \lim_{x\to q_{l}}  D^{(d_{l}-j)} (x-q_{l})^{d_{l}}f(x)=0,\label{6.65mm}
  \end{equation} 
 for all $1\le j\le d_{l}$, and all $l=1,\ldots, k'$. This shows that the rational function $f(x)$ has no finite poles, which implies that  $f(x)$ is a polynomial,  and since $\lim_{x\to\ff}f(x)=0$, we must have $f(x)\equiv 0$.  Using (\ref{got.1a})  we get (\ref{got.1}).

Let    \begin{equation}
   B_{j}(q_{l})={a_{l,j}(- {1})^{j} \over q^{j}_{l}}.
   \end{equation} 
Then    if all   the $|q_{l}|>1$  it follows from  (\ref{got.1}) that for all $|x|\le 1$, 
\bea
{x\over P(x)}&=&\sum_{l=1}^{k'}\sum_{j=1}^{d_{l}}{a_{l,j}(- {1})^{j} \over q^{j}_{l}(1-x/q_{l})^{j}} \label{got.3}\\
&=&\sum_{n=0}^{\ff}\sum_{l=1}^{k'}\sum_{j=1}^{d_{l}}  B_{j}(q_{l}) \times {j-1+n \choose j-1}\( \frac{1 }{q_{l}}\)^{n}x^{n}.\nn
\eea 
Therefore, using (\ref{gf.2a})  and (\ref{gf.3}) we see that for all $|x|<1$ 
\begin{equation}
\sum_{n=0}^{\ff} \phi_{n}x^{n}=\sum_{n=0}^{\ff}\sum_{l=1}^{k'}\sum_{j=1}^{d_{l}}  B_{j}(q_{l}) \times {j-1+n \choose j-1}\( \frac{1 }{q_{l}}\)^{n}x^{n}.\label{got.3v}
\end{equation}
  This proves (\ref{get.10}). Since   all   $|q_{l}|>1$ it is clear that (\ref{got.3v}) converges for $x=1$
so that by combining the last two displays  we see that,
\begin{equation}
{1\over P(1)}= \sum_{n=0}^{\ff} \phi_{n}.
\end{equation}
Since by (\ref{gf.1j}), $\phi_{n}\ge 0$  for all $n\in\mN$, we  get  (\ref{6.58nn}). 
 
\medskip	For $(ii)$	we see that  as in  (\ref{got.3})   for all $|x|<1$, 
\be 
{x\over P(x)} 
=\sum_{n=0}^{\ff}\(B_{1}(1)+\sum_{l=2}^{k'}\sum_{j=1}^{d_{l}}  B_{j}(q_{l}) \times {j-1+n \choose j-1}\( \frac{1 }{q_{l}}\)^{n} \)x^{n},
 \ee 	
 where
 \begin{equation}
   B_{1}(1)=-\lim_{x\to 1}\frac{x(x-1)}{P(x)}=-\frac{1}{P'(1)}=c_{1},\label{6.68mm}
    \end{equation}	
by   L'Hospital's Rule and  (\ref{gf.4a}).	This gives (\ref{6.59ww}), in which
\begin{equation}
   \psi_{n}=\sum_{l=2}^{k'}\sum_{j=1}^{d_{l}}  B_{j}(q_{l}) \times {j-1+n \choose j-1}\( \frac{1 }{q_{l}}\)^{n}.
   \end{equation}
Since $ |q_{l}|>1$, $2\le l\le k' $, it is clear  that  $  \|\psi\|_{1}<\ff$.     \qed

 	 \begin{example}\label{ex-6.2a} {\rm  

 	 Suppose that   $P(x)$ has  real   roots, $a$ and   $-b$, where   $a$ has multiplicity 1 and $-b$ has multiplicity 2,  and $a \ge 1$.  In this case 
\bea 
 P(x)&=&-\frac{1}{ab^{2}} (x-a)(x+b) ^{2}  \\
 &=& 1-\frac{1}{ab^{2}} \(x^{3}+(2b-a)x^{2}+(b^{2}-2ab)x\)
 \eea 
 Therefore,
 \begin{equation}
   p_{1}=\frac{b^{2}-2ab}{ab^{2}},\qquad   p_{2}=\frac{2b-a}{ab^{2}},\qquad 
p_{3}=\frac{1}{ab^{2}}.   \end{equation}
 When $b>2(a+1)$, $p_{1}> p_{2}> p_{3}$. (We know from Lemma \ref{lem-6.3nn} that we must have $a\ge 1$ and that $\sum_{j=1}^{3}p_{j}\le 1$ and is equal to 1 if and only if $a=1$.)
 
 We have
  \begin{equation}
   B_{1}(a)={(-1)  \over  a }\lim_{x\to a}    {x(-ab^{2})  \over (x+b)^{2}}=\frac{ab^{2}}{ (a+b)^{2}} . \label{6.5s}
   \end{equation}
  \be 
   B_{1}( -b)={1 \over  b} (-ab^{2})\lim_{x\to -b}   D^{(1)}{x  \over (x-a)} = \frac{a^{2} b}{ (a+b)^{2}} \label{6.5t}     
   \ee 
   \begin{equation}
   B_{2}(-b)={-ab^{2} \over  b^{2}} \lim_{x\to -b}   {x   \over (x-a)}=-\frac{ab}{(a+b)}\label{6.5u}
   \end{equation}
   Therefore,
   \bea
   \phi_{n}&=&\frac{ab^{2}}{ (a+b)^{2}}\(\frac{1}{a}\)^{n}-\(\frac{ab}{(a+b)}(n+1)-\frac{a^{2} b}{ (a+b)^{2}}\)\(\frac{1}{-b}\)^{n}\\
   &=&\frac{ab^{2}}{ (a+b)^{2}}\(\frac{1}{a}\)^{n}+(-1)^{n+1}\(\frac{abn}{(a+b)} +\frac{a  b^{2}}{ (a+b)^{2}}\)\(\frac{1}{ b}\)^{n}
   \eea
When $a=1$ this is,
   \be 
   \phi_{n}=\frac{ b^{2}}{ (1+b)^{2}} +(-1)^{n+1}\(\frac{ bn}{(1+b)}+\frac{ b^{2}}{ (1+b)^{2}}\)\(\frac{1}{ b}\)^{n}. 
   \ee 
   One can check that in this case,
   \begin{equation}
   \sum_{j=1}^{3}j p_{j}=\frac{(1+b)^{2}}{b^{2}}.
   \end{equation}

     }\end{example}

  	 \noindent{\bf Proof of    Theorem \ref{theo-4.2mmfx}, (\ref{7.86}) }
We use Theorem \ref{theo-mchain}. To begin   we obtain the denominator in (\ref{121.4a}).  
    Let $  \xi=\{  \xi_{n},n\in\mN \}$ be a Gaussian sequence defined exactly  as $\wt\xi$ is defined  in   (\ref{ark.1}) but with the additional  conditions that $\sum_{l=1}^{k}p_{l}=1$.  We now show that  
 \begin{equation}
   \limsup_{n\to\ff}\frac{\xi_{n}}{(2n\log\log n)^{1/2}}=\frac{1}{ \sum_{l=1}^{k}lp_{l} }\qquad a.s.\label{7.5mm}
   \end{equation}

 	\medskip	It follows from  Lemma \ref{lem-6.7} $(ii)$   that,  
  \begin{equation}
 \phi_{n}=c_{1}+\psi_{n}, \hspace{.2 in} \mbox{where} \hspace{.2 in} \psi \in \ell_{1}^{+}.\label{kov.1}
 \end{equation}
  By (\ref{ark.2}) we can write   
\be
\xi_{n}=c_{1}S_{n}+\rho_{n},\label{7.31mm}
\ee where    
 \begin{equation}
S_{n}=\sum_{j=1}^{n}g_{j} \hspace{.2 in}\mbox{and}\hspace{.2 in} \rho_{n}=\sum_{j=1}^{n} \psi_{n+1-j}g_{j}.\label{kov.2}
 \end{equation}
Note that $E\(\rho_{n}^{2}\)\le \|\psi\|_{2}^{2}$ for all $n\in\mN$. It follows from the Borel-Cantelli Lemma that,
  \begin{equation}
 \varlimsup_{j\to \ff} {|\rho_{j}| \over \sqrt{2 \log j}}\le \|\psi\|_{2}  \qquad 
a.s.\label{kov.3}
\end{equation}
It now follows from  (\ref{7.31mm})  and the standard law of the iterated logarithm for $S_{n}$ that (\ref{7.5mm}) holds.  

\medskip	We now show that (\ref{condreal}) holds.  Let $\VV=\{\VV_{j,k};j,k\in \mN\}$ be as in   (\ref{7.27mm}).  We now find an 
  estimate for the row sums of   $ \(  \VV(l,n)\)^{-1}$. 
 	For $n\geq k$ set  
 \be
 \Xi(l,n)=(\xi_{l+1},\xi_{l+2},\ldots, \xi_{l+n})\label{5.90}
 \ee
 and
 \be  
G(l,n)=(   \eta_{l+1},\ldots, \eta_{l+k},  g_{l+k+1},g_{l+k+1},\ldots, g_{l+n}),
 \ee 
   where
   \begin{equation}
  \eta_{l+j}= \xi_{l+j}-\sum_{i=1}^{j-1}p_{i}\xi_{l+j-i}, \hspace{.2 in}j= 1\ldots, k. \label{5.92}
   \end{equation}
Note that this is similar in form to (\ref{ark.1}), but starting from $l+1$.

 \medskip	  We use several lemmas. The first one is easy to verify. 
   \begin{lemma} \label{lem-5.8}
  \begin{equation}
 G(l,n)^{T}=L (l,n)\Xi(l,n) ^{T},\label{autoj.8k}
 \end{equation}
 where $L$ is given in (\ref{ark.4}).
  \end{lemma}
 
 	  It follows  from (\ref{autoj.8k}) that, 
 \begin{equation}
G(l,n)^{T}G(l,n)=L(l,n)\Xi(l,n)^{T} \Xi(l,n) L(l,n)^{T}. \label{autoj.10k}
 \end{equation}
We take the  expectation of each side and get the vector equation,  
 \begin{equation}
B\otimes I_{n-k}=L(l,n) W(l,n)L(l,n)^{T},\label{autoa.11k}
 \end{equation}
 where  
 \begin{eqnarray}
 &&B=\mbox{Cov}(  \eta_{l+1},\ldots, \eta_{l+k}).   \nonumber
 \end{eqnarray} 
 
\begin{lemma} \label{lem-6.8mm} 
 \begin{equation}
\( \VV(l,n)\)^{-1}{\bf  1}_{n}= \begin{pmatrix}
\( \VV(l,k)\)^{-1}{\bf  1}_{k}    \\
0        
 \end{pmatrix},\label{ntok.1}
 \end{equation}
 where ${\bf  1}_{m}$ denotes an $m$ dimensional column vector  with all  of its components equal to 1.
 \end{lemma}

 Note that $\( \VV(l,n)\)^{-1}{\bf  1}_{n}$ is an $n$ dimensional vector with components that are the row sums of $\( \VV(l,n)\)^{-1}$. Therefore, (\ref{ntok.1})
states that the first $k$ row sums of $ \( \VV(l,n)\)^{-1}$ are   equal to the row sums of  $ \( \VV(l,k)\)^{-1}$, and the remaining row sums are equal to $0$.

\medskip	\Proof 
Using (\ref{autoa.11k}) we see that    

 \begin{equation}
 \( \VV(l,n)\)^{-1}= L(l,n)^{T} \(B^{-1} \otimes I_{n-k}\) L(l,n). \label{autoa.11kv}
 \end{equation}
 In addition, since $L(l,n)$ is a lower triangular matrix  we can write it in the block form,
\be  L(l,n)=\begin{pmatrix}
 F &0    \\
G&H       
 \end{pmatrix},\label{7.70mm}
 \ee
 where $F$ is a $k\times k$ matrix.
It is easy to check that
 \be \hspace{.2 in}\(L(l,n)\)^{-1}=\begin{pmatrix}
 F^{-1} &0    \\
-H^{-1}G F^{-1} &H^{-1}       
 \end{pmatrix}.
 \ee
 We also note that since all row sums of $L(l,n)$ after the $k$-th row are equal to zero,
   \be 
L(l,n){\bf  1}_{n}=\begin{pmatrix}
 F{\bf  1}_{k}    \\
0      
 \end{pmatrix}.\label{7.72mm}
 \ee
  It  follows from (\ref{autoa.11kv}) that 
  \begin{equation}
 \( \VV(l,n)\)^{-1}{\bf  1}_{n}= L(l,n)^{T} \(B^{-1} \otimes I_{n-k}\) L(l,n){\bf  1}_{n}. \label{pp.11kv}
 \end{equation}
  Using   (\ref{7.72mm}) we see that,
  \begin{equation}
    \(B^{-1} \otimes I_{n-k}\) L(l,n){\bf  1}_{n}=\begin{pmatrix}
B^{-1} F{\bf  1}_{k}    \\
0      
 \end{pmatrix}.\label{7.75mm}
   \end{equation}
Consequently,
 \begin{equation}
 \( \VV(l,n)\)^{-1}{\bf  1}_{n}=\begin{pmatrix}
 F^{T} &G^{T}    \\
0&H^{T}       
 \end{pmatrix}\begin{pmatrix}
B^{-1} F{\bf  1}_{k}    \\
0      
 \end{pmatrix}= \begin{pmatrix}
F^{T}B^{-1} F{\bf  1}_{k}    \\
0      
 \end{pmatrix}.\label{ntok.1a}
 \end{equation}
  On the other hand, 
 by
 (\ref{autoa.11k}) , 
 \begin{gather*}
\begin{pmatrix}
 F^{-1} &0    \\
-H^{-1}G F^{-1} &H^{-1}       
 \end{pmatrix} \begin{pmatrix}
 B &0    \\
0&I      
 \end{pmatrix}\begin{pmatrix}
 F^{-1} &0    \\
-H^{-1}G F^{-1} &H^{-1}       
 \end{pmatrix}^{T} = \VV (l,n),
 \end{gather*}
from which we obtain 
\begin{equation}
F^{-1}B(F^{T})^{-1}= \VV(l,k).\label{autofd.1}
\end{equation}
Consequently,
\be
F^{T}B^{-1}F=\( \VV(l,k)\)^{-1}.\label{7.77mm}
\ee 
Using this and (\ref{ntok.1a}) we get (\ref{ntok.1}).\qed

We now consider $\VV(l,k)$.  
 
 	\begin{lemma}\label{lem-6.10mm} When $\sum_{l=1}^{k}p_{l}=1$,
	  \begin{equation}
    E\(\xi_{m}\xi_{n}\)=c_{1}^{2}(m\wedge n)+a_{n,m},\label{6.104mm}
   \end{equation}
where $|a_{m,n}|\le D<\ff$, for all $m,n\in\mN$.	
	
\end{lemma}	 

 \Proof   
    By (\ref{7.27mm}) and Lemma  \ref{lem-6.7} $(ii)$,  when $m\leq n$, we have  
 \bea
  E\(\xi_{m}\xi_{n}\)&=&\sum_{j=0}^{m -1}\phi_{m-j}\phi_{n-j} \label{auto.16rk}\\
&= & c^{2}_{1}m+ c_{1}  \sum_{j=0}^{m -1}   \psi_{ m-j} +c_{1}  \sum_{j=0}^{m -1}   \psi_{ n-j}    \nn\\
&  &\quad + c_{1}^{2}\sum_{j=0}^{m -1} \psi_{ m-j}  \psi_{ n-j}  . \nn
\eea
Clearly, for all $  p\ge m$,  
  \begin{equation}
   \bigg |\sum_{j=0}^{m-1}  \psi_{ p-j} \bigg| \le     \sum_{j=1}^{\ff}  |\psi_{  j}|=\|\psi \|_{1},
   \end{equation}
and,  
 \bea 
 \bigg| \sum_{j=0}^{m -1} \psi_{m-j}  \psi_{n-j} \bigg |
   &\le&  \sum_{j=1}^{\ff  }|\psi_{ j}|^{2}= \|\psi \|_{2}^{2}\nn,\eea
where we use the Schwartz Inequality.
  Combining all these inequalities we   see that for $m\le n$,
  \begin{equation}
    E\(\xi_{m}\xi_{n}\)=c_{1}^{2}(m\wedge n)+a_{n,m},\label{6.104mmo}
   \end{equation}
   where, 
  \begin{equation}
| a_{m,n}|\le 2  c_{1}  \|\psi \|_{1}+\(c_{1} \|\psi \|_{2}\)^{2}:= D<\ff.\label{6.89mh}
   \end{equation}\qed

 The next lemma is used to obtain (\ref{condreal}).

 	\begin{lemma} \label{lem-7.8}   For all $1\le i\le k$,
\begin{equation}
 \sum _{j=1 }^{k}  \VV(l,k )^{i,j} =O\(1/l\).\label{ntok.2}
   \end{equation}
   \el 

 	\Proof    
  It follow from  Theorem \ref{lem-invA} that $ \VV$ is the potential  density   of a Markov chain. Therefore so is $ \VV(l,k)$. Consequently, $\VV(l,k)^{-1}$ is an M-matrix with positive row sums.  This gives the first inequality in (\ref{7.88b})    below,
 \begin{equation}
| \VV(l,k)^{j,i}|\leq  \VV (l,k)^{j,j}\leq  A_{j,j}\leq 2.\label{7.88b}
 \end{equation}
 The second inequality in (\ref{7.88b}) follows from Lemma \ref{lem-dom}, below.    The third inequality in (\ref{7.88b}) is given in   (\ref{ark.14}).
  
Clearly,
 \begin{equation}
\sum_{j=1 }^{k}  \VV(l,k)_ {i,j}   \VV(l,k)^{j,i}=1,\qquad 1\le i\le k.\label{7.89}
 \end{equation}
Furthermore, by Lemma \ref{lem-6.10mm},  
 \bea
  1&=&\sum _{j=1 }^{k} \VV(l,k)_ {i,j}   \VV(l,k)^{j,i}\label{7.90}\\
  & =& c_{1}^{2}  \sum _{j=1 }^{k}  ((l+i)\wedge (l+j))\VV(l,k)^{i,j}+   \sum _{j=1 }^{k}  a_{l+i,l+j}  \VV(l,k)^{j,i}\nn\\
   &=&c_{1}^{2} l \sum _{j=1 }^{k} \VV(l,k)^{i,j}+c_{1}^{2}  \sum _{j=1 }^{k}  (i\wedge j)\VV(l,k)^{i,j}+   \sum _{j=1 }^{k}  a_{l+i,l+j}  \VV(l,k)^{j,i}.\nn
   \eea
Therefore,
\bea
   c_{1}^{2} l \sum _{j=1 }^{k}  \VV(l,k)^{i,j}&\le &1+(c_{1}^{2}k+D) \sum _{j=1 }^{k}  |\VV(l,k)^{i,j}|\\
   &\le &1+2k(c_{1}^{2}k+D) ,\nn
   \eea
 where we use  (\ref{7.88b}).   This  gives (\ref{ntok.2}).\qed

   \bl\label{lem-dom}
  Let $X\!=\!
(\Om,  \FF_{t}, X_t,\th_{t},P^x
)$ be a transient symmetric Borel right process with state space $ \mN$,  and  potential   densities    $U=\{U_{j,k}$, $j,k\in\mN\}$ and $Q$-matrix, $Q$.	
Assume that 
\begin{equation}
   U_{j,k}>0\qquad\mbox{and}\qquad  |Q_{j,j}|<\ff,\qquad \forall\, j,k\in\mN. 
   \end{equation}
Then for any distinct sequence  $l_{1},l_{2},\ldots, l_{n}$ in $ \mN$, the matrix $K =\{U_{l_{i},l_{k}}\}_{i,j=1}^{n}$ is invertible and,  
\begin{equation}
K^{j,j}\leq |Q_{l_{j},l_{j}}|,\qquad \forall\,1\leq j\leq n. \label{dom.1}
\end{equation}
\el

 	\Proof  We follow the proof of   \cite[Lemma A.1]{MRnec}.  For all $k\in \mN$  set,  
\be
L_{t}^{k}=\int_{0}^{t}1_{\{X_{s}=k\}}\,ds. \label{dom.0}
\ee
It follows   from this that for all or all $j,k\in \mN$  we have,
\begin{equation}
U_{j,k}=E^{j}\(L_{\ff}^{k}\).
\end{equation}
Define the   stopping time,
\begin{equation}
 \si =\inf \{ t\geq 0\,|\, X_{ t}\in \{ l_{1},l_{2},\ldots, l_{n}, \De\}\cap \{X_{0}
\}^{ c}\}\label{ap.1}
\end{equation}
which may be infinite. By    \cite[(A.5)]{MRnec},
\begin{equation}
K^{j,j}\leq {1 \over E^{l_{j}}\(L^{l_{j}}_{\si}\)}.\label{dom.2}
\end{equation}
On the other hand,    the amount of time  $X_{t}$, starting at $l_{j}$, remains at $l_{j}$  is, 
\begin{equation}
 \si_{j}:=\inf \{ t\geq 0\,|\, X_{ t}\in   \{l_{j}
\}^{ c}\},\label{ap.2}
\end{equation}
which implies, by (\ref{dom.0}) that,
\begin{equation}
L^{l_{j}}_{\si_{j}}=\si_{j}.\label{6.121}
  \end{equation}
In addition, $ \si_{j}\leq  \si$, so that $E^{l_{j}}\(L^{l_{j}}_{\si_{j}}\)\leq E^{l_{j}}\(L^{l_{j}}_{\si}\)$. Therefore, it follows from (\ref{dom.2}) and (\ref{6.121}) that
\begin{equation}
K^{j,j}\leq {1 \over E^{l_{j}}\(L^{l_{j}}_{\si_{j}}\)}={1 \over E^{l_{j}}\(\si_{j}\)}.\label{dom.3}
\end{equation}
 Since $\si_{j}$  is an exponential random variable with mean  $1/|Q_{l_{j},l_{j}}|$; (see  \cite[Section 2.6]{Norris}), we get (\ref{dom.1}).\qed

   We now consider the potentials corresponding to $ \VV$.

  \begin{lemma}\label{lem-4.3ja}  Let    $f= \VV h$, where  
  $h\in \ell_{1}^{+}$.
Then
\begin{equation}
   f_{j}= g( {j})+\rho_{j} ,\qquad \forall j\in \mN,\label{4.18jja}
   \end{equation}
   where $g$ is an increasing strictly concave function and $\sup_{j}\rho_{j}= d\|h\|_{1}$ for some finite constant $ d$.
    \end{lemma} 

\Proof  We show in (\ref{6.104mm}) that, 
\begin{equation}
  \VV_{j,k}=c_{1}^{2}(j\wedge k)+  a_{j,k},\label{7.42mm}
   \ee 
where $|a_{j,k}|\le d$.
Therefore
\begin{equation}
   f_{j}=c_{1}^{2}\sum_{k=1}^{\ff}(j\wedge k)h_{k}+\sum_{k=1}^{\ff} h_{k}a_{j,k}.\label{7.43mm}
   \end{equation}
The lemma   now follows from  Theorem \ref{theo-3.1f}.\qed
  
   The next lemma shows that   (\ref{condreal}) holds.   
 
 \begin{lemma} \label{lem-6.11} Let    $f= \VV h$, where  
  $h\in \ell_{1}^{+}$. Then
\begin{equation}
\sum_{ j,p=1}^{n}( \VV(l,n))^ {j,p}f_{p+l}=o_{l}\(1\), \quad\mbox{ uniformly in }n.\label{6.94x}
\end{equation}

 \end{lemma}
 
 \Proof It follows from Lemmas \ref{lem-6.8mm} and \ref{lem-7.8} that    for all $l$ sufficiently large, there exists a constant $C$ such that,
  \bea
\sum_{ j,p=1}^{n}( \VV(l,n))^ {j,p}f_{p+l}&=&\sum_{  p=1}^{n}f_{p+l}\sum_{j=1}^{n}( \VV (l,n))^ {p,j} \\
&=&\sum_{  p=1}^{k}f_{p+l}\sum_{j=1}^{k}( \VV (l,k ))^ {p,j} \nn\le C\frac{f_{l+k}}{  l }.
\eea 
 By Lemma \ref{lem-4.3ja}, $f(j)=o(j)$ and since   $k$ is a fixed number, we get (\ref{6.94x}).\qed

		\medskip	\noindent{\bf Proof of Theorem \ref{theo-4.2mmfx}  (\ref{7.86}) continued}  
This  follows from Theorem \ref{theo-mchain}. Lemma \ref{lem-6.11} shows that (\ref{condreal}) holds.  The limit result  in  (\ref{7.5mm}) identifies the denominator in (\ref{121.4a}),   and  
Lemma \ref {lem-4.3ja} gives (\ref{1.5mm}).  \qed

\noindent{\bf Proof of Theorem \ref{theo-4.2mmfx},  (\ref{7.87}) } This follows from Theorem \ref{theo-1.8mm}. We show that the hypotheses in (\ref{1.39}) are satisfied.  
It follows from    (\ref{6.6new}) that $\inf_{j\geq 1}  \VV_{j,j}\newline =1$.  
 Therefore, the first condition in (\ref{1.39}) is satisfied.  
In addition, by (\ref{7.27mm}), when $n\ge m$,
\begin{equation} 
\VV_{m,n}=\sum_{j=0}^{m -1}\phi_{m -j}\phi_{n -j}=\sum_{j=1}^{m }\phi_{ j}\phi_{n-m+j}\label{6.118}
   \end{equation}
Therefore,
\begin{equation} 
\sum_{n=m}^{\ff} \VV_{m,n} =\sum_{j=1}^{m }\phi_{ j}\sum_{n=m}^{\ff}\phi_{n-m+j}\label{6.119}\le \|\phi\|_{1}^{2}.
   \end{equation}
Obviously, this holds when $n<m$ so we see that the second condition in (\ref{1.39}) is also satisfied. 

Furthermore, we see that
\begin{equation}
   \lim_{n\to\ff}    \VV_{n,n}=  \|\phi\|_{2}^{2}:=c^{*}.\label{6.116}
   \end{equation}Therefore, (\ref{7.87})  follows from Theorem \ref{theo-1.8mm}.

To obtain the upper bound in (\ref{1.55mm})
 we note that by (\ref{ark.1}),  
 \bea
 E\(\wt \xi_{n}^{\,2}\)&=&E\(\sum_{l=1}^{k}p_{l} \wt\xi_{n-l}\)^{2}+1\label{ark.1q}\\
 &=& \sum_{l,l'=1}^{k}p_{l} p_{l'}E(\wt\xi_{n-l}\wt\xi_{n-l'}) +1\nn\\
 &\le & \(\sum_{l =1}^{k}p_{l}  \)^{2}  E\(\wt \xi_{n}^{\,2}\)   +1\nn.
\eea
 Here we use the   Cauchy-Schwarz Inequality and  the fact that $E(\xi_{n}^{2})\uparrow$ to get, 
 \begin{equation}
   E(\wt\xi_{n-l}\wt\xi_{n-l'})\le \(E( \wt \xi_{n-l}^{\,2})E( \wt \xi_{n-l'}^{\,2})\)^{1/2}\le  E\(\wt \xi_{n}^{\,2}\) .
     \end{equation}
      The lower bound is obtained from (\ref{gf.1j}). We can add additional terms in situations where it is useful.
      
  The fact that $f\in c^{+}_{0}$ if and only if  $f= \VV  h$, where $   h  \in c_{0}^{+}$ follows from Lemma \ref{lem-8.2nn} once we show that (\ref{8.14}) holds.   
To see this we note that.
\bea
 \sum_{m=1}^{n/2}   \VV_{m,n}  &= & \sum_{m=1}^{n/2}\sum_{j=1}^{m }\phi_{ j}\phi_{n-m+j} 
=   \sum_{j=1}^{m }\phi_{ j}\sum_{m=1}^{n/2}\phi_{n-m+j}\nn\\
& 
   \leq& \|\phi\|_{1} \sum_{k=n/2}^{\ff}\phi_{k}.
\eea
 \qed
  
 \vspace{-.3 in}
 \begin{remark} \label{rem-6.2} {\rm It follows from Lemma  \ref{lem-6.7} $(i)$ that 
\be 
c^{*}=\|\phi\|_{2}^{2}= \sum_{l,l'=1}^{k'}   \sum_{j=1}^{d_{l}} \sum_{j'=1}^{d_{l'}}B_{j}(q_{l})B_{j'}(q_{l'})   F_{j,j'}(q_{l}q_{l'}),   \label{7.4nn}
\ee
where $  B_{j}(q_{l})$ is given in (\ref{6.57mm})
and
\begin{equation}
 F_{j,j'}(q_{l}q_{l'})= \sum_{n=0}^{\ff }{j-1+n \choose j-1}{j'-1+n \choose j'-1} \(\frac{1}{q_{l}q_{l'}}\)^{n}.\label{get.15}
\end{equation}
 }\end{remark}

\begin{example} {\rm Suppose that   
\be 
 P(x)=-\frac{1}{ab} (x-a)(x +b)  =1-\frac{1}{ab}(x^{2}+(b-a)x)
 \ee 
 where $a>1$ and $b\ge a$. This assures us that $p_{1}$ and $p_{2}>0$ and that $p_{1}+p_{2}<1.$
 
 We have
 \begin{equation}
   B_{1}(a)=\frac{-1}{a}\lim_{x\to a}\frac{x (x-a)(-ab)}{P(x)}=\frac{ab}{a+b}.
   \end{equation}
 Similarly,
  \begin{equation}
   B_{1}(-b)=\frac{ 1}{ b}\lim_{x\to -b}\frac{x (x+b)(-ab)}{P(x)}=-\frac{ab}{a+b}.
   \end{equation}
 	 Consequently,
 \bea 
c^{*} &=& \(\frac{ab}{a+b}\)^{2}   \(    F_{1,1}(a^{2})+ F_{1,1}(b^{2})-2F_{1,1}(a(-b))\),   \label{7.4d}\\
&=& \(\frac{ab}{a+b}\)^{2}   \(   \frac{a^{2}}{a^{2}-1}+ \frac{b^{2}}{b^{2}-1}-2\frac{ab}{ab+1}\), 
\eea

 For a concrete example suppose that 
  $a=-1+\sqrt5$ and $-b= -(1+\sqrt5)$. (These are the roots of
$
  x^{2}/4+ {x/2}-1
$.) Then,
    \be 
      \|\phi\|^{2}_{2}=  \frac{4}{5}\(\frac{6-2\sqrt5}{5-2\sqrt5}+\frac{6+2\sqrt5}{5+2\sqrt5}-\frac{8}{5}\)=\frac{48}{25}\approx 1.92.\nn
   \ee 
  (The bound in (\ref{1.55mm}) is 16/7$\approx$ 2.28.)
   }\end{example}

 	\medskip	\noindent{\bf Proof of Theorem \ref{theo-rest}   }    Consider $ \{\wt \YY_{\al,t_{j}} ,j\in \mN\}$.  This is an $\al$-permanental sequence with kernel,	 
   \begin{equation}
 \wt \VV_{t_{j},t_{k}}=  \VV_{t_{j},t_{k}}+f_{t_{k}},\qquad j,k\in\mN. \label{7.85q}
   \end{equation}
 It follows from (\ref{6.104mm}) that for an increasing sequence $\{t_{j}\}$,  
 \begin{equation}
 \wt \VV_{t_{j},t_{k}}=c_{1}^{2}( t_{j}\wedge t_{k})+O(1)+f_{t_{k}},\qquad  j,k\in\mN.\label{5.144mm}
   \end{equation}
   Set  
   \begin{equation}
      \wh\VV_{t_{j},t_{k}}=\frac{   \wt \VV_{t_{j},t_{k}}}{  ( \wt \VV_{t_{j},t_{j}})^{1/2} (\wt \VV_{t_{k},t_{k}})^{1/2}}.
   \end{equation}
    For $t_{j}\le t_{k}$ we have,
  \begin{equation}
    \wh\VV_{t_{j},t_{k}}   +   \wh\VV_{t_{k},t_{j}}=\frac{2t_{j}+O(1)+f_{t_{k}}+f_{t_{j}}}{(t_{j} t_{k})^{1/2}} .
   \end{equation}
Using the hypothesis that $f_{j}=o(j^{1/2})$  we see that for $t_{j}\le t_{k}$, 
 \begin{equation}
      \wh\VV_{t_{j},t_{k}}   + \wh\VV_{t_{k},t_{j}}= 2\(\frac{t_{j}}{t_{k}}\)^{1/2}  +o(1),\qquad\mbox{as $t_{j}\to \ff$}. \end{equation}
	In particular if $t_{j}=\th^{j}$ for some $\th>1$, for all $j\in\mN$, we have
 \begin{equation}
     \wh\VV_{t_{j},t_{k}}   +  \wh\VV_{t_{k},t_{j}}=  2\th^{-|k-j|/2} +o(1),\qquad\mbox{as $j,k\to \ff$}.\label{6.112}, \end{equation}	
Also, it is easy to see that, 	
 \begin{equation}
    \wh\VV_{\th^{j},\th^ {k}}   -  \wh\VV_{{\th^{k},\th^ {j}} } =  o(1),\qquad\mbox{as $j,k\to \ff$}\label{6.113}.
      \end{equation} 
 
The estimates  in (\ref{6.112})  and (\ref{6.113}) enable us to show that the hypotheses in \cite[Lemma 7.1] {MRejp} are satisfied.
Therefore, by taking $\th$  sufficiently large we have that any $\ep>0$,  
\begin{equation}
   \lim\sup_{j\to \ff}\frac{\wt \YY_{\al,\th^{j}} }{\th^{j} \log j}\ge 1-\ep.
   \end{equation}
This gives the lower bound in (\ref{7.86}) for all $\al>0$.\qed

Extending  the genealizaton of first order linear regressions in (\ref{auto.1qa}), we generalize the class of higher order Gaussian autoregressive sequences and find their covariances. In the beginning of this section we 
 consider  a class of  $k-$th  order autoregressive Gaussian sequences, $\wt \xi=\{\wt \xi_{n},n\in\mN \}$,  for $k\ge 2$.   
Let  $g_{1}, g_{2},\ldots$ be independent standard normal random variables and let  $ p_{i}>0$, $i=1,\ldots,k$, with  $\sum_{l=1}^{k}p_{l}\le 1$.    We define the Gaussian sequence  $\ov \xi=\{\ov \xi_{n},n\in\mN \}$  by,
\begin{equation}
\ov \xi_{1}=\frac{ g_{1}}{a}, \hspace{.2 in}\mbox{and}\quad\ov \xi_{n}=\sum_{l=1}^{k}p_{l}\ov\xi_{n-l}+g_{n},\hspace{.2 in}n\geq 2,\label{ark.1qq}
\end{equation}
where $\ov \xi_{i}=0$ for  all $i\le 0$ and    $a\ne 0$.  

\begin{lemma}
\bea
\VV^{[a^{2}]}_{m,n}\,:= E(\ov\xi_{m}\ov\xi_{n})&=&\frac{ 1-a^{2}}{a^{2} }    \phi_{m}\phi_{n}+E\(\wt\xi_{m}\wt\xi_{n}\)\label{covov.1}\\
&=&\frac{ 1-a^{2}}{a^{2} } \VV_{m,1}\VV_{1,n}+\VV_{m,n}.\nn
\eea
Furthermore,
  for all $ j\in\mN$,
\begin{equation}
\lim_{j\to\ff} \frac{\VV^{[a^{2}]}_{j,j}}{\VV_{j,j}}=1  .\label{8.3mm}
   \end{equation}
\end{lemma}

   \Proof
Generalizing (\ref{ark.2}) in Lemma \ref{lem-8.1mm}     we have,
\begin{equation}
\ov\xi_{n}=\phi_{n}\frac{g_{1}}{a}+\sum_{j=2}^{n} \phi_{n+1-j}g_{j}, \label{ark.2qq}
 \end{equation}
 where the $\phi_{n}$ are defined in (\ref{7.2mm}) for $\wt\xi_{n}$, not $\ov\xi_{n}$.
 The only difference between this and (\ref{ark.2}) is that $g_{1}$ is replaced by $g_{1}/a$. Therefore,  it follows from this and (\ref{7.27mm}) that,
\bea
 E(\ov\xi_{m}\ov\xi_{n})&=&\frac{\phi_{n}\phi_{m}}{a^{2}}+\sum_{j=2}^{m\wedge n}\phi_{m+1-j}\phi_{n+1-j}  \label{7.27mmqq}\\
 &=& \frac{ 1-a^{2}}{a^{2} }    \phi_{m}\phi_{n}+E\(\wt\xi_{m}\wt\xi_{n}\). \nn
   \eea
 The last equation in (\ref{covov.1}) follows from (\ref{7.27mm}).

 To obtain  (\ref{8.3mm}) we note that by (\ref{covov.1}),  
 \be 
  \VV^{[a^{2}]}_{j,j} =\VV_{j,j}+\frac{ 1-a^{2}}{a^{2}}   \phi ^{2}_{j}. \ee 
If $\sum_{j=1}^{k}p_{j}<1$ it follows from 
  (\ref{6.58nn}) that  $\{\phi_{j}\}\in \ell_{1.}$ This gives (\ref{8.3mm}) in this case. When $\sum_{j=1}^{k}p_{j}=1$ it follows from  (\ref{6.59ww}) and (\ref{get.10ii}) that $\phi_{j}=c_{1}+\psi_{j}$ where $\{\psi_{j}\}\in \ell_{1}$.
Since $\lim_{j\to\ff }\VV_{j,j}=\ff$ in this case we also get (\ref{8.3mm}). 
\qed

  We now show that $\VV^{[a^{2}]}_{m,n}$ is the potential density of a transient Markov chain. For the reason given in Remark \ref{rem-dec} we assume   that   $p_{i}\downarrow$. 

 \medskip		Consider the  matrix $ A$ defined in Lemma \ref{lem-7.5mm}. 
We generalize this matrix 
 by replacing $A_{1,1}=1+\sum_{i=1}^{k}p ^{2}_{i}$ by $a^{2}+\sum_{i=1}^{k}p ^{2}_{i}$.  Denote the generalized   matrix by $A^{[a^{2}]}$. In this notation $A=A^{[1]}$.

  \begin{theorem}\label{theo-8.1mm} If
   \begin{equation}
  a^{2}\geq    \frac{1}{2}\(\sum_{i=1}^{k}p_{i} \(2-\sum_{i=1}^{k}p_{i}\)-\sum_{i=1}^{k}p_{i} ^{2}\),\label{8.17}
   \end{equation}
  then $|a|>0$ and $-A^{[a^{2}]}$ is the $Q$-matrix of a transient Markov chain $\YY^{[a^{2}]}$ with potential density $\{ \VV^{[a^{2}]}_{m,n};\,m,n   \in\mN\}$.

 \end{theorem}
 
 \Proof    We show in    Lemma \ref{lem-6.5mm}, (\ref{ark.15v}), that the $m-$th row sums of the $-A^{[a^{2}]}$, $2\le m\le k$ are strictly greater than 0. Therefore, to see that  $-A^{[a^{2}]}$ 
    is a Q-matrix of a transient Markov chain, it suffices to check that the first row sum of $A^{[a^{2}]}$
 is greater than or equal to 0. We write this row sum as, 
   \begin{equation}
  \sum_{j=1}^{\ff}  A^{[a^{2}]}_{1,j}=a^{2}+\sum_{i=1}^{k}p_{i}^{2}+\ga, \label{ark.15sa2}
   \end{equation} 
   where $\ga$ is the sum of all terms to the right of the diagonal. It follows from  (\ref{ark.15s}) that,
   \be 
1+\sum_{i=1}^{k}p_{i}^{2} +  2\ga  =\(1-\sum_{i=1}^{ k}   p_{i}\)^{2}. 
 \label{ark.15sa} 
  \ee
Therefore,
\begin{equation}
   \ga=-\frac{1}{2}\(1+\sum_{i=1}^{k}p_{i}^{2} -  \(1-\sum_{i=1}^{ k}   p_{i}\)^{2} \).
   \end{equation}
By (\ref{ark.15sa2}) 
we see that the first row sum of $A^{[a^{2}]}$ is  strictly greater than  zero if,
 \begin{equation}
   a^{2}\geq  -\ga-\sum_{i=1}^{k}p_{i}^{2},\label{5.161l}
   \end{equation}
   which gives (\ref{8.17}). 
   
  Note that
    \begin{equation}
   -\ga-\sum_{i=1}^{k}p_{i}^{2}=   \frac{1}{2}\(\sum_{i=1}^{k}p_{i} \(2-\sum_{i=1}^{k}p_{i}\)-\sum_{i=1}^{k}p_{i} ^{2}\).
   \end{equation}
  It is easy to see that unless $p_{1}=1$ the right-hand side of (\ref{5.161l}) is strictly greater than 0. Since this is not possible by hypothesis, we see that $|a|>0$. \qed

Assume that (\ref{8.17}) holds. As in the proof of Theorem \ref{lem-invA}, 
to show that $\VV^{[a^{2}]}$ is the potential density for the Markov chain with Q-matrix $-A^{[a^{2}]}$ it suffices to show that
  \begin{equation}
  \VV^{[a^{2}]}A^{[a^{2}]}=I.\label{8.17r}
  \end{equation} 
    Using   (\ref{covov.1})   we see that (\ref{8.17r}) can be writen as,
  \begin{equation}
  \sum_{j=1}^{\ff}\(\frac{ 1-a^{2}}{a^{2} }   \VV_{m,1}\VV_{1,j}+\VV_{m,j}\)\(A_{j, n}+     (a^{2}-1)\de_{1}(j)\de_{1}(n)\)=\de_{m,n}.\label{}
  \end{equation}
  Since $\VV A =I$  by  (\ref{matrixj.1}) and $ \VV_{1,1}=1$, we need only show that for all $m$,
  \begin{equation}
 \frac{ 1-a^{2}}{a^{2} }   \VV_{m,1}   \sum_{j=1}^{\ff}\VV_{1,j}A_{j, n}+ \frac{ 1-a^{2}}{a^{2} } (a^{2}-1)  \VV_{m,1}\de_{1}(n)+  (a^{2}-1)  \VV_{m,1}\de_{1}(n)=0,\nn
  \end{equation}
which follows easily  since $ \sum_{j=1}^{\ff}\VV_{1,j}A_{j, n}=\de_{1}(n)$.
   \qed

We use Theorem \ref{theo-8.1mm}      
to extend  Theorem \ref{theo-4.2mmfx} to potentials of the form  $\VV^{[a^{2}]}$.
  
 \begin{theorem}\label{theo-4.2mmfxq} 
 Suppose that  $a^{2}$  satisfies (\ref{8.17}). Then Theorem \ref{theo-4.2mmfx} holds with 
$\YY$ and    $\VV$ replaced by $\YY^{[a^2]}$ and    $\VV^{[a^2]}$. 
   \end{theorem}

 	\Proof   The analogue of (\ref{7.87}) follows from Theorem \ref{theo-1.8mm} as in the  proof of Theorem \ref{theo-4.2mmfx},   $(i)$.  We now verify that the conditions for  Theorem \ref{theo-1.8mm} are satisfied.    By (\ref{covov.1})  
		 \be 
 \VV^{[a^{2}]}_{i,j}= \VV _{i,j}+\frac{ 1-a^{2}}{a^{2} }   \phi_{i}\phi_{j}   \label{hars.1aw}  .
   \ee  
Therefore, by (\ref{6.119}) and the fact that   $ \phi_{i} \le 1$ for all $i\in\mN,$
 \bea
   \sum_{j=1}^{\ff} \VV^{[a^{2}]}_{i,j}&=&   \sum_{j=1}^{\ff}  \VV _{i,j}+\frac{ 1-a^{2}}{a^{2} }     \phi_{i}  \sum_{j=1}^{\ff} \phi_{j}\\
   & \le&\nn 2 \|\phi\|^{2}_{1}+\frac{ 1-a^{2}}{a^{2} } \|\phi\|^{2}_{1}=\frac{ 1+a^{2}}{a^{2} }  \|\phi\|^{2}_{1}.
   \eea 
 Therefore, $\VV^{[a^{2}]}$ satisfies the second condition in (\ref{1.39}).  
 
Since $-A^{[a^{2}]}$ 
    is a Q-matrix,  $\VV^{[ a^{2}]}_{j,j}>0$ for each $j\in\mN$. In addition it follows from   (\ref{6.119}) and (\ref{8.3mm}) that 
\begin{equation}
\lim_{j\to\ff} \VV^{[a^{2}]}_{j,j}= c^{*}.\label{8.3nn}
   \end{equation}
 Therefore, $\VV^{[a^{2}]}$ also satisfies the first condition in (\ref{1.39}).  
  Using (\ref{8.3nn}) and Theorem \ref{theo-1.8mm} we get the analogue (\ref{7.87}). 

 	\medskip	The proof of the analogue of (\ref{7.86})  follows from   a slight generalization of the proof of Theorem \ref{theo-4.2mmfx}, $(ii)$.  We find an estimate for the row sums of $ \(  \VV^{[a^{2}]}(l,n)\)^{-1}$. Consider the terms defined in (\ref{5.90})--(\ref{5.92}) but with $\xi$ replaced by $\ov\xi$ defined in (\ref{ark.1qq}).
Lemmas \ref{lem-5.8} and  \ref{lem-6.8mm}  continue to hold with this substitution. The next lemma gives an analogue of Lemma \ref{lem-6.10mm}.

	\begin{lemma}\label{lem-6.10nn} When $\sum_{l=1}^{k}p_{l}=1$,
	  \begin{equation}
    E\(\ov \xi_{m}\ov \xi_{n}\)=c_{1}^{2}(m\wedge n)+a'_{n,m},\label{6.104nn}
   \end{equation}
where $|a'_{m,n}|\le D'<\ff$, for all $m,n\in\mN$.	
	
\end{lemma}	 

 \Proof    This follows immediately from (\ref{7.27mmqq}), (\ref{6.104mm}) and then   the fact that $ \phi_{i} \le 1$ for all $i\in\mN $.\qed
    
\noindent{\bf Proof of Theorem \ref{theo-4.2mmfxq} continued}    Using Lemma \ref{lem-6.10nn}  and following the proof of Lemma \ref{lem-7.8} we see that 
    (\ref{ntok.2}) holds for $\VV^{[a^{2}]}$. Similarly, Lemmas \ref{lem-4.3ja}
    and \ref{lem-6.11}  hold  for $\VV^{[a^{2}]}$. Consequently the proof of  the analogue (\ref{7.86}) follows immediately from the proof of  Theorem \ref{theo-4.2mmfx}.\qed  
    
  Theorem \ref{theo-rest}      also holds 
for potentials of the form  $\VV^{[a^{2}]}$.
    
   \begin{theorem} \label{theo-restaa} Under the hypotheses of Theorem \ref{theo-4.2mmfxq} assume in addition that $f_{j}=o(j^{1/2})$ as $j\to\ff$. Then   the analogue (\ref{7.86}) holds for all $\al>0$.
 \end{theorem}
 
  \Proof  This is immediate since Lemma \ref{lem-6.10nn} gives (\ref{5.144mm}).\qed

\begin{remark} {\rm Similar to what we point out in Remark \ref{rem-5.2} the condition in   (\ref{8.17}) is necessary for $\VV^{[a^{2}]}$ to be the potential of a Markov chain whereas (\ref{covov.1}) holds for all $a\ne 0$.    }\end{remark}

  	 	 \section{Proof of Theorem \ref{theo-mchain}}\label{sec-thm1.2}

Let $H=\{H_{j,k};j,k=1,\ldots,n\}$ be an $n\times n$ matrix with positive entries. We define,
\begin{equation}
H_{Sym}=\{(H_{i,j}H_{j,i})^{1/2}\}_{i,j=1}^{n} .   \end{equation}
Let  $K$ be an $n\times n$ inverse $M$-matrix and let $A=K^{-1}$.  We define  
\begin{equation}
A_{sym}= \left\{
\begin{array}{cl}
A_{j,j}&j=1,\ldots,n\\
-(A_{i,j}A_{j,i})^{1/2}&i,j=1,\ldots,n, i\ne j
\end{array}  \right.,\label{2.15}
\end{equation}
and 
\begin{equation}
  K_{isymi}= (A_{sym})^{-1}.\label{con7}
\end{equation}
(The notation $isymi$ stands for: take the inverse, symmetrize and take the inverse again.)    Obviously, when $K$ is symmetric, $A_{sym}=A$ and $ K_{isymi}=K$, but when $K$ is not symmetric, $ K_{isymi}\ne K$.

\begin{lemma}\label{lem-2.2mm}   The matrix     $K_{ isymi}$ is an inverse $M$-matrix and, consequently   is the kernel of $\al$-permanental random variables.
\end{lemma}

\Proof The matrix  $A=K^{-1}$ is a non-singular $M$-matrix. Therefore, by \cite[Lemma 3.3]{MRHD},  $A_{sym}$ is a non-singular  $M$-matrix.   We   denote  its inverse by $K_{isymi}$. 
The fact that  $K_{ isymi}$ is  the kernel of $\al$-permanental random variables  follows from  \cite[Lemma 4.2]{EK}.\qed

\medskip	
Theorem \ref{theo-mchain} is an application of   the next lemma which is    \cite[Corollary 3.1]{MRHD}.

	 \begin{lemma}\label{cor-7.1}  For any $\al>0$ let $\wt X_\al =(\wt X_{\al,0},\wt X_{\al,1},\ldots,\wt X_{\al,n})$ be an $\al$-permanental random variable with  kernel $K(n+1)$ that is an inverse   $M$-matrix and set $A(n+1)=K(n+1)^{-1}$. Let $\wt  Y_\al =(\wt Y_{\al,0},\wt Y_{\al,1},\ldots,\wt Y_{\al,n})$   be the $\al$-permanental random  variable determined by   $K(n+1)_{isymi}$.  Then for all functions $g $ of $\wt X_\al(n+1)$ and $\wt Y_\al(n+1)$ and sets $\BB$ in the range of $g $,
\bea
\lefteqn{ \frac{|A(n+1)|^{\al}}{|  A(n+1)_{sym}|^{\al}}  P\(g (\wt Y_\al(n+1))\in \BB \)  \le  P\( g (\wt   X_\al(n+1))\in\BB \)\hspace{.1in}\label{3.16w}}\\
& & \le \(1- \frac{|A(n+1)|^{\al}}{|  A(n+1)_{sym}|^{\al}} \) +\frac{|A(n+1)|^{\al}}{|  A(n+1)_{sym}|^{\al}}P\(g ( \wt Y_\al(n+1))\in\BB \).\nn
\eea
\el
It is clear that for this lemma to be useful we would like to have ${|A(n+1)|^{\al}}/\newline	{|  A(n+1)_{sym}|^{\al}} $  close to 1. 

\medskip	To obtain limit theorems we apply this lemma to sequences $\wt X_{\al}(l,n+1)=(\wt X_{\al,l },\wt X_{\al,l+1},\ldots,\wt X_{\al,l+n})$ 
with kernels $K(l,n+1)$ and consider
\begin{equation}
\nu_{l,n}:={ |A(l,n+1)_{sym} | \over  |A(l,n+1)  |}, \label{desc.12}
\end{equation}
where $A(l,n+1)=(K(l,n+1))^{-1}$. ($\wt Y_{\al}(l,n+1)$ is the $\al$-permanental random  variable determined by   $K(l,n+1)_{isymi}$.)

Here is how we obtain the matrices  $K(l,n+1)$. We start with a transient symmetric Borel right process, say $X$,   with state space $ \mN$,   and potential density $U =\{U_{j,k}\}_{j,k=1}^{\ff}$.  Then by    \cite[Lemma A.1]{MRnec},    
\be
U(l,n)=\{U_{l+j,l+k}\}_{j,k=1}^{n},\label{posrowsum}
\ee is the potential density of a transient  symmetric Borel right process,   say $\wh X$ on    $\{1,\ldots, n\}$. This implies that $U(l,n)$ is3
a symmetric inverse $M$ matrix with  positive row sums,  i.e., $\sum_{ k=1}^{n}(U(l,n))^{j,k}\geq 0$, for all $1\leq j\leq n$.

\medskip	Let    $f=\{f_{n}\}_{n\in\mN}$ be an excessive function with respect to $X$.   It follows from Theorem \ref{theo-borelN} that,  
\begin{equation}
 \wt U(l,n)=  \{U_{l+j,l+k}+f_{l+k}\}_{j,k=1}^{n}
   \end{equation}
is the kernel of an $\al$-permanental vector.
 	We define $K(l,n+1)$ to be an extension of $ \wt U(l.n)$ in the following way: 
\bea
\begin{array}{ lcll }
K(l,n+1)_{j,0}&=&1,&    j=0,\ldots,n, \label{10.6mm}\\
K(l,n+1)_{0,k}&=& f_{l+k},&  k=1,\ldots,n, \\
K(l,n+1)_{j,k}&=&U_{l+j,l+k}+f_{l+k},\quad& j,k=1,\ldots,n. 
\end{array}
\eea
Written out this is,
 \bea \hspace{-.3 in}K(l,n+1)& =&\left (
\begin{array}{ cccc } 1 &f_{l+1}&\ldots&f_{l+n}  \\
1   &U_{l+1,1}+f_{l+1}&\ldots&U_{l+1,n}+f_{l+n}  \\
\vdots& \vdots &\ddots &\vdots  \\
1   &U_{l+n,1}+f_{l+1}&\ldots&U_{l+n,n}+f_{l+n}  
\end{array}\right ).\label{19.39}
  \eea

  It is clear from (\ref{19.39}), by subtracting the first row from all other rows, that,
  \begin{equation}
 |K(l,n+1) |=  |U(l,n)|.\label{19.21e}
  \end{equation}
Therefore $K(l,n+1) $ is invertible. Let  $A(l,n+1)=K(l,n+1)^{-1}$.  
By multiplying the following matrix on the right by $K(l,n+1)$ one can check that,   

\[
\hspace{-4.3in}A(l,n+1) =
\]\be   \left(\!\!\!
\begin{array}{ cccc  }  1+\rho_{l,n}&-\sum_{ j=1}^{n}(U(l,n))^{j,1} f_{l+j}  &\dots&-\sum_{ j=1}^{n}(U(l,n))^{j,n} f_{l+j}  \\
-\sum_{ k=1}^{n}(U(l,n))^{1,k}& U(l,n)^{1 ,1 } & \dots &  U(l,n)^{ 1 ,n }    \\
\vdots&\vdots& \ddots&\vdots  \\
- \sum_{ k=1}^{n}(U(l,n))^{n,k} & U(l,n)^{n ,1 }&  \dots & U(l,n)^{n ,n }   \end{array}\!\!\!\right)  \label{19.40}
  \ee 
  where 
\begin{equation}
\rho_{l,n}=\sum_{ j,k=1}^{n}(U(l,n))^ {j,k}f_{l+k} .\label{condrealww}
\end{equation}

 	Note that all the row sums of $A(l,n+1)$ are equal to $0$, except for the first row sum which is equal to  $1$. Also the terms $U(l,n)^{j,k }$, $j,k=1,\ldots,n$, $j\ne k$ are negative because $U(l,n)$ is an  inverse $M$ matrix. Therefore, to show that $A(l,n+1)$ is an  M-matrix with positive row sums we need only check that 
\begin{equation}
\sum_{ j=1}^{n}(U(l,n))^{j,k} f_{l+j}\geq 0,\qquad \forall\,1\leq k\leq n.\label{localexc.1}
\end{equation}

We first consider the case in which,  
\begin{equation}
f =Uh,\qquad h\in l^{+}_{1}.\label{localexc.2}
\end{equation}
We point out in the second paragraph after the statement of Theorem \ref{theo-borelN} that in this case $f_{j}<\ff$, for all $j\in\mN.$

\medskip	
It follows from  \cite[Theorem 6.1]{MRejp}  applied to the transient symmetric Borel right process $Z$, with state space $\mN$ and potential $f$ in (\ref{localexc.2}) that we can obtain a transient symmetric Borel right process $\wt Z$, with state space $\mN\cup *$, where $*$ is an isolated point, such that $\wt Z$ has potential densities
\bea
   \wt U_{j,k}&=&U_{j,k}+f_{k},\qquad j,k\in\mN\\
     \wt U_{*,k}&=& f_{k},\quad\mbox{and}\quad    \wt U_{j,*}= \wt U_{*,*}=1\nn.
   \eea
It then follows from    \cite[Lemma A.1]{MRnec} that $K(l,n+1)$, defined in (\ref{10.6mm}), is  invertible and its inverse, $A(l+n)$ is a nonsingular $M$ matrix,  
 so   (\ref{localexc.1}) holds.   The inequality in (\ref{localexc.1}) can be extended to hold for all excessive functions because any  excessive function  is the  increasing limit of   potentials $\{f^{\(m\)}\}$ of the form (\ref{localexc.2}). (See  the proof of \cite[Theorem 1.11]{MRejp}.)  

\begin{remark} {\rm  The reader may wonder why we work with $K(l,n+1)$   instead of  simply $\{U_{l+j,1}+f_{l+k}\}_{j,k=1}^{n}$.  It is because it is easy to find $(K(l,n+1))^{-1}$ and it turns out to be a simple modification of $U(l,n)^{-1}$.  This is not the case for the inverse of $\{U_{l+j,1}+f_{l+k}\}_{j,k=1}^{n}$. 
 }\end{remark}	

	The next lemma is the critical estimate in the proof of Theorem \ref{theo-mchain}.

 \begin{lemma}\label{lem-10.3mm} For the   matrices $A(l,n+1)$ and $A(l,n+1)_{sym}$,
\be
1 \le \nu_{l,n} \le  1+\rho_{l,n}.\label{10.11}
\ee
\end{lemma}

\Proof 
  It follows from   (\ref{19.21e}) that
  \begin{equation}
 |A(l,n+1)|= |(U(l,n))^{-1}|.\label{19.21f}
  \end{equation} 
Also, since  $U$ is symmetric,
\be  A(l,n+1)_{sym} =\left (
\begin{array}{ cccc  }   1+\rho_{l,n}&-m(l,n)_{1 } &\dots&-m(l,n)_{n } \\
- m(l,n)_{1 }& U(l,n)^{1 ,1 } & \dots &  U(l,n)^{ 1 ,n }   \\
\vdots&\vdots& \ddots&\vdots  \\
- m(l,n)_{n } & U(l,n)^{n ,1 }&  \dots & U(l,n)^{n ,n }  \end{array}\right ) ,\label{19.41}
  \ee 
where 
\be
m(l,n)_{k }=\(c(l,n)_{k }r(l,n)_{k }\)^{1/2},
\ee 
and
\begin{equation}
c(l,n)_{k  }=\sum_{ j=1}^{n}(U(l,n))^{j,k} f_{l+j},\quad\mbox{and}\quad r(l,n)_{k }=\sum_{ j=1}^{n}(U(l,n))^{k,j}.\label{19.42}
\end{equation}
We  write this in block form,   
\be A(l,n+1)_{sym} =\left (
\begin{array}{ ccccc }(1+\rho_{l,n}) &-{\bf  m}(l,n)\\
-{\bf  m}(l,n)^{T}&U(l,n)^{-1}  \end{array}\right ),\label{19.43}
  \ee 
where ${\bf  m}(l,n)=(m(l,n)_{1 },\ldots,m(l,n)_{n })$.  Therefore,
\begin{equation}
|A(l,n+1)_{sym} |=|U(l,n)^{-1} |\,\,\( (1+\rho_{l,n})  -{\bf  m}(l,n)U(l,n){\bf  m}(l,n)^{T}\).\label{19.43a}
\end{equation}
(See, e.g.,  \cite[Appendix B]{DMM}.)

Using this and   (\ref{19.21f}) we see that  
\begin{equation}
\nu_{l,n}= (1+\rho_{l,n}) - {\bf  m}(l,n)U(l,n){\bf  m}(l,n)^{T}.\label{19.43c}
\end{equation}
It follows from    \cite[Lemma 3.3]{MRHD} that $\nu_{l,n}\ge 1$. Furthermore, since $U(l,n)$ is positive, ${\bf  m}(l,n)U(l,n){\bf  m}(l,n)^{T}\geq 0$. This gives (\ref{10.11}).\qed

The next lemma gives another critical estimate. Recall that      $K(l,n+1)_{isymi}$ is defined to be $(A(l,n+1))^{-1}_{sym})^{-1}$. It is an $(n+1)\times (n+1)$ matrix indexed by $j,k=0,\ldots,n$. In the next lemma we consider the $n\times n$ matrix $\{ K_{isymi}(l,n+1)\}_{j,k=1}^{n}$.

\bl\label{lem-10.4mm}  
\begin{equation}
\{ K_{isymi}(l,n+1)\}_{j,k=1}^{n}  =\{U(l,n)_{j ,k} +a(l,n)_{j } a(l,n)_{k }\}_{j,k=1}^{n}\label{22.30aaj}
\end{equation}
where
\begin{equation}
a(l,n)_{j} =\nu^{-1/2}_{l,n}\({\bf  m}(l,n) U(l,n)\)_{j}\leq f^{1/2}_{ l+j},\qquad  1\le j\le n.\label{abound}
\end{equation}    
\el

	\Proof  
By (\ref{19.43}) and the formula for  the inverse of  $A(l,n+1)_{sym}$ written as a block matrix; (see, e.g., \cite[Appendix B]{DMM}), we have,
\bea \lefteqn{K(l,n+1)_{isymi} \label{19.44}}\\
&& =\left (
\begin{array}{ ccccc } \nu^{-1}_{l,n} &\nu^{-1}_{l,n} {\bf  m}(l,n)U(l,n) \\
\nu^{-1}_{l,n} U(l,n){\bf  m}(l,n)^{T}  &U(l,n)+\nu^{-1}_{l,n} U(l,n){\bf  m}(l,n)^{T}{\bf  m}(l,n)U(l,n)  \end{array}\right ).\nn
  \eea 
%where 
%\[\ga_{l,n}= \( (1+\rho_{l,n})  -M(l,n)G(l,n)M(l,n)^{T}\)^{-1}= \nu_{l,n}.\] 
%Hence by (\ref{ga.1a})
%\begin{equation}
%{1 \over 1+\rho_{l,n}}\leq  \ga_{l,n}\leq 1.\label{ga.1h}
%\end{equation}
Note that for $i,j=1\ldots,n$,
\be \(U(l,n)+\nu^{-1}_{l,n} U(l,n){\bf  m}(l,n)^{T}{\bf  m}(l,n)U(l,n) \)_{i,j}\label{19.46} = U(l,n)_{i,j} +a(l,n)_{i} a(l,n)_{j}. 
\ee
Using the fact that $U(l,n)\geq 0$, we see that,
\bea
\({\bf  m}(l,n) U(l,n)\)_{j}&=&\sum_{i=1}^{n}m(l,n)_{i}U(l,n)_{i,j}=\sum_{i=1}^{n}\({c(l,n)_{i}r(l,n)_{i}}\)^{1/2}\,\,U(l,n)_{i,j}\nn \\ 
&  \leq &\(\sum_{i=1}^{n} c(l,n)_{i} U(l,n)_{i,j}\)^{1/2}\(\sum_{i=1}^{n} r(l,n)_{i}  U(l,n)_{i,j}\)^{1/2}.\nn\\
\label{abound.11}
\eea
Furthermore,
\bea
\sum_{i=1}^{n} c(l,n)_{i} U(l,n)_{i,j}&=&   \sum_{i=1}^{n} \sum_{ k=1}^{n}(U(l,n))^{k,i} f_{l+k} U(l,n)_{i,j}\label{abound.12}  \\
&=&     \sum_{ k=1}^{n}f_{l+k}  \sum_{i=1}^{n}(U(l,n))^{k,i}  U(l,n)_{i,j}\nn\\
&=&     \sum_{ k=1}^{n}f_{l+k}  \de_{k,j}=f_{l+j}\nn,
\eea
and, similarly,
\be 
\sum_{i=1}^{n} r(l,n)_{i} U(l,n)_{i,j}=1. \label{abound.13} 
\ee
Therefore,
\begin{equation}
\({\bf  m}(l,n) U(l,n)\)_{j}\le f^{1/2}_{l+j}.\label{abound.14} 
\end{equation}
Using this and (\ref{10.11}) we get (\ref{abound}).\qed

We can now give a concrete corollary of  Lemma \ref{cor-7.1}.

		 	 \begin{theorem}\label{theo-10.1mm}         For any $\al>0$, let  $ \wt X_{\al}(l,n )=(   \wt X_{\al,l+1},\ldots,  \wt X_{\al,l+n})$ be an $\al$-permanental random  variable determined by  the kernel 
\begin{equation}
\{U(l,n)_{j,k}+  f_{l+k}\}_{j,k=1}^{n} .
\end{equation}
Let   $   \wt  Y_\al (l,n)=(  \wt   Y_{\al,l+1},\ldots,   \wt  Y_{\al,l+n})$   be an $\al$-permanental random  variable determined by  the symmetric kernel,
\begin{equation}
\{ U(l,n)_{j ,k} +a(l,n)_{j } a(l,n)_{k }\}_{j,k=1}^{n} , \label{10.28}
\end{equation}
where    $ a(l,n)_{j},\;j=1,\ldots,n$, is given in (\ref{abound}).

Suppose that
\begin{equation}
 \rho_{l,n}=\sum_{ j,k=1}^{n}(U(l,n))^{-1}_{j,k}f_{l+k}\le \de_{l},\quad \mbox{ where } \quad\de_{l}=o(l).
\end{equation}
Then for all functions $  g $ of   $\wt X_\al(l,n )$ and $ \wt Y_\al(l,n )$,  and sets $\BB$ in the range of $  g $, and all $l$ sufficiently large,
\bea
P\(  g (  \wt  Y_\al(l,n))\in \BB \)-2\al\de_{l} & \le  &P\(   g ( \wt X_\al(l,n ))\in\BB \)\hspace{.1in}\label{3.16wa} \\
& \le&   2\al\de_{l} + P\(  g (  \wt  Y_\al(l,n))\in\BB \).\nn
\eea
\et

\medskip	\Proof  This follows from Lemma \ref{cor-7.1} and Lemmas \ref{lem-10.3mm} and \ref{lem-10.4mm}, with $  K(l,n+1)$ as defined in (\ref{10.6mm}). However we take   $g$ in Lemma \ref{cor-7.1} restricted to $( \wt  Y_{\al,1},\ldots,  \wt Y_{\al,n})$  and $ ( \wt X_{\al,1},\ldots,\wt X_{\al,n})$.
We also use the inequality
\begin{equation}
\(\frac{1}{1+\rho_{l,n}}\)^{\al}>1-2\al\de_{l,} 
\end{equation}
all $l$ sufficiently large.\qed

	\medskip	\noindent{\bf Proof of Theorem \ref{theo-mchain}} This is a direct application of Theorem \ref{theo-10.1mm}. We continue with the notation in Theorem \ref{theo-10.1mm} but   initially  we restrict ourselves to the cases where  $\al=k/2$, for integers $k\ge 1.$ We use (\ref{3.16wa}) with the event 
\begin{equation}
\{ g (   \wt   Y_ {k/ 2}(l,n))\in\BB\}= \left\{\sup_{1\leq j\leq n}{ \wt Y_{k/2,l+j}\over \phi_{l+j}} \leq 1\right\},
\end{equation}
and similarly for   $\{\wt X_{k/2,l+j}\}_{j=1}^{n}$. We have that for all $l$ sufficiently large  and $M>0$,  
\bea
P\(     \sup_{1\leq j\leq n}{  \wt  Y_{k/2,l+j}\over \phi_{l+j}} \leq M  \)-k\de_{l} & \le  &P\(      \sup_{1\leq j\leq n}{ \wt X_{k/2,l+j}\over \phi_{l+j}} \leq M \)\hspace{.1in}\label{3.16wb} \\
& \le&    k\de_{l} + P\(  \sup_{1\leq j\leq n}{ \wt  Y_{k/2,l+j}\over \phi_{l+j}} \leq M  \).\nn
\eea
The key point here is that  
\begin{equation}
\{ \wt  Y_{k/2,l+j}\}_{j=1}^{n}\stl \left \{ \sum_{i=1}^{k} \frac{( \eta_{i,l+j}+a(l,n)_{j}\xi_{i})^{2}}{2}\right\}_{j=1}^{n},
\end{equation}
where $
\{ \eta_{i,l+j}+a(l,n)_{j}\xi_{i}\}_{j=1}^{n}
$, $i=1,\ldots,k$, are independent copies of  $\{ \eta_{ l+j}+a(l,n)_{j}\xi \}_{j=1}^{n}$.  This follows from the definition of permanental processes in  (\ref{int.1}).   

We write 
\bea
&& \hspace{-.2in} \sum_{i=1}^{k} ( \eta_{i,l+j}+a(l,n)_{j}\xi_{i})^{2} =   \sum_{i=1}^{k}   \eta^{2}_{i,l+j} +2 a(l,n)_{j}\sum_{i=1}^{k}   \eta _{i,l+j}\xi_{i}+a^{2}(l,n)_{j}\sum_{i=1}^{k}\xi_{i}^{2} \nn\\
&&\le \sum_{i=1}^{k}   \eta^{2}_{i,l+j} +2 f^{1/2}_{l+j}\(\sum_{i=1}^{k}   \eta _{i,l+j}^{2}\)^{1/2}\(\sum_{i=1}^{k}\xi_{i}^{2}\)^{1/2}+ f _{l+j}\sum_{i=1}^{k}\xi_{i}^{2},  
\eea
by (\ref{abound}). Therefore,
\bea
&& \sup_{1\leq j\leq n}{  \wt  Y_{k/2,l+j}\over 2\phi_{l+j}}\label{9.37mm}\\
&&\qquad\le  \sup_{1\leq j\leq n}{\sum_{i=1}^{k}   \eta^{2}_{i,l+j} \over 2\phi_{l+j}} +2\ep^{1/2}_{l}\(\rho_{k}\sup_{1\leq j\leq n}{\sum_{i=1}^{k}   \eta^{2}_{i,l+j} \over 2\phi_{l+j}}\)^{1/2}+\ep _{l}\chi_{k},\nn\\
& &\qquad :=  \sup_{1\leq j\leq n}{\sum_{i=1}^{k}   \eta^{2}_{i,l+j} \over 2\phi_{l+j}}+A_{l,n}+{B_{l}}\nn,
\eea
where $\ep_{l}=\sup_{j\ge 1}(f_{l+j}/\phi_{l+j})$, which by (\ref{1.5mm}), goes to zero as $l\to \ff$, and $\chi_{k}=\sum_{i=1}^{k}\xi_{i}^{2}$. 

Consider the first inequality in (\ref{3.16wb}) and take the limit as $n\to\ff$. For all $\ep>0$ we have,
\bea
&&  P\(      \sup_{1\leq j\leq \ff}\frac{  \wt X_{k/2,l+j}}{\phi_{l+j} } \leq 1+\ep \)\label{10.45a}\\
&&\qquad\nn\ge  P\(      \sup_{1\leq j\leq \ff}{\sum_{i=1}^{k}   \eta^{2}_{i,l+j} \over 2\phi_{l+j}} \leq 1+\ep-A_{ l,\ff}-B_{l } \)-k\de_{l}.
\eea
Similarly, it follows from the second inequality in (\ref{3.16wb}) and the analogue of (\ref{9.37mm}) for the lower bound,     that for all $\ep>0$ we have,
\bea
&&  P\(      \sup_{1\leq j\leq \ff}\frac{   \wt X_{k/2,l+j}}{\phi_{l+j} } \leq 1-\ep \)\label{10.45}\\
&&\qquad\nn\le  P\(      \sup_{1\leq j\leq \ff}{\sum_{i=1}^{k}   \eta^{2}_{i,l+j} \over 2\phi_{l+j}} \leq 1-\ep+A_{ l,\ff}+B_{l } \)+k\de_{l}.
\eea

	It follows from (\ref{121.4a}) and  Lemma \ref{lem-Rev} below that,
\begin{equation}
\varlimsup_{j\to\ff } {\sum_{i=1}^{k} \eta^{2}_{i,j} \over   2\,\phi _{j}}=1, \hspace{.2 in}a.s.\label{rev.2o}
\end{equation} 
Therefore, if we take the limits in  (\ref{10.45a}) and (\ref{10.45}) as $l\to\ff$   we get that for all $\ep'>0.$ 
\begin{equation}
1-\ep'\le    \varlimsup_{j\to\ff}{  \wt X_{k/2,j}\over \phi_{j }} \leq 1+\ep', \qquad a.s.    \label{10.47}\end{equation}
and since this holds for all $\ep>0$     we get,
\begin{equation}
    \varlimsup_{j\to\ff}{   \wt X_{k/2,j}\over \phi_{j }} = 1, \qquad a.s.    \label{10.48}\end{equation}

Now, suppose that $1/2\le \al\le k'$ for some integer $k'$. Since (\ref{10.48}) holds for $k=1/2$ and $k=k'$ we can use the property that $\al$-permanental processes are infinitely divisible and positive to see that (\ref{gena.ka}) holds.\qed

\bl\label{lem-Rev}

Let    $\{\eta_{j};j\in\mN\}$ be a Gaussian sequence and for each $i\in \mN$,    
let $\{\eta_{i,j}, j\in\mN\}$ be an independent copy of $\{\eta_{j};j\in\mN\}$.
Let $\{\phi_{j}\}$ be a   sequence     such that,   
\be
\varlimsup_{j\to \ff}\frac{  |\eta_{j}|}{( 2\,\phi_{j} )^{1/2}}=1 \qquad 
a.s.\label{rev.1},
\ee
then for any integer $k\ge 1$,
\begin{equation}
\varlimsup_{j\to\ff } {\sum_{i=1}^{k} \eta^{2}_{i,j} \over   2\,\phi _{j}}=1, \hspace{.2 in}a.s.\label{rev.2}
\end{equation} 
(This also holds for a Gaussian process $\{\eta_{t};t\in R^{+}\}.$  \el

\Proof We follow the proof of the law of the iterated logarithm for Brownian motion in \cite[Theorem 18.1]{Revesz}.  Clearly, we only need to prove the upper bound.

Fix $k$
and $\ep>0$ and $X_{j}=(\eta_{1,j},\ldots, \eta_{k,j}) $ and   $u$ be a unit vector in $R^{k}$.  By checking their covariances we see that
\be
\{(u\cdot X_{j}),\,\, j=1,\ldots \}\stackrel{law}{=}\{\eta_{j},\,\, j=1,\ldots \}.\label{9.13mm}
\ee
Therefore,  by (\ref{rev.1}),  
\be
\varlimsup_{j\to \ff}\frac{ |(u\cdot X_{j})|}{(2\,\phi_{j} )^{1/2}}=1 \qquad 
a.s.\label{rev.3} 
\ee
Note that $\| X_{j}\|_{2}=(\sum_{i=1}^{k} \eta^{2}_{i,j})^{1/2}$ and,
\be 
\varlimsup_{j\to \ff}\frac{ \| X_{j}\|_{2}}{(2\,\phi_{j} )^{1/2} }=\varlimsup_{j\to \ff}\sup_{\|u\|_{2}=1}\frac{ |(u\cdot X_{j})|}{(2\,\phi_{j} )^{1/2} }.
\label{rev.4}
\ee
For any $\ep>0$ we can find a finite  set of unit vectors $\UU(m)=(u_{1},\dots,u_{m}\}$  in $R^{k}$ with the property that for any unit vector $u$ in $R^{k}$, $\inf_{1\le l\le m} \|u-u_{l}\|_{1}\leq \ep/k$. 

Let $u$ be a unit vector in $R^{k}$.
For  all  $u_{l}\in\UU(m)$,
\bea
|(u\cdot X_{j})|&\le&|((u-u_{l})\cdot X_{j})|+|(u_{l}\cdot X_{j})|\\
&\le  &\|u-u_{l}\|_{1} \left |\frac{ (u-u_{l}) }{\|u-u_{l}\|_{1}}\cdot X_{j} \right |+\sup_{l=1,\ldots, M}|(u_{l}\cdot X_{j})|\nn\\
&\le  &\|u-u_{l}\|_{1} \left |\frac{ (u-u_{l}) }{\|u-u_{l}\|_{1}}\cdot X_{j} \right |+\sup_{l=1,\ldots, M}|(u_{l}\cdot X_{j})|\nn\\
&\le  &\|u-u_{l}\|_{1} \sum_{i=1}^{k}|\eta_{i,j}|+\sup_{l=1,\ldots, M}|(u_{l}\cdot X_{j})|\nn.
  \eea
	Since this holds for all $u_{l}\in\UU(m)$, we see that,
\bea
  |(u\cdot X_{j})|&\le&\min_{1\le l\le m} \|u-u_{l}\|_{1} \sum_{i=1}^{k}|\eta_{i,j}|+\sup_{l=1,\ldots, M}|(u_{l}\cdot X_{j})|\\
  &\le&\frac{\ep}{k} \sum_{i=1}^{k}|\eta_{i,j}|+\sup_{l=1,\ldots, M}|(u_{l}\cdot X_{j})|\nn.
  \eea
Consequently,
\be 
 \varlimsup_{j\to \ff}\sup_{\|u\|_{2}=1}\frac{ |(u\cdot X_{j})|}{(2\,\phi_{j} )^{1/2}}\le\varlimsup_{j\to \ff}\sup_{l=1,\ldots, M}\frac{ |u_{l}\cdot X_{j}|}{(2\,\phi_{j} )^{1/2} }+{\ep \over k}\sum_{i=1}^{k}\varlimsup_{j\to \ff}\frac{ |\eta_{i,j}|}{(2\,\phi_{j} )^{1/2} }  .\label{9.18}
  \ee 
It follows from  (\ref{rev.1}) that the last term is bounded by $\ep$. 

 \medskip	Let $\Om'$ be the event that equality holds in (\ref{rev.3})  with $u=u_{l}$ for all $l=1,\ldots, m$. 
It follows   that for any $\ep$ and any $\om\in \Om'$ we can find $j_{0}(\om)$ such that 
\begin{equation}
\frac{ |u_{l}\cdot X_{j}(\om)|}{(2\,\phi_{j} )^{1/2} }\leq 1+\ep,\qquad \forall j\geq j_{0}(\om) \mbox{ and all }1\le l\le m.
\end{equation}
Since $P\(\Om'\)=1$, it now follows from (\ref{9.18}) that 
\be 
 \varlimsup_{j\to \ff}\sup_{\|u\|_{2}=1}\frac{ |(u\cdot X_{j})|}{(2\,\phi_{j} )^{1/2}}\le 1+2\ep,\qquad a.s.\label{9.18a}
  \ee 
   Since this holds for all $\ep>0$,    
 the upper bound for (\ref{rev.2}) follows from (\ref{rev.4}).\qed

\section{Proof of Theorem \ref{theo-1.8mm}}\label{sec-8}

      Let  $M$  be an $ \mN\times  \mN$ matrix and consider the operator norm on $\ell_{\ff}\to\ell_{\ff}$,  
  \begin{equation}
\|M\|=\sup_{\|x\|_{\ff}\leq 1}  \|Mx\|_{\ff}=\sup_{j} \sum_{k} |M_{j,k}|.\label{unif.1}
  \end{equation}
  
\begin{lemma}\label{lem-11.8} Let $M=\{M_{j,k},j,k\in\mN\}$ be a     positive matrix and assume that   both   $\|M\|$ and $\|M^{T}\| <\ff$. Then for all $\ep>0$,  there exists a sequence $\{i_{n},n\in\mN\}$ such that $i_{n}\le n(\|M\|+\|M^{T}\|)/ \ep $, for all $n\in\mN$, and 
\begin{equation}
   M_{ i_{j},i_{k}}\le \ep,\qquad \forall\,j,k\in \mN,\hspace{.1in}j\ne k.\label{small.51}
   \end{equation}

\end{lemma}

 \Proof  Assume to begin that $M$ is symmetric.   Fix $\ep>0$,  and consider $\{M_{1,k}\} _{k=1}^{\ff}$. Not  more than $\|M\|/\ep$ of these terms can be greater than $\ep$. Let   $\{M_{1,k_{1}(  p_{i})}, i=1,\dots C_{1} \}$ denote the  terms in $\{M_{1,k}\} _{k=1}^{\ff}$ which are greater than $\ep$ and set
 $R_{1}=\{k_{1}(  p_{i}) , i=1,\dots C_{1}\}$. As we just pointed out  $ |R_{1}|\le { \|M\|}/{\ep}$. 
 
 Note that
  \begin{equation}
   M_{1,k}\le \ep\qquad \forall\,   k\in R^{ c}_{1}.
   \end{equation}
 Set $i_{1}=1$  and  
set $i_{2}$ equal to the smallest index in $R^{ c}_{1}$ that is greater than $i_{1}$. 

We repeat this procedure starting with  $M_{i_{2},k}$ with $k\in R^{ c}_{1} $ to get  $R_{2}$ where $ |R_{2}|\le { \|M\|}/{\ep}$ and,
  \begin{equation}
   M_{i_{2},k}\le \ep,\qquad \forall\,   k\in R^{ c}_{2}.
   \end{equation}
 Therefore, for $j=1,2$,
   \begin{equation}
   M_{i_{j},k}\le \ep\qquad \forall\,   k\in \(   R_{1}\cup R_{2}\)^{ c}=R^{ c}_{1}\cap R^{ c}_{2}.
   \end{equation}
We continue this procedure  setting $i_{3}$ equal to the smallest integer in $ \(   R_{1}\cup R_{2}\)^{ c}$ that is greater than $i_{2}$, and so on,  to get  $\{i_{n},n\in\mN\}$.
  This completes the proof  when $M$ is symmetric. 
  
    More generally, assume only that   both   $\|M\|$ and $\|M^{T}\| <\ff$.  We use a construction similar to the one above but we work alternately with both  $M$ and $M^{T}$. Therefore, we can obtain $\{i_{l}, l=1,\ldots, n\}$ and a set $S_{n}\subset \mN$ such that 
$|S_{n}|\le n(\|M\|+\|M^{T}\|)/ \ep $ and for $l=1,\ldots, n$,
   \begin{equation}
   M_{i_{l},k}\le \ep, \qquad \mbox{and}\qquad M_{k, i_{l}}\le \ep, \qquad\forall\,   k\in \(S_{n}\)^{ c}.
   \end{equation}
Choose $i_{n+1}$ equal to the smallest integer in $\(S_{n}\)^{ c}$ that is greater than $i_{n}$. 

We continue the  above procedure starting with  $M_{i_{n+1},k}$ and $M_{k, i_{n+1}}$, with $k\in S^{ c}_{n} $, to get  $R_{n+1}$ where $ |R_{n+1}|\le {(\|M\|+\|M^{T}\|)}/{\ep}$ and,
  \begin{equation}
   M_{i_{n+1},k}\le \ep, \qquad \mbox{and}\qquad   M_{k, i_ {n+1}}\le \ep,\qquad\forall\,   k\in R^{ c}_{n+1}.
   \end{equation}
 Therefore, for $j=1,\ldots, n+1$,
   \begin{equation}
   M_{i_{j},k}\le \ep,\qquad   M_{k, i_{j}}\le \ep\qquad\forall\,   k\in \(   S_{n}\cup R_{n+1}\)^{ c}=S_{n}^{ c}\cap R^{ c}_{n+1}.
   \end{equation}
This shows that   for all $\ep>0$,  there exists a sequence $\{i_{n},n\in\mN\}$ such that $i_{n}\le n(\|M\|+\|M^{T}\|)/ \ep $, for all $n\in\mN$, and in particular that (\ref{small.51})
holds.

 \qed

\noindent{\bf  Proof of  Theorem \ref{theo-1.8mm} }    It follows from (\ref{int.1}) that for all $j\in\mN$,  
\begin{equation}
   \frac{\wt X_{\al,j} }{\wt U_{j,j}}\stl \xi_{\al},
   \end{equation}
where $\xi_{\al}$ has probability density function ${x^{\al-1} e^{-x}}/{|\Ga(\al)}$. Using  the Borel-Cantelli Lemma, we get,
 \begin{equation}
   \limsup_{n\to\ff} \frac{\wt X_{\al, n} }{\wt U_{n,n}\log  {n}}\le 1    \qquad  {a.s.}\label{small.53mm}
   \end{equation} 
   This gives the upper bound in (\ref{121.4mb}) because  since  $f \in c^{+}_{0}$  and    $\inf   U_{n,n} >0$, we have
\begin{equation}
   \lim_{n\to\ff}\frac{  U_{n,n}}{\wt U_{n,n}}=1.\label{9.8mm}
   \end{equation}

\medskip	     To get the lower bound in  (\ref{121.4mb}) consider, 
\begin{equation}
\wh U_{j,k}=\frac{\wt U_{j,k}}{(\wt U_{j,j}\wt U_{k,k})^{1/2}}, \qquad j,k\in\mN .
   \end{equation}
  It follows from Lemma \ref{lem-11.8} that for all $\ep>0$  there exists a sequence $\{i_{n},n\in\mN\}$ with  
  \be
  i_{n}\le 2n\| U\|/ (\ep\de),\qquad \forall n\in\mN,\label{small.51j}
  \ee 
  such that,
  \begin{equation}
  U_{ i_{j},i_{k}}\le \frac{\ep\de}{2},\qquad \forall\,j,k\in \mN,\hspace{.1in}j\ne k,\label{small.51a}
   \end{equation} 
  Therefore, 
   \begin{equation}
\wh U_{i_{j},i_{k}}+\wh U_{i_{k},i_{j}}\le  \frac{2  U_{i_{j},i_{k}}+f(i_{k})+f(i_{j})}{(  U_{i_{j},i_{j}}  U_{i_{k},i_{k}})^{1/2}} \le   \ep+\frac{ f(i_{k})+f(i_{j})}{\de}.
    \end{equation}
   Using the fact that $f\in c^{+}_{0}$ we see that we can find an $n_{0}$ such that 
     \begin{equation}
\wh U_{i_{j},i_{k}}+\wh U_{i_{k},i_{j}}\le   2\ep,\qquad\forall j,k\ge n_{0}.    \end{equation}
Therefore, by  \cite[Lemma 7.1]{MRejp},  
   \begin{equation}
   \limsup_{n\to\ff} \frac{\wt X_{\al, i_{n}}}{ \wt U_{i_{n},i_{n}}  \log (n-n_{0})}\ge 1  -6\ep \qquad  {a.s.}\label{small.52},
   \end{equation} 
or, equivalently,
   \begin{equation}
   \limsup_{n\to\ff} \frac{\wt X_{\al, i_{n}}}{\wt U_{i_{n},i_{n}}  \log n  }\ge 1  -6\ep \qquad  {a.s.}\label{small.52a}
   \end{equation} 
Using (\ref{9.8mm}) and   (\ref{small.51j}), we get, 
  \begin{equation}
   \limsup_{n\to\ff} \frac{\wt X_{\al, i_{n}}}{  U_{i_{n},i_{n}}\log i_{n}}\ge 1  -6\ep \qquad  {a.s.}\label{small.53}
   \end{equation}
which gives (\ref{121.4mb}). \qed

  If  $f=Uh$, where $\|U\|<\ff$, and $h\in\ell_{1}$ then   since $\sum_{i}f_{i}=\sum_{i}\sum_{j}U_{i,j}h_{j}$, it follows using the symmetry of $U$ that $f\in\ell^{+}_{1}$ and consequently in $ c^{+}_{0}$. However, when $U$ has some regularity,  $f\in c^{+}_{0}$   if and only if $h\in c^{+}_{0}$.
  
 	 	\begin{lemma}\label{lem-8.2nn} Let $f=Uh$, where $\|U\|<\ff$, $h\in c_{0}^{+}$ and   there exists a   sequence  $\{k_{j}\}$, $\lim_{j\to\ff}k_{j}=\ff$, such that
\begin{equation}
   \lim_{j\to\ff}\sum_{k=1}^{k_{j}}U_{j,k}=0.\label{8.14}
   \end{equation}
Then $f\in c^{+}_{0}.$

  If  $\inf _{j}U_{j,j}>0$  then $f\in c^{+}_{0} $ implies that $h\in c^{+}_{0} $.
 \end{lemma}

\Proof The first statement  follows from the inequality,
\begin{equation}
f_{j} =  \sum_{k=1}^{\ff }U_{j,k} h_{k} \le \|h\|_{\ff}\sum_{k=1}^{k_{j} }U_{j,k  }+\|U\|\sup_{k\ge k_{j}}|h_{k}|.   \end{equation}
   The second statement is obvious, since,
   \begin{equation}
   f_{j}=\sum_{k=1}^{\ff}U_{j,k}h_{k}\ge U_{j,j}h_{j}.
   \end{equation}
\qed

 \noindent	{\bf  Proof of  Theorem \ref{theo-1.8gen} }     We show in \cite[Theorem 6.1]{MRejp}  that $\wt U$ is the kernel of an $\al$-permanental sequence.    Also, it follows from  (\ref{1.39gen}) that   $f\in \ell^{+}_{1}\subset c^{+}_{0}$ and that the hypotheses of Lemma \ref{lem-11.8} are satisfied. Consequently, the proof   follows as in the proof of Theorem \ref{theo-1.8mm}.\qed

 	\noindent{\bf  Proof of  Theorem \ref{theo-lev} }  
	The L\'evy process  $X$   is obtained by killing a  L\'evy process   say $\wh X$ on $\mathbb Z $   at the end of an independent exponential time with mean $1/\bb$.  Let  $\{\wh p_{t}(i,j);j,k\in\mathbb Z\}$ denote  the transition densities for $\wh X$
 and $\{ p_{t}(i,j);j,k\in\mathbb Z\}$   the transition densities for  $ X$. We have
 \begin{equation}
 p_{t}(i,j)=e^{-\bb t}\wh p_{t}(i,j).\label{lev.2}
 \end{equation}
Consequently,  
 \begin{equation}
 U_{j,k}=\int_{0}^{\ff} e^{-\bb t}\wh p_{t}(j,k)\,dt,\qquad \forall j,k\in \mathbb Z.\label{lev.3}
 \end{equation}
Since $\wh X$ is a Levy process we have
\begin{equation}
   \wh p_{t}(i,j)=\wh p_{t}(0,j-i):=\wh p_{t}(j-i).\label{lev.3a}
   \end{equation}
   Therefore, for all $j\in\mathbb Z,$
     \be
  U_{j, j}=U_{0,0} =\int_{0}^{\ff} e^{-\bb t}\wh p_{t}(0)\,dt=\int_{0}^{\ff}   p_{t}(0)\,dt.\label{10.29mm}
  \ee

  To see that   (\ref{lev.0}) is the kernel of an $\al$-permanental sequence we  first note that since $X$ is an exponentially killed  L\'evy process on $\mathbb Z $ with potential density   $U$, then   $\ov X=-X$ is also a  L\'evy process on $\mathbb Z $, the dual of $X$,  with transition densities
    \begin{equation}
\ov p_{t}(i,j)=  p_{t}(-i,-j):=p_{t}( i-j)=p_{t}(j, i), \label{lev.0j1}
  \end{equation}
  and   consequently, potential densities  
  \begin{equation}
\ov U _{i, j}=  U_{j,i},\qquad i,j\in\mathbb Z. \label{lev.0j2}
  \end{equation}

  The proof that   (\ref{lev.0}) is the kernel of an $\al$-permanental process for all functions $f$ that are finite excessive functions for $X$ proceeds in three steps.   
   
\medskip	 We first show that for any   $g=\{g_{k}\}$ where $g_{k}=\sum_{j=-\ff}^{\ff}\ov U _{k, j} h_{j}$, and  $h\in \ell^{+}_{1}(\mathbb Z) $,
      \begin{equation}
 U_{j,k}+g_{k},\qquad j,k\in \mathbb Z,\label{lev.0m1}
   \end{equation} 
 is the kernel of an $\al$-permanental process. To see this note that by (\ref{lev.0j2}),
 $g_{k}=\sum_{j=-\ff}^{\ff}\ov U _{k, j} h_{j} =\sum_{j=-\ff}^{\ff} h_{j}U _{j, k}  $. Therefore it follows from  \cite[Theorem 6.1]{MRejp} that (\ref{lev.0m1})   is the restriction to $ \mathbb Z\times  \mathbb Z$ of  the  potential densities of  a transient Borel right process
 $\wt X\!=\!
(\Om,  \FF_{t},\wt  X_t,\th_{t},\wt P^x
)$ with state space $ \mathbb Z\cup \{\ast\}$, where $\ast$ is  an isolated point.   Consequently, (\ref{lev.0m1}) is  the kernel of an $\al$-permanental sequence.

\medskip	We  show next that (\ref{lev.0m1}) is the kernel of an $\al$-permanental process for any    $g $ that  is a finite excessive function for    $\ov  X$. We use the following lemma  which is Lemma 6.2 in   \cite{MRejp}.
  \begin{lemma} \label{lem-6.1} Assume that   for each $n\in \mathbb N$,   $u^{ ( n)}(s,t),\,s,t\in S $, is the kernel of an $\al$-permanental process.  If
  $  u^{ ( n)}(s,t)\to   u(s,t)$ for all $s,t\in S$, then $  u(s,t) $ is the kernel of an $\al$-permanental process.
  \end{lemma}
  
  We now use arguments from  the proof of  \cite[Theorem 1.11]{MRejp}.  Consider  a general   function $g=\{g_{k}\}$ that is a finite excessive function for $\ov X$.  It follows from \cite[II, (2.19)]{BG} that  there exists a sequence of functions $h^{(n)}=\{h^{(n)}_{k}\}  \in \ell_{\ff}^{+}( \mathbb Z)$    such that $g^{(n)}$ defined by,
\begin{equation}
g^{(n)}_{k} =\sum_{j=-\ff}^{\ff}\ov U _{k, j} h^{(n)}_{j}\label{rp.11},
\end{equation}
is also in $\ell_{\ff}^{+}( \mathbb Z)$ and is such that for each $k\in\mN$, $g^{(n)}_{k}\uparrow g_{k}$.

 \medskip	 If   $h^{(n)}\in\ell_{1}^{+}$ then by the first step in this proof we have that   $\{ U_{j,k}+g^{(n)}_{k}),  j,k\in \mathbb Z \}$ are kernels of  $\al$-permanental   processes.  Consequently,  by Lemma \ref{lem-6.1}, (\ref{lev.0m1})  is the kernel  of  $\al$-permanental process.

\medskip	If 
$h^{(n)}\notin\ell_{1}^{+}$ we first consider  $h^{(n)}1_{[-m,m]}$  which clearly is in $\ell_{1}^{+}$ for each $m<\ff$. We then set
\be
g^{(n,m)}_{k}= \sum_{j=-\ff}^{\ff}\ov U _{k, j} h^{(n)}_{j}1_{\{-m\leq j\leq m\}}.
\ee Therefore, as in the previous paragraph, we have  that  $\{ U_{j,k}+g^{(n,m)}_{k},  j,k\in \mathbb Z \}$ is the kernel of an $\al$-permanental process. Taking the limit as $m\to \ff$, it follows from  Lemma \ref{lem-6.1} that  $\{ U_{j,k}+g^{(n)}_{k},  j,k\in \mathbb Z \}$  is the kernel of an $\al$-permanental process. Since $g^{(n)}_{k}\to g_{k}$ we  use Lemma \ref{lem-6.1} again to see that  (\ref{lev.0m1})   is the kernel  of  an $\al$-permanental process for all finite excessive functions $g$ for $\ov X$.

 \medskip	 The last step in the proof that   (\ref{lev.0}) is the kernel of an $\al$-permanental process is  to show that when  $f_{k}$ is  a finite excessive function for $ X$, then $f_{-k}$ is   a finite excessive function for  $\ov X$.
 To see this, note that    if  $f_{k}$ is  a finite excessive function for  $ X$,  
  then, by definition,   
 \begin{equation}
\sum_{k=-\ff}^{\ff} p_{t}( k-i)f_{k}=\sum_{k=-\ff}^{\ff} p_{t}(i,k)f_{k}\uparrow f_{i},\qquad \mbox{as }  t\downarrow0.
 \end{equation}
It follows from this that  as $t\downarrow 0$,  
  \begin{equation}
\sum_{k=-\ff}^{\ff} p_{t}(k,i)f_{-k}= \sum_{k=-\ff}^{\ff} p_{t}(-k-(-i))f_{-k}\uparrow f_{-i}.
  \end{equation}
Consequently $f_{-k}$ is  a finite excessive function  for $\ov X$.

\medskip	 This  completes the proof that (\ref{lev.0}) is the kernel of an $\al$-permanental process. Using this and the fact that $\lim_{k\to\ff}f_{-k}=0$, proceeding exactly as in the proof of Theorem \ref{theo-1.8mm}, we get the upper bound in (\ref{lev.1}).

 	  To obtain the lower bound in (\ref{lev.1}) we note that 
by  (\ref{lev.3}) (\ref{lev.3a}) and  Fubini's Theorem,
   \begin{equation}
\sum_{i=-\ff}^{\ff} U_{i, j}=\sum_{j=-\ff}^{\ff} U _{i, j}=\int_{0}^{\ff} e^{-\bb t} \,dt=\frac{1}{\bb}.\label{lev.4}
 \end{equation}
  Using this and (\ref{10.29mm}) we see that 
  the conditions in (\ref{1.39gen}) and Lemma \ref{lem-11.8} are all satisfied for $\{U_{j,k};j,k\in\mN\}$. Therefore, as in the proofs of Theorems \ref{theo-1.8mm} and \ref{theo-1.8gen}, the lower bound in (\ref{lev.1})  follows from  \cite[Lemma 7.1]{MRejp}.

 	\medskip	   To verify the last statement in this theorem we see from  the proof of Lemma \ref{lem-8.2nn}, that we need only show that 
   there exists a   sequence  $\{k_{j}\}$, $\lim_{j\to\ff}k_{j}=\ff$, such that   \begin{equation}
   \lim_{|j|\to\ff}\sum_{k=-k_{|j|}}^{k_{|j|}}U_{j,k}=0.\label{8.14o}
   \end{equation}
  We have,
  \begin{equation}
   \lim_{|j|\to\ff} \sum_{k=-k_{|j|}}^{k_{|j|}} U_{j,k}=    \lim_{|j|\to\ff} \sum_{k=-k_{|j|}}^{k_{|j|}} U_{0,-j+k}=   \lim_{|j|\to\ff}  \sum_{l=-j-k_{|j|}}^{-j+k_{|j|}} U_{0,l} . \end{equation}
It follows from (\ref{lev.4}) that this last term goes to zero when $k_{j}=j/2$.
 \qed

   \section{Uniform Markov chains}\label{sec-uniform}

\begin{lemma} \label{lem-9.1mm}
Let $X\!=\!
(\Om,  \FF_{t}, X_t,\th_{t},P^x
)$ be a transient  Borel right process with state space $ \mN$, finite Q-matrix  $Q$, and  strictly positive potential   densities  $U=\{U_{j,k}$, $j,k\in\mN\}$. Then,  
 \begin{equation}
-\de_{i,l}=\sum_{j=1}^{\ff} Q_{i,j}U_{j,l},\qquad  \mbox{for all $i, l\in\mN$}\label{inv.34a}.
\end{equation}

\end{lemma} 

\Proof    
Set $q(i)=-Q_{i,i}$. Without loss of generality we can take $q(i)>0$.  For any function $h$ we have,  
\be
U h(i) = {h(i)\over q(i)}+\sum_{j=1,j\not=i}^{\ff}  {Q_{i,j}\over q(i)}U\!h(j).\label{11.58}
\ee
To see this, let $\tau_{i}$ be the time of the first exit from state $i$ and  note that, 
 \be 
 U h(i)=E^{i}\(\int_{0}^{\tau_{i}}h(X_{t})\,dt\)+E^{X_{\tau_{i}}}\(\int_{0}^{\ff}h(X_{t})\,dt\).
\ee
Using the facts that the  exit time is an exponential random variable with expectation $1/q(i)$, and 
 the probability that upon exit the process jumps from $i$ to $j$ is $Q_{i,j}/q(i)$, we get  the two terms in (\ref{11.58}).

 It follows from (\ref{11.58}) that 
 \begin{equation}
  -h(i)=-q(i)Uh(i)+\sum_{j=1,j\not=i}^{\ff}  Q_{i,j} Uh(j)=\sum_{j =1}^{\ff}  Q_{i,j} U h(j).\label{inv.34x}
  \end{equation}
  Take $h(k)=\de_{l,k}$ and note that,
  \begin{equation}
   \sum_{j =1}^{\ff}  Q_{i,j} U h(j)=   \sum_{j =1}^{\ff}  Q_{i,j}   \sum_{k=1}^{\ff} U_{j,k} \de_{l,k}=\sum_{j=1}^{\ff} Q_{i,j}U_{j,l}.
   \end{equation}
  Therefore, by (\ref{inv.34x}),
    \begin{equation}
 -\de_{l,i}=\sum_{j=1}^{\ff} Q_{i,j}U_{j,l},
   \end{equation}
  which is (\ref{inv.34a}).
\qed

Lemma  \ref{lem-9.1mm}  gives the following useful inequality:

\begin{lemma} \label{lem-9.2} Let $X$, $Q$ and $U$ be as defined in Lemma \ref{lem-9.1mm}. 
Then,
  \be
 U_{i,i}\ge  \frac{1}{ |Q_{i,i}|}   ,\qquad \forall\, i\in\mN.\label{uniflb.9}
   \ee
\end{lemma}
\Proof
Since  $Q(i,i)<0$ it follows from (\ref{inv.34a}) that   
\begin{equation}
  1=|Q_{i,i}|U_{i,i}-\sum_{j\ne i} Q_{i,j}U_{j,i},
  \end{equation}
and since  $Q(i,j)\geq 0$ for $i\neq j$   we get (\ref{uniflb.9}). 
\qed

 The inequality in (\ref{uniflb.9}) can also be obtained  from the facts  that $1/  |Q_{i,i}|=1/q(i)$  is the expected amount of time     the process spends at $i$
during each visit to $i$, whereas $  U_{i,i}$ is the total expected amount of time spent at $i$ when the process starts at $i$.

\medskip	  We  say   that   a Markov chain $X$  is uniform when  it's $Q$ matrix has the property that $\|Q\|<\ff$.  When    a Markov chain is uniform we  can give  additional relationships between it's $Q$ matrix and its potential. 
   Since all the row sums of $Q$ are negative,  
     \begin{equation}
\sup_{j}  |Q_{j,j}|\le \|Q\|\leq 2\sup_{j}  |Q_{j,j}|.\label{unif.4}
 \end{equation}

   \begin{lemma} \label{lem-11.7}    Let $X$, $Q$ and $U$ be as defined  in Lemma \ref{lem-9.1mm} and assume that $X$ is a uniform Markov chain.  
   \begin{itemize}
   \item[(i)]If the row sums of $Q$ are bounded away from $0$ then 
 $\|U\| <\ff$.  
    \item[(ii)] If in addition  if $Q$ is  a $(2m+1)-$diagonal   matrix for some $m\geq 1$,
    \begin{equation}
U_{i,k}\leq C e^{-\la |i-k|},\qquad\forall i,k\in\mN,\label{u8.10m}
\end{equation}
for some constants $C,\la>0$. 
    \end{itemize}
 \end{lemma}

\medskip	 \Proof	$(i)$ If $ \|Q\|<\ff$ and the row sums of $Q$ are bounded away from $0$, then there exists $\bb>0$ such that , 
\begin{equation}
   \de=\|I+\bb Q\|<1.
  \end{equation} 
    It then follows from \cite[Section 5.3]{F} that   $\|e^{t\(I+\bb Q\)}\|\leq e^{\de t}$, or  equivalently, $\|e^{t Q}\|\leq e^{-\(1-\de\) t/\bb}$. Using   \cite[Section 5.3]{F} again, and the fact that  the transition semi-group,  $P_{t}=e^{t Q}$, we have
\begin{equation}
\Big\|\int_{0}^{\ff}  P_{t}\,dt\Big\|<\ff.
\end{equation}
Since $U=\int_{0}^{\ff}  P_{t}\,dt$, we have  $\|U\| <\ff$.

\medskip	\noindent $(ii)$
Let $\si=\inf \{t\,|\,X_{t}\neq X_{0}\}$,  the time of the first jump of $X$. Then for all $n\in\mN,$
\bea
P^{n}\(X_{\si}\in \mN\)&=& \sum_{i\neq n}P^{n}\(X_{\si}=i\)= \sum_{i\neq n} { Q_{n,i}  \over |Q_{n,n }| }.\label{9.12}
\eea
Note that since  the row sums   of $Q$ are bounded away from $0$ there exists a $\de>0$ such that,
\begin{equation}
|Q_{n,n }| -\sum_{i\neq n}   Q_{n,i}\geq \de, \label{u8.25} 
\end{equation}
uniformly in $n$. Furthermore, since  
$\sup_{n}|Q_{n,n }| \le \|Q\|$, we have, 
\begin{equation}
\sum_{i\neq n} { Q_{n,i}  \over |Q_{n,n }| }\le 1- \frac{\de}{   \|Q\|}.\label{u8.25a} 
\end{equation}
Therefore, by (\ref{9.12}),
\begin{equation}
P^{n}\(X_{\si}\in \mN\) \leq  1- \frac{\de}{   \|Q\|}, \qquad\forall n\in\mN.\label{u8.26a}
\end{equation}

We show immediately below that for all $i<k$, $i,k\in\mN$,    
\begin{equation}
   P^{i}\(T_{k}<\ff\)\leq e\(1- \frac{\de}{   \|Q\|}\)^{(k-i)/m}.\label{9.15}
  \end{equation}
Since,  
\begin{equation}
U_{i,k}=P^{i}\(T_{k}<\ff\)U_{k,k}\leq P^{i}\(T_{k}<\ff\)\|U\|, \label{u8.27}
\end{equation}
it follows that 
\begin{equation}
U_{i,k}\le e \|U\|  \(1- \frac{\de}{   \|Q\|}\)^{(k-i)/m}, \hspace{.2 in}i< k.\label{u8.27s}
\end{equation}
This gives (\ref{u8.10m}) with $C=e \|U\| $ and $\la=\de/\|Q\|$.

\medskip	We now obtain (\ref{9.15}).
Let $[(k-i)/m]=l$ and 
\begin{equation}
L_{j}=\{j,j+1, \ldots, j+m-1\}.\label{u8.20}
 \end{equation}
Since the Markov chain $X$ can move at most  $m$ units at each jump, 
\begin{equation}
\{X_{0}=i, T_{k}<\ff\}=\{X_{0}=i\} \cap_{j=1}^{l-1} \{S_{j}<\ff\} \cap \{T_{k}\circ S_{l-1}<\ff\},\label{u8.28}
\end{equation}
where $S_{1}=T_{L_{i+1}}$ and $S_{j}=T_{L_{i+(j-1)m+1}}\circ S_{j-1}$, $ j=2,\ldots, l-1$.
Then by the Markov property and (\ref{u8.26a})
\bea
P^{i}\(T_{k}<\ff\)&=&E^{i}\(  \cap_{j=1}^{l-1} \{S_{j}<\ff\}     E^{X_{S_{l-1}}}\(   T_{k}<\ff     \) \)\label{u8.29}\\
&\leq &  \(1- \frac{\de}{   \|Q\|}\)E^{i}\(  \cap_{j=1}^{l-1} \{S_{j}<\ff\}     \)\nn\\
&\leq& \(1- \frac{\de}{   \|Q\|}\)E^{i}\(  \cap_{j=1}^{l-2} \{S_{j}<\ff\}  E^{X_{S_{l-2}}}\(  T_{L_{i+(l-2)m+1}}<\ff  \)   \)\nn\\
&\leq& \(1- \frac{\de}{   \|Q\|}\)^{2}E^{i}\(  \cap_{j=1}^{l-2} \{S_{j}<\ff\}     \).\nn  
\eea
Continuing this procedure we get 
\begin{equation}
P^{i}\(T_{k}<\ff\) \leq  \(1- \frac{\de}{   \|Q\|}\)^{l}, 
  \end{equation}
which gives (\ref{9.15}).\qed

\medskip	\noindent{\bf  Proof of Theorem \ref{theo-1.10} }  
To show that the first condition in (\ref{1.39}) holds we use Lemma \ref{lem-9.2} and (\ref{unif.4}) to see that,
\begin{equation}
  U_{i,i}\ge  \frac{1}{ |Q_{i,i}|} \ge  \frac{1}{ \sup_{j}|Q_{j,j}|} \ge  \frac{1}{ \|Q\|} .
  \end{equation}
The second condition in (\ref{1.39}) is given in  Lemma \ref{lem-11.7}.

\medskip	  Now suppose that  $Q$ is  a $(2m+1)-$diagonal   matrix for some $m\geq 1$. It follows from (\ref{u8.10m})  and Lemma \ref{lem-8.2nn}    that $f\in c^{+}_{0}$  if and only if    $f= U  h$ for some $   h  \in c_{0}^{+}$. \qed 
 
 \begin{remark} {\rm 
When $X$ in Theorem \ref{theo-borelN} is a uniform Markov chain   with Q-matrix $Q$  and $f=Uh$ with $h\in \ell^{+}_{1}$, then it follows from the proof of the theorem that $\wt U$ is the restriction to  $\mN$ of the potential density of a uniform Markov chain on $\{0\}\cup\mN$   with Q-matrix
\begin{equation}
Q_{j,k},\qquad j,k\in\mN,\label{3.56nnx}
\end{equation} \vspace{-.2 in}
\begin{equation}
Q_{0,0}=1+\|h\|_{1}, \quad Q_{j,0}=-\sum_{k=1}^{\ff}Q_{j,k}, \quad j\in\mN, \quad \mbox{and }\quad  Q_{0,k}=-h_{k},\qquad k\in\mN.      \nn
\end{equation}
 It is clear that all the row sums of this Q-matrix are equal to 0, except for the first row sum which is equal to 1.
 
 }\end{remark}

  At the ends of Sections \ref{sec-gen}, \ref{sec-6}  and \ref{sec-7} we examine  the effects on the   covariances of certain Gaussian sequences that are also potentials of Markpov chains when we shift a parameter $\mathbf s$ by $\De>s_{1}$. We show that when the `shifted' covariance is itself a potential, all the elements of the $Q$ matrix of this new potential is are equal to  the elements  $Q$ matrix of the original potential, except for  $(1,1)$ coordinate with is a function of $\De.$ (See page \pageref{page40}).)
   The next lemma reverses and generalizes this proceedure. It   examines the effects on the potential  densities of Markov chains when we change any one term of their $Q$-matrices.
    We consider the matrix
  $E(k,l)=\{E(k,l)_{i,j};i,j\in\mN\}$,  with one non-zero element, where 
    \begin{equation}
 E(k,l)_{i,j}=\de_{(k,l)}(i,j).\label{3.86}
   \end{equation}

	 \begin{lemma}\label{lem-5.15} Let $Q$ be the Q-matrix 
  of a symmetric transient uniform Markov chain   on $ \mN$ with potential density $U$ satisfying,   \begin{equation}
  0< U_{j,k}Q_{j,k}<U_{k,k} Q_{j,k} ,\qquad\mbox{for some } j\ne k\in\mN,\label{5.161}
   \end{equation}
   and assume that for some   real number $b$ the matrix,
\begin{equation}
   Q+bE(k,l), \label{8.5}
   \end{equation} is the Q-matrix 
  of a transient   Markov chain $X$ on $ \mN$.
  
 Then if either,
  \begin{itemize} 
  \item[(i)] $Q$ is a   $(2n+1)-$diagonal matrix for some $n\ge 1$,
  
or    \item[(ii)]  $\sum_{j=1}^{\ff}U_{i,j} <\ff$ for each $i\ge 1$,
    \end{itemize}
     we have,   
     $b< 1/U_{l,k}$ and  the potential of $X$ is given by $W=\{W_{i,j};i,j\in\mN\}$  where,   
  \begin{equation}
 W_{i,j}= U_{i,j}+{b U_{i,k}  U_{ l,j }\over 1-bU_{k,l}} .\label{hars.1}
  \end{equation}
 
   \end{lemma}

 \medskip	 \Proof     In order for (\ref{8.5}) to   be  the Q-matrix 
  of a transient   Markov chain $X$ on $ \mN$,  we must have,
\begin{equation}
   Q_{k,k}+b+\sum_{j=1,j\ne k}^{\ff}Q_{k ,j}\le 0.
   \end{equation}
Therefore,
\begin{equation}
   Q_{k,k}U_{k,k}+bU_{k,k}+U_{k,k}\sum_{j=1,j\ne k}^{\ff}Q_{k,j}\le 0.\label{8.8}
   \end{equation}
It follows from (\ref{5.161}) that
\be
 \sum_{j=1,j\ne k}^{\ff}Q_{k,j}U_{j,k}<U_{k,k}\sum_{j=1,j\ne k}^{\ff}Q_{k,j}.
\ee
Consequently,
\begin{equation}
   Q_{k,k}U_{k,k}+bU_{k,k}+ \sum_{j=1,j\ne k}^{\ff}Q_{k,j}U_{j,k}< 0.
   \end{equation}
Therefore, by Lemma \ref{lem-9.1mm},
\begin{equation}
   -1+bU_{k,k} <0,
   \end{equation}
or,  equivalently,
\begin{equation}
       b<\frac{1}{U_{k,k}}.
   \end{equation}
   Since $U_{l,k}\le U_{k,k}$, this implies that $  b< 1/{U_{l,k}}$ for all $l\in\mN.$
We also note that $b\ge -Q_{j,k}$ when $j\ne k$.

\medskip We now obtain (\ref{hars.1}).	 Let  
  \be 
s:=s(l,k)=\frac{b}{ 1-bU_{l,k}}.\label{8.7}
  \ee
     By Lemma \ref{lem-invpot}
  it suffices to show that for each $i,n$  in   $\mN.$  
    \begin{equation}
 \sum_{j=1}^{\ff} \(   U_{i,j} +s U_{i,k}  U_{ \ell,j } \)\(Q_{j,n}+b E(k,l)_{j,n} \)=-I_{i,n}.\label{hars.2} 
  \end{equation}
We first note that,
 \begin{equation}
  \sum_{j=1}^{\ff} U_{i,j}E(k,l)_{j,n}= \left\{
    \begin{array}{cl}
    U_{i,k}&\qquad\mbox{when } n=l \\  
    0&\qquad \mbox{otherwise.} \\ 
  \end{array}\right.
     \end{equation}
We  write this as 
\begin{equation}
 \sum_{j=1}^{\ff} U_{i,j}E(k,l)_{j,n}=U_{i,k}I_{l,n}\label{4.16}.
   \end{equation}
It follows from this and  
  Lemma \ref{lem-9.1mm} that, 
  \begin{equation}
 \sum_{j=1}^{\ff}  U_{i,j}    \(Q_{j,n}+b E(k,l)_{j,n}\)=-I_{i,n}+bU_{i,k}I_{l,n}.\label{hars.3}
  \end{equation}
Using    
  Lemma \ref{lem-9.1mm} again we also see that,  
 \begin{equation}
 \sum_{j=1}^{\ff}    U_{i,k}  U_{ \ell,j }  Q_{j,n}= -U_{i,k}I_{l,n}, \label{hars.4}
 \end{equation}
 and by (\ref{4.16})    \begin{equation} 
  \sum_{j=1}^{\ff}    U_{i,k}  U_{ \ell,j } E(k,l)_{j,n}=U_{ \ell,k }U_{i,k}I_{l,n}. \label{hars.5}
 \end{equation}
 It follows from the last four equations that to get (\ref{hars.2}) we must have,  
\begin{equation}
    \(b-s+sbU_{l,k}\)U_{i,k}I_{l,n}=0,\qquad \forall i,n\in\mN,
   \end{equation}
which follows from (\ref{8.7}), since, $b-s+sbU_{l,k}=0$. \qed

\begin{remark} {\rm   Consider (\ref{hars.1}) with $k=l$. Then we can write
      \begin{equation}
   W_{i,j}=U_{i,j}+ c_{i}c_{j},\qquad \forall j,k\in\mN,\label{5.181q}
   \end{equation}
   where $c=\{c_{i}\}\in\mathbf Z$ is a sequence of real numbers.
     
       \medskip	If $k\ne l$, unless
\begin{equation}
   U_{i,k}U_{l,j}= U_{j,k}U_{l,i},
   \end{equation}
$W$ is not symmetric. Furthermore, unless
\begin{equation}
   U_{i,k}U_{l,j}=f_{j},\qquad\forall i\in\mN
   \end{equation}
 $W$ does not have the form of  (\ref{1.10}). It has the form
\begin{equation}
   W_{i,j}=U_{i,j}+ c_{i}d_{j},\qquad \forall j,k\in\mN.\label{5.181}
   \end{equation}
 where $d=\{d_{j}\}\in\mathbf Z$  and $d\ne c$.   In these cases $W$ is  a  new class of non-symmetric kernels for permanental processes.
 
  }\end{remark} 
    
 	   \begin{example} \label{ex-8.1}{\rm Consider the matrices $Q$ and $\wh W$ in (\ref{harr2.211x}) and (\ref{harr2.21dd}) and create the $Q$ matrix, \begin{equation}
   \wt Q=Q+bE(1,2), \label{8.5ss}
   \end{equation}
where
\begin{equation}
   b=\frac{b'}{1-r^{2}},\label{8.49}
   \end{equation}
   so that this first row of $\wt  Q$ is
   \begin{equation}
   \frac{1}{1-r^{2}}(-1, r+b',0,0,\ldots ),
   \end{equation}
   and all the other rows are unchanged. Since $r<1$ there are values of $b$ for which $\wt  Q$ is a $Q$ matrix.

 Using   (\ref{8.5}) and (\ref{hars.1}) we see that the potential corresponding to $\wt  Q$ is
   $\wt W=\{\wt W_{i,j};i,j\in\mN\}$  where,   
  \begin{equation}
 \wt W_{i,j}= \wh W_{i,j}+{b  \wh W_{i,1}   \wh W_{ 2,j }\over 1-b \wh W_{1,2}} = \wh W_{i,j}+{b  \wh W_{i,1}   \wh W_{ 2,j }\over 1-br}.\label{hars.1o}
  \end{equation}
  In particular 
 \be 
  \wt  W_{1,1}=\wh W_{1,1}+\frac{b\wh W_{1,1}\wh W_{2,1}}{1-br}=1+\frac{br}{1-br},
   \ee
   and for $j\ge 2$,
   \be  \wt W_{1,j}= \wh W_{1,j}+\frac{b\wh W_{1,1}\wh W_{2,j}}{1-br}=r^{j-1}+\frac{b r^{j-2}}{1-br} 
   \ee
   and
   \be \wt W_{j,1}=  \wh W_{j,1}+\frac{b\wh W_{j,1}\wh W_{2,1}}{1-br}=r^{j-1}+\frac{br^{j} }{1-br}\nn,
   \ee 
   which shows that $\wt W$ is not symmetric. 
      For $j,k\ge 2$,
  \begin{equation}
    \wt  W_{j,k}=r^{|j-k|}+\frac{br^{j+k-3} }{1-br},
     \end{equation}  
   so, for these values, $\wt W_{j,k}=\wt W_{k,j}$.
      
   Note that  for  $\wt  Q$ to be a $Q$ matrix we must have $r+b' \le 1$. Therefore, by (\ref{8.49}), we must have $ b\le 1/(1+r) $. Consequently we see that for $j\ge 2$,   $\wt W_{j,1} <  \wt W_{1,j}\le \wt W_{1,1}=\wt W_{2,2}$, although  obviously,
   \begin{equation}
      \lim_{j\to\ff}\wt  W_{j,j}=1.\label{8.55}
      \end{equation}

 \medskip	 Let $  X_{\al}=\{  X_{\al}(n),n\in\mN\}$ be an $\al$-permanental sequence with kernel $\wt W$. It follows from  Theorem \ref{theo-1.8gen} and (\ref{8.55}) that for all $\al>0$,
    \begin{equation}
       \limsup_{n\to\ff}\frac{X_{\al}(n)}{ \log n}=1,\qquad a.s.\label{8.54}
       \end{equation}

   }\end{example}


\begin{thebibliography} {10}


 
   \bibitem{BG}
R. Blumenthal,    and  R. Getoor,   
 {\em Markov {P}rocesses and {P}otential {T}heory}.
  New York: Academic Press, (1968). Page 95.

 
   \bibitem{DMM}
 C. Dellacherie, S. Martinez and  J. San Martin,
 {\em Inverse M-Matrices and Ultrametric Matrices},
Lecture Notes in Mathematics 2118, Springer,
 NY, (2014). Pages 85, 85.
 
% \bibitem {EK1}
%N. Eisenbaum  and  H. Kaspi,
%\newblock A charcterization of the infinitely divisible squared Gaussian processes, Annals of %Probability,\,
%{\em 34},  (2006),  728-742.
 



\bibitem {EK}
N. Eisenbaum  and  H. Kaspi,
\newblock On permanental processes, {\em Stochastic Processes and their Applications},
{\em 119},  (2009),  1401-1415. Pages 2, 27, 82.  

   \bibitem{F}
D. Freedman,
 {\em Markov chains},
Second Edition, Springer,
 NY, (1983). Pages 6, 99.

  
%     \bibitem{GKP}
%R. Graham, D. Knuth and O. Patashnik,
% {\em Concrete Mathematics: A Foundation for Computer Science},
%Second Edition, Addison-Wesley, 
% NY, (1994).



 \bibitem {Kovalb}
V, A. Koval,
\newblock The law of the iterated logarithm for   Gaussian sequences and its applications, {\em Theor. Probability and Math. Stat.},
{\em 54},  (1997),  69-75. Pages 8, 24.




\bibitem{book} M. B.  Marcus and J.~Rosen, {\em Markov Processes,
Gaussian Processes and Local Times}, Cambridge University Press, New
York,  (2006). Pages 3, 4.
 
\bibitem{MRsuf} M. B. Marcus and J.~Rosen, A sufficient condition for the continuity of permanental processes with
 applications to local times of Markov processes,     Annals of Probability,   (2013), Vol. 41, No. 2, 671-698. Page  15.
 

\bibitem{MRnec} M. B. Marcus and J.~Rosen,
Conditions for permanental processes to be unbounded, Annals of Probability,\, 45, (2017), 2059--2086.  Pages  15, 73, 73, 82, 84.



\bibitem{MRHD} M. B. Marcus and J.~Rosen, Permanental random variables, $M$-matrices and  $\al$-permanents,  High Dimensional Probability VII,  The Carg\`ese Volume, Progress in Probability, 71 (2016), 363--379, Birkhauser, Boston. Pages 82, 85.

\bibitem{MRejp} M. B. Marcus and J.~Rosen, Sample path properties of permanental processes,  EJP, (2018), No. 58, 1-49. Pages  2, 4, 5, 5, 14, 15, 27, 77, 84, 84, 95, 95, 96.


\bibitem{MRns} M. B. Marcus and J.~Rosen, Permanental sequences with kernels that are not equivalent to a symmetric matrix. Proceedings of the Eighth Conference on High Dimensional Probability, Oaxaca, Mexico. Page 2.


\bibitem{Norris} J. R.  Norris, {\em Markov Chains}, Cambridge University Press, New
York,  (1999). Pages 17, 17, 73.

 
 \bibitem{RY}
Daniel Revuz and Marc Yor.
\newblock {\em Continuous martingales and {B}rownian motion}, volume 293.
\newblock Springer-Verlag, Berlin, third edition, 1999. Pages 17, 17.

\bibitem{Revesz} P.  Revesz, {\em Random walk in random and non-random environments}, Cambridge University Press, New
York,  (1990). Page  90.



   \bibitem{TK}
H. Taylor and S. Karlin,
 {\em An Introduction to Stochastic Modeling},
Third Edition, Academic Press,
 NY, (1998). Pages 8, 16.

%\bibitem{VJ}  D.  Vere-Jones,
 %    {\em Alpha-permanents},
 %   New Zealand J. of Math.,  (1997), 26, 125--149.



\end{thebibliography}
  \end{document}